\documentclass[12pt]{amsart}

\usepackage[T1]{fontenc}
\usepackage{amssymb}
\usepackage{bm}
\usepackage{amsmath}
\usepackage{stmaryrd}
 \usepackage{mathpple}

\usepackage{amsmath,amsfonts, amscd,amsthm, graphicx}
 \usepackage{amsrefs}
\usepackage[all]{xy}
\usepackage{color}
\usepackage{quiver}
\usepackage{tikz-cd}
\usepackage{tikz}
\usetikzlibrary{arrows,decorations.markings}

\usepackage{hyperref}

\def\XXint#1#2#3{{\setbox0=\hbox{$#1{#2#3}{\int}$}
\vcenter{\hbox{$#2#3$}}\kern-.5\wd0}}

\numberwithin{equation}{section}

\newcommand{\mbf}{\mathbf}
\newcommand{\mbb}{\mathbb}
\newcommand{\mf}{\mathfrak}
\newcommand{\mc}{\mathcal}

\theoremstyle{plain}
\newtheorem{thm}{Theorem}[section]

\newtheorem{thm-defn}{Theorem/Definition}[section]
\newtheorem{lem}[thm]{Lemma}
\newtheorem{lem-defn}[thm]{Lemma/Definition}
\newtheorem{prop}[thm]{Proposition}
\newtheorem{cor}[thm]{Corollary}

\newtheorem{prop-defn}[thm]{Proposition-Definition}

\newtheorem{thm-alg}[thm]{Theorem/Algorithm}

\allowdisplaybreaks[4] 

\begin{document}

\tikzset{
 ->-/.style={decoration={markings, mark=at position 0.2 with {\arrow{stealth}}},
    postaction={decorate}},
}

 \title{Quantum Difference Equations for Affine Type $A$ Quiver Varieties I: General Construction}
 \author{Tianqing Zhu}
 \address{Department of Mathematics, Columbia University and Yau Mathematical Sciences Center, Tsinghua University}
\email{tz2611@columbia.edu}

 \date{}

 \maketitle
\begin{abstract}
Following the philosophy of \cite{OS22}, we construct quantum difference equations for affine type $A$ quiver varieties in terms of the quantum toroidal algebra $U_{q,t}(\hat{\hat{\mf{sl}}}_{n})$. As part of the construction, we define the wall set of each quiver variety through the action of the universal $R$-matrix and show that, on each interval in $H^2(X,\mbb{Q})$, this wall set is essentially equivalent to the one determined by the $K$-theoretic stable envelope. We illustrate the construction by deriving explicit quantum difference operators for the instanton moduli space $M(n,r)$ and the equivariant Hilbert scheme $\text{Hilb}_{n}([\mbb{C}^2/\mbb{Z}_{r}])$.
\end{abstract}

\tableofcontents

\section{\textbf{Introduction}}

\subsection{Quantum difference equation and $K$-theoretic enumerative geometry}
It is an important problem to establish a connection between the enumerative geometry and the geometric representation theory. For the case of the Nakajima quiver varieties, Maulik and Okounkov showed \cite{MO12} that the quantum differential equation in the equivariant cohomology of the Nakajima quiver variety is equivalent to the trigonometric Casimir connection of the corresponding quiver type of the geometric Yangian. 

The solution to the quantum differential equation is connected to the genus zero Gromov-Witten theory curve counting and quantum cohomology for the Nakajima quiver varieties. It has deep connection to many cohomological curve counting theory. For example, In the case of $\text{Hilb}_n(\mbb{C}^2)$ the Hilbert scheme of points on $\mbb{C}^2$, its equivariant quantum cohomology can be described in terms of the genus zero Gromov-Witten/Donaldson-Thomas theory on $\mbb{C}^2\times\mbb{P}^1$ as in \cite{BP08}\cite{OP10}.

In the equivariant $K$-theory, one
could also carry out the similar story by replacing the quantum differential equation
with the quantum difference equation. The enumerative invariants here are replaced by
the capped operators, capping operators and vertex functions. In \cite{O15}, Okounkov shows that the capping operator $\mbf{J}(z)$ is the regular solution to the $q$-difference equations:
\begin{align*}
\mbf{J}(zq^{\mc{L}})\mc{L}=\mbf{M}_{\mc{L}}(z)\mbf{J}(z)
\end{align*}
with $\mc{L}$ the line bundle over $M(\mbf{v},\mbf{w})$. Moreover, 
Okounkov and Smirnov \cite{OS22} proved that the quantum difference equation for the Nakajima quiver variety is equivalent to a difference analog of the trigonometric Casimir connection of the corresponding quiver type of the geometric quantum loop group. In \cite{AO21}, it was also shown that the monodromy for the capping operator is described by the elliptic stable envelopes.

The foundational work on quantum difference equations was developed by Etingof
and Varchenko \cite{EV02}, who constructed the dynamical Weyl group, subsequently employing it to construct difference operators for the intertwiners for the $U_{q}(\hat{\mf{g}})$-modules of
finite ADE type. These difference equations are the difference analog of the trigonometric Casimir connections, and will be referred to as quantum difference equations.

One important question in the construction of the quantum difference equation is to get a systematic way to derive the formula for the quantum difference operator. Though it is known that we can use the geometric quantum loop group to construct the quantum difference operator. It is not direct to know the explicit formula for the elements of the geometric quantum loop groups and its corresponding wall subalgebras.

On the side of the quantum differential equation, since the quantum difference equation can be degenerated to the quantum differential equation as described in \cite{BM15}\cite{Z24}, the quantum difference equation will give more intrinsic understanding on the structure for the quantum differential equation. For example, the monodromy representation for the quantum differential equation turns out to be connected to the study of the structure for the quantum difference equation as shown in \cite{S21}\cite{Z24}. 

\subsection{Main result of the paper}
In this paper, we construct a holonomic system of difference equations on $\text{Pic}(M(\mbf{v},\mbf{w}))\otimes\mbb{C}^*$ with values in the endomorphism of the equivariant $K$-theory of affine type $A$ quiver varieties by using the quantum toroidal algebra $U_{q,t}(\hat{\hat{\mf{sl}}}_{n})$. This system is an algebraic generalisation of the quantum difference equations constructed in \cite{OS22} for Nakajima quiver varieties. Our construction relies on the slope factorisation of the quantum toroidal algebra established in \cite{N15}. 

Using this factorisation, we systematically define a monodromy operator $\mbf{B}_{\mbf{m}}(z)$ for each slope subalgebra $\mc{B}_{\mbf{m}}$ and then construct the algebraic quantum difference operator $\mbf{B}_{\mc{L}}^s(z)$. 
The family of the difference operators $\{\mbf{B}_{\mc{L}}^s(z)\}_{\mc{L}\in\text{Pic}(M(\mbf{v},\mbf{w}))}$ forms a holonomic system over $\text{Pic}(M(\mbf{v},\mbf{w}))\otimes\mbb{C}^{*}$ for the slope point $s\in\mbb{R}^n$. The following two theorems are the main results of the paper.

\begin{thm}\label{main-theorem}(Theorem \ref{holonomic-theorem})
Let $M(\mbf{v},\mbf{w})$ be an affine type $A$ quiver variety. For line bundles $\mc{L}\in\text{Pic}(M(\mbf{v},\mbf{w}))$, there is a holonomic system of operators $\mbf{B}_{\mc{L}}^s(z)\in\text{End}(K_{T}(M(\mbf{v},\mbf{w}))(z^{\mc{L}})_{\mc{L}\in\text{Pic}(M(\mbf{v},\mbf{w}))})$; that is,
\begin{align*}
T_{p,\mc{L}}^{-1}\mbf{B}_{\mc{L}}^s(z)T_{p,\mc{L}'}^{-1}\mbf{B}_{\mc{L}'}^s(z)=T_{p,\mc{L}'}^{-1}\mbf{B}_{\mc{L}'}^s(z)T_{p,\mc{L}}^{-1}\mbf{B}_{\mc{L}}^s(z).
\end{align*}

\end{thm}

In addition, we give a universal formula for the monodromy operators $\mbf{B}_{\mbf{m}}(z)$ of affine type $A$ quiver varieties under the expansion around $z=0$.

\begin{thm}\label{universal-simplified version}(Theorem \ref{universal-formula-for-difference-operator})
On the representation space $\text{End}(K_T(M(\mbf{v},\mbf{w})))$, the monodromy operator $\mbf{B}_{\mbf{m}}(z)$ admits the expansion around $z=0$:
\begin{equation*}
\begin{aligned}
&\mbf{B}_{\mbf{m}}(z)\\
=&\prod_{h=1}^{g}:(\textbf{Heisenberg algebra part})\prod_{k=0}^{\substack{\rightarrow\\\infty}}\\
&\prod_{\substack{\gamma\in\Delta(A)\\m\geq0}}^{\leftarrow}(\exp_{q^{2}}(-(q-q^{-1})z^{k(-\mbf{v}_{\gamma}+(m+1)\bm{\delta}_h)}p^{-k\mbf{m}\cdot(-\mbf{v}_{\gamma}+(m+1)\bm{\delta}_h)}q^{\frac{k}{2}(\mbf{v}^TC\mbf{v}-\mbf{w}^T\mbf{w})\bm{\theta}\cdot(-\mbf{v}_{\gamma}+(m+1)\bm{\delta}_h))}(\\
&q^{-2k(-\mbf{v}_{\gamma}+(m+1)\bm{\delta})^TC(-\mbf{v}_{\gamma}+(m+1)\bm{\delta})}\varphi^{-(k+1)(-\mbf{v}_{\gamma}+(m+1)\bm{\delta})}e_{(\delta-\gamma)+m\delta}f_{(\delta-\gamma)+m\delta}'\varphi^{-(k+1)(-\mbf{v}_{\gamma}+(m+1)\bm{\delta})})\\
&\times\exp(-(q-q^{-1})\sum_{m\in\mbb{Z}_{+}}\sum_{i,j=1}^{l_h}z^{km\bm{\delta}_h}p^{-km\mbf{m}\cdot\bm{\delta}_{h}}u_{m,ij}q^{-2km^2\bm{\delta}^TC\bm{\delta}}\varphi^{-(k+1)m\bm{\delta}}e_{m\delta,\alpha_i}f_{m\delta,\alpha_i}'\varphi^{-(k+1)m\bm{\delta}})\\
&\times\prod_{\substack{\gamma\in\Delta(A)\\m\geq0}}^{\rightarrow}\exp_{q^{2}}(-(q-q^{-1})z^{k(\mbf{v}_{\gamma}+m\bm{\delta}_h)}p^{-k\mbf{m}\cdot(\mbf{v}_{\gamma}+m\bm{\delta}_h)}q^{\frac{k}{2}(\mbf{v}^TC\mbf{v}-\mbf{w}^T\mbf{w})\mbf{\theta}\cdot(\mbf{v}_{\gamma}+m\bm{\delta}_h))}\\
&q^{-2k(\mbf{v}_{\gamma}+m\bm{\delta})^TC(\mbf{v}_{\gamma}+m\bm{\delta})}(\varphi^{-(k+1)(\mbf{v}_{\gamma}+m\bm{\delta})}e_{\gamma+m\delta}f_{\gamma+m\delta}'\varphi^{-(k+1)(\mbf{v}_{\gamma}+m\bm{\delta})})):.
\end{aligned}
\end{equation*}
Here $::$ denotes normal ordering.
\end{thm}

As illustrative examples, we compute the algebraic quantum difference operators for the equivariant Hilbert scheme $\text{Hilb}_{n}([\mbb{C}^2/\mbb{Z}_{r}])$ associated with the $A_{r-1}$ surface and for the instanton moduli space $M(n,r)$ on $\mbb{P}^2$. In the latter case, the resulting quantum difference equation agrees exactly with the equation constructed in \cite{OS22}.

In the subsequent paper \cite{Z24}, we use this construction to solve the quantum difference equations for affine type $A$ quiver varieties and study their degeneration to quantum differential equations. We will show that the limiting behaviour of the corresponding connection matrices admits a particularly simple description in terms of the monodromy operators $\mbf{B}_{\mbf{m}}(z)$.

The relationship between our construction and the approach of \cite{OS22} may be described as follows. Okounkov and Smirnov use the $K$-theoretic stable envelope to construct the geometric quantum loop group $U_{q}^{MO}(\hat{\mf{g}}_{Q})$. Moreover, the $K$-theoretic  stable envelopes of different generic slopes $\mbf{m}\in\mbb{R}^n$ then determine the wall subalgebras $U_{q}^{MO}(\mf{g}_{\mbf{m}})$. The comparison between $U_{q}^{MO}(\hat{\mf{g}}_{Q})$ and $U_{q}(\hat{\mf{g}}_{Q})$, and between their respective wall subalgebras, is governed by the following theorem in \cite{Z25}.

\begin{thm}\label{conjecture-stable-envelope}
There are isomorphisms of algebras:
\begin{align*}
U_{q}^{MO}(\hat{\mf{g}}_{Q})\cong U_{q}(\hat{\mf{g}}_{Q}), \qquad U_{q}^{MO}(\mf{g}_{\mbf{m}})\cong\mc{B}_{\mbf{m}}.
\end{align*}
\end{thm}

In this paper, we also prove that the geometric quantum difference equation for the affine type $A$ quiver variety is the algebraic quantum difference equation of the slope $s$ infinitesimally close to $0\in\mbb{R}^n$:

\begin{thm}[See Theorem \ref{geometric-qde=algebraic-qde:label}]
The algebraic quantum difference operator for the affine type $A$ Nakajima quiver variety of the slope $s$ infinitesimally close to $0\in\mbb{R}^n$ coincides with the geometric quantum difference operator up to a constant operator.
\end{thm}

Moreover, we expect this construction to reveal concrete structures in the quantum difference and quantum differential equations of quiver varieties, including their monodromy representations. The degeneration from quantum difference equations to quantum differential equations has been established in Dynkin ADE type \cite{BM15}, and we expect an analogous relationship in the quantum toroidal setting \cite{BT19}. A detailed analysis of this degeneration is given in \cite{Z24}.

The construction also suggests applications to enumerative geometry. In particular, under Theorem \ref{geometric-qde=algebraic-qde:label}, the regular solution to the quantum difference equation corresponds to the capping operator in the equivariant $K$-theory of the Nakajima quiver variety.

\subsection{Connection with qKZ}
We also prove that the quantum difference operators commute with the two-point qKZ operator.
\begin{thm}(See Theorem \ref{qKZ-commute})
For arbitrary line bundles $\mc{L},\mc{L}'\in\text{Pic}(M(\mbf{v},\mbf{w}))$ and a slope $s\in\mbb{R}^n$, the qKZ operators commute with the $q$-difference operators:
\begin{align*}
\Delta_{s}(T_{p,\mc{L}}^{-1}\mbf{B}^s_{\mc{L}}(z))T_{u}^{-1}z_{(1)}\mc{R}^s(u)=T_{u}^{-1}z_{(1)}\mc{R}^s(u)\Delta_{s}(T_{p,\mc{L}}^{-1}\mbf{B}^s_{\mc{L}}(z)).
\end{align*}
Here $\mc{R}^s(u)$ is the image of the universal $R$-matrix $R^s$ in the representation $\bigoplus_{\mbf{v}_1,\mbf{v}_2\in\mbb{N}^I}K_{T}(M(\mbf{v}_1,\mbf{w}_1))\otimes K_{T}(M(\mbf{v}_2,\mbf{w}_2))$.
\end{thm}

More generally, one expects the quantum difference equations to commute with arbitrary $n$-point qKZ equations for the quantum toroidal algebra $U_{q,t}(\hat{\hat{\mf{sl}}}_n)$ in the sense of \cite{FR92}. Analogous commutativity results are known for the trigonometric Casimir connection and the difference Casimir connection in ADE type \cite{EV02,L11}. Such commutativity provides another perspective on the connection matrices and solutions of the quantum difference equations.

From the viewpoint of enumerative geometry, the qKZ equation in the setting of \cite{OS22} corresponds to the difference equation in the equivariant variables satisfied by the capping operators in \ref{capping-difference}. The $R$-matrix $z_{(1)}\mc{R}^{(s)}(u)$ corresponds to the shift operator $S(u,z)$. We therefore expect these constructions to clarify the relationship between $K$-theoretic counts and intertwiners of modules over quantum algebras.

\subsection{Further directions}
This paper is a first step toward understanding quantum difference equations for Nakajima quiver varieties, and more generally for quiver varieties arising in $K$-theoretic enumerative geometry, through the theory of shuffle algebras. Particularly rich structures emerge when $K$-theoretic quasimap theory and shuffle algebras are studied in the framework of $K$-theoretic Hall algebras.

In \cite{Z24}, we analyse the quantum difference equations for affine type $A$ quiver varieties. A central object is the connection matrix relating solutions in different asymptotic regions; this matrix is closely related to an ordered product of the monodromy operators $\mbf{B}_{\mbf{m}}(z)$. By choosing a suitable generic path representing the quantum difference operator $\mbf{B}_{\mc{L}}^s(z)$, one may arrange for the individual monodromy operators to take their simplest forms, namely those of $U_{q}(\mf{sl}_2)$ type or $U_{q}(\hat{\mf{gl}}_1)$ type. This reflects a ``root decomposition'' of the slope subalgebra $\mc{B}_{\mbf{m}}$, analogous to the root decomposition of a classical Lie algebra. These results imply that the degeneration of the quantum difference equation yields the Dubrovin connection for the quantum cohomology of the affine type $A$ quiver variety, and that its monodromy representation can be recovered from the operators $\mbf{B}_{\mbf{m}}(z)$.

The construction is expected to extend to more general quivers and weighted modules. In \cite{Z26}, we generalise it to arbitrary types of quiver varieties and modules. Algebraically, this extension may be viewed as a generalisation of qKZ equations to other classes of modules and quantum groups, with potential applications in both representation theory and enumerative geometry.

\subsection{Organization of the paper}
The paper is organised as follows. In Section $2$, we review the quantum toroidal and quantum affine algebras, the shuffle-algebra realisation of the quantum toroidal algebra, and its slope subalgebras. We also describe how the coproduct and universal $R$-matrix of the quantum toroidal algebra relate to those of the slope subalgebras; these results are used in the construction of the quantum difference operators. 

Section $3$ contains the main construction of the quantum difference operators and the proofs of Theorems \ref{holonomic-theorem} and \ref{universal-simplified version}. In Section $4$, we treat the instanton moduli space and the equivariant Hilbert scheme $\text{Hilb}_{n}([\mbb{C}^2/\mbb{Z}_r])$ as examples.

\subsection*{Acknowledgments.}
The author is grateful to Nicolai Reshetikhin and Si Li for introducing him to the rich interplay between mathematical physics and representation theory, for their sustained support throughout this project, and for their insightful suggestions on the structure of the paper. The author also thanks Andrei Okounkov, Andrei Negu\c{t}, and Andrey Smirnov for helpful discussions on quantum groups and stable envelopes and for carefully reading a draft of the paper. Further thanks are due to Peng Shan and Changjian Su for valuable discussions and excellent seminars on geometric representation theory. The author was supported by the BMSTC--ACZSP International Collaboration Grant (Grant No. Z221100002722017).

\section{\textbf{Quantum toroidal algebra and geometric actions}}

In this section, we review the basic notions concerning quantum toroidal algebras and quantum affine groups that will be used throughout the paper. For details, see \cite{N15}\cite{N20}.

\subsection{Quantum groups $U_{q}(\hat{\mf{sl}}_{n})$}
The quantum group $U_{q}(\hat{\mf{sl}}_{n})$ is a $\mbb{Q}(q)$-Hopf algebra:
\begin{align*}
U_{q}(\hat{\mf{sl}}_{n})=\mbb{Q}(q)\langle x_{i}^{\pm1},\psi^{\pm1}_{s},c^{\pm1}\rangle^{i\in\mbb{Z}/n\mbb{Z}}_{s\in\{1,\cdots,n\}}
\end{align*}
modulo the fact that $c$ is central, as well as the following relations:
\begin{align*}
&\psi_{s}\psi_{s'}=\psi_{s'}\psi_{s}\\
&\psi_{s}x_{i}^{\pm}=q^{\pm(\delta^{i+1}_{s}-\delta^{i}_{s})}x_{i}^{\pm}\psi_{s}\\
&[x_{i}^{\pm},x_{j}^{\pm}]=0,\qquad\text{if }j\notin\{i-1,i+1\}\\
&[x_{i}^{\pm},[x_{i}^{\pm},x_{j}^{\pm}]_{q}]_{q^{-1}}=0\qquad\text{if }j\in\{i-1,i+1\}\\
&[x_{i}^{+},x_{j}^{-}]=\frac{\delta^{j}_{i}}{q-q^{-1}}(\frac{\psi_{i+1}}{\psi_{i}}-\frac{\psi_{i}}{\psi_{i+1}})
\end{align*}
for all $i,j\in\mbb{Z}/n\mbb{Z}$ and $s,s'\in\{1,\cdots,n\}$. Here $[-,-]_{q}$ is the $q$-bracket:
\begin{align*}
[a,b]_{q}=ab-qba.
\end{align*}

We further require the following relation:
\begin{align*}
\psi_{s+n}=c\psi_{s},\forall s\in\mbb{Z}.
\end{align*}

The counit is given by $\epsilon(x_{i}^{\pm})=0$ and $\epsilon(\psi_{s})=1$, while the coproduct is given by
\begin{align*}
&\Delta(c)=c\otimes c,\Delta(\psi_{s})=\psi_{s}\otimes\psi_{s}\\
&\Delta(x_{i}^{\pm})=\frac{\psi_{i+1}}{\psi_{i}}\otimes x_{i}^{+}+x_{i}^{+}\otimes1\\
&\Delta(x_{i}^{-})=1\otimes x_{i}^{-}+x_{i}^{-}\otimes\frac{\psi_{i}}{\psi_{i+1}}.
\end{align*}

We also introduce the two half-subalgebras of $U_{q}(\hat{\mf{sl}}_{n})$:
\begin{align*}
&U_{q}^{\geq}(\hat{\mf{sl}}_{n})=\mbb{Q}(q)\langle x_{i}^{+},\psi_{s}^{\pm1},c^{\pm1}\rangle^{i\in\mbb{Z}/n\mbb{Z}}_{s\in\{1,\cdots,n\}}\subset U_{q}(\hat{\mf{sl}}_{n})\\
&U_{q}^{\leq}(\hat{\mf{sl}}_{n})=\mbb{Q}
(q)\langle x_{i}^{-},\psi_{s}^{\pm1},c^{\pm1}\rangle^{i\in\mbb{Z}/n\mbb{Z}}_{s\in\{1,\cdots,n\}}\subset U_{q}(\hat{\mf{sl}}_{n}).
\end{align*}

Both are bialgebras, and they are equipped with the following bialgebra pairing, called the Drinfeld pairing:
\begin{align}\label{quantum-sln-bialgebra-pairing}
\langle-,-\rangle:U_{q}^{\geq}(\hat{\mf{sl}}_{n})\otimes U_{q}^{\leq}(\hat{\mf{sl}}_{n})\rightarrow\mbb{Q}(q)
\end{align}
such that:
\begin{align*}
\langle x_{i}^{+},x_{j}^{-}\rangle=\frac{\delta^{i}_{j}}{q^{-1}-q},\qquad\langle\psi_{s},\psi_{s'}\rangle=q^{-\delta^{s}_{s'}}
\end{align*}

\subsection{The Heisenberg algebra}
We consider the $q$-deformed Heisenberg algebra:
\begin{align*}
U_{q}(\hat{\mf{gl}}_{1})=\mbb{Q}(q)\langle p_{\pm k},c^{\pm1}\rangle_{k\in\mbb{Z}-\{0\}}
\end{align*}
where $c$ is central and the generators $p_{\pm k}$ all commute, except for:
\begin{align*}
[p_{k},p_{-k}]=k\cdot\frac{c^{k}-c^{-k}}{q^{k}-q^{-k}}
\end{align*}

The counit $\epsilon$ is given by $\epsilon(p_k)=0$ and $\epsilon(c)=1$. The coproduct is determined by $\Delta(c)=c\otimes c$ and
\begin{align*}
&\Delta(p_k)=c^k\otimes p_k+p_k\otimes1\\
&\Delta(p_{-k})=1\otimes p_{-k}+p_{-k}\otimes c^{-k}
\end{align*}

As in the previous subsection, the $q$-deformed Heisenberg algebra $U_{q}(\hat{\mf{gl}}_{1})$ can be realised as the Drinfeld double of its two halves:
\begin{align*}
&U_{q}^{\geq}(\hat{\mf{gl}}_1)=\mbb{Q}(q)\langle p_{k},c^{\pm1}\rangle_{k\in\mbb{N}}\subset U_{q}(\hat{\mf{gl}}_1)\\
&U_{q}^{\leq}(\hat{\mf{gl}}_1)=\mbb{Q}(q)\langle p_{-k},c^{\pm1}\rangle_{k\in\mbb{N}}\subset U_{q}(\hat{\mf{gl}}_1)
\end{align*}
with respect to the bialgebra pairing:
\begin{equation}\label{Heisenberg-algebra-pairing}
\begin{tikzcd}
\langle-,-\rangle: U_{q}^{\geq}(\hat{\mf{gl}}_1)\otimes U_{q}^{\leq}(\hat{\mf{gl}}_1)\arrow[r]&\mbb{Q}(q)
\end{tikzcd}
\end{equation}
such that
\begin{align*}
\langle p_{k},p_{-k}\rangle=\frac{k}{q^{-k}-q^k}
\end{align*}

\subsection{Quantum affine group $U_{q}(\hat{\mf{gl}}_{n})$}

The quantum affine group $U_{q}(\hat{\mf{gl}}_{n})$ can be defined by the RTT formalism (see Section 5.1 of \cite{N15} for a detailed definition):
\begin{align*}
U_{q}(\hat{\mf{gl}}_{n}):=\mbb{Q}(q)\langle e_{\pm[i;j)},\psi_{s}^{\pm1},c^{\pm1}\rangle^{s\in\{1,\cdots,n\}}_{(i,j)\in\mbb{Z}^2}/(\text{relations})
\end{align*}
where $[i,j)$ is defined in \eqref{cyclic-vectors} such that $e_{\pm[i+n,j+n)}=e_{\pm[i,j)}$, and it satisfies the following relations:
\begin{align*}
\psi_{k}\cdot e_{\pm[i;j)}=q^{\pm(\delta^{i}_k-\delta^j_k)}e_{\pm[i;j)}\psi_k,\qquad\forall\text{arcs }[i;j)\text{ and }k\in\mbb{Z}.
\end{align*}

For arbitrary $a\leq c$ and $b\leq d$, we also have
\begin{align}
\frac{e_{\pm[a,c)}e_{\pm[b,d)}}{q^{\delta^b_a-\delta^d_b+\delta^d_a}}-\frac{e_{\pm[b,d)}e_{\pm[a,c)}}{q^{\delta^b_c-\delta^d_b+\delta^d_c}}=(q-q^{-1})[\sum_{a\leq x<c}^{x\equiv d}e_{\pm[b,c+d-x)}e_{\pm[a,x)}-\sum_{a<x\leq c}^{x\equiv b}e_{\pm[x,c)}e_{\pm[a+b-x,d)}]
\end{align}

\begin{align}\label{gln-relation}
[e_{[a,c)},e_{-[b,d)}]=(q-q^{-1})[\sum^{x\equiv b}_{a\leq x<c}\frac{e_{[c+b-x,d)}e_{[a,x)}}{q^{\delta^b_c+\delta^a_c-\delta^a_b}}\frac{\psi_x}{\psi_c}-\sum_{a<x\leq c}^{x\equiv d}\frac{e_{[x,c)}e_{-[b,a+d-x)}}{q^{-\delta^b_a+\delta^d_b-\delta^d_a}}\frac{\psi_x}{\psi_a}].
\end{align}

Their standard coproduct is given by
\begin{align*}
&\Delta(e_{[i;j)})=\sum_{s=i}^{j}e_{[s;j)}\frac{\psi_{s}}{\psi_{i}}\otimes e_{[i;s)}\\
&\Delta(e_{-[i;j)})=\sum_{s=i}^{j}e_{-[i;s)}\otimes e_{-[s;j)}\frac{\psi_{i}}{\psi_{s}}.
\end{align*}

There is an isomorphism of Hopf algebras:
\begin{align}\label{quantum gln decomposition}
U_{q}(\hat{\mf{gl}}_{n})\cong U_{q}(\hat{\mf{sl}}_{n})\otimes U_{q}(\hat{\mf{gl}}_{1})/(c\otimes1-1\otimes c)
\end{align}
which is given by:
\begin{align*}
&e_{[i;i+1)}=x_{i}^{+}(q-q^{-1})\\
&e_{-[i;i+1)}=x_{i}^{-}(q^{-2}-1).
\end{align*}

The isomorphism identifies the quantum group $U_{q}(\hat{\mf{gl}}_{n})$ with $U_{q}(\hat{\mf{sl}}_{n})\otimes U_{q}(\hat{\mf{gl}}_1)/(c\otimes 1-1\otimes c)$. Here the Heisenberg generators $\{p_{\pm k}\}_{k\in\mbb{Z}}\in U_{q}(\hat{\mf{gl}}_1)$ in $U_{q}(\hat{\mf{gl}}_{n})$ satisfy the following relation:
\begin{align*}
[p_{k},p_{l}]=k\delta^0_{k+l}\frac{(c^k-c^{-k})(q^{nk}-q^{-nk})}{(q^k-q^{-k})^2}.
\end{align*}

In this paper, we use a slightly smaller subalgebra of $U_{q}(\hat{\mf{gl}}_{n})$, obtained by restricting the set of Cartan elements:
\begin{align*}
U_{q}(\hat{\mf{gl}}_{n}):=\mbb{Q}(q)\langle e_{\pm[i;j)},\varphi_{s}^{\pm1},c^{\pm1}\rangle^{s\in\{1,\cdots,n\}}_{(i,j)\in\mbb{Z}^2/(n,n)\mbb{Z}},\qquad\varphi_s:=\frac{\psi_{s+1}}{\psi_s}.
\end{align*}

We continue to denote this subalgebra by $U_{q}(\hat{\mf{gl}}_{n})$. The isomorphism $U_{q}(\hat{\mf{gl}}_{n})\cong U_{q}(\hat{\mf{sl}}_{n})\otimes U_{q}(\hat{\mf{gl}}_{1})/(c\otimes1-1\otimes c)$ remains valid, with the Cartan elements $\psi_s$ of $U_{q}(\hat{\mf{sl}}_{n})$ replaced by $\varphi_s$.

\subsection{Quantum toroidal algebra $U_{q,t}(\hat{\hat{\mf{sl}}}_{n})$}

The quantum toroidal algebra $U_{q,t}(\hat{\hat{\mf{sl}}}_{n})$ is defined by:
\begin{align*}
U_{q,t}(\hat{\hat{\mf{sl}}}_{n})=\mbb{C}\langle\{e_{i,d}^{\pm}\}_{1\leq i\leq n}^{d\in\mbb{Z}},\{\varphi_{i,d}^{\pm}\}_{1\leq i\leq n}^{d\in\mbb{N}_{0}}\rangle/(\sim)
\end{align*}
with $\varphi_{i,d}^{\pm}$ commuting among themselves and the relation $(\sim)$ is given by:
\begin{equation}
\begin{aligned}
&e_i^{ \pm}(z) \varphi_j^{ \pm^{\prime}}(w) \cdot \zeta\left(\frac{w^{ \pm 1}}{z^{ \pm 1}}\right)=\varphi_j^{ \pm^{\prime}}(w) e_i^{ \pm}(z) \cdot \zeta\left(\frac{z^{ \pm 1}}{w^{ \pm 1}}\right) \\
&e_i^{ \pm}(z) e_j^{ \pm}(w) \cdot \zeta\left(\frac{w^{ \pm 1}}{z^{ \pm 1}}\right)=e_j^{ \pm}(w) e_i^{ \pm}(z) \cdot \zeta\left(\frac{z^{ \pm 1}}{w^{ \pm 1}}\right) \\
&{\left[e_i^{+}(z), e_j^{-}(w)\right]=\delta_i^j \delta\left(\frac{z}{w}\right) \cdot \frac{\varphi_i^{+}(z)-\varphi_i^{-}(w)}{q-q^{-1}}}\\
& e_i^{ \pm}\left(z_1\right) e_i^{ \pm}\left(z_2\right) e_{i \pm^{\prime} 1}^{ \pm}(w)+\left(q+q^{-1}\right) e_i^{ \pm}\left(z_1\right) e_{i \pm^{\prime} 1}^{ \pm}(w) e_i^{ \pm}\left(z_2\right)+e_{i \pm^{\prime} 1}^{ \pm}(w) e_i^{ \pm}\left(z_1\right) e_i^{ \pm}\left(z_2\right)+\\
& +e_i^{ \pm}\left(z_2\right) e_i^{ \pm}\left(z_1\right) e_{i \pm^{\prime} 1}^{ \pm}(w)+\left(q+q^{-1}\right) e_i^{ \pm}\left(z_2\right) e_{i \pm^{\prime} 1}^{ \pm}(w) e_i^{ \pm}\left(z_1\right)+e_{i \pm^{\prime} 1}^{ \pm}(w) e_i^{ \pm}\left(z_2\right) e_i^{ \pm}\left(z_1\right)=0.
\end{aligned}
\end{equation}
where the generating series are defined by
\begin{equation*}
e_i^{ \pm}(z)=\sum_{d \in \mathbb{Z}} e_{i, d}^{ \pm} z^{-d} \quad \varphi_i^{ \pm}(z)=\sum_{d=0}^{\infty} \varphi_{i, d}^{ \pm} z^{\mp d}.
\end{equation*}

Here $i,j\in\{1,\cdots,n\}$, and $z$ and $w$ are variables of colours $i$ and $j$, respectively. The rational function $\zeta(x_i/x_j)$ is defined by
\begin{equation*}
\zeta\left(\frac{x_i}{x_j}\right)=\frac{\left[\frac{x_j}{q t x_i}\right]^{\delta_{j-1}^i}\left[\frac{t x_j}{q x_i}\right]^{\delta_{j+1}^i}}{\left[\frac{x_j}{x_i}\right]^{\delta_j^i}\left[\frac{x_j}{q^2 x_i}\right]^{\delta_j^i}},\qquad [x]:=x^{1/2}-x^{-1/2}. 
\end{equation*}

The standard coproduct structure $\Delta: U_{q, t}\left(\hat{\hat{\mathfrak{sl}}}_n\right) \longrightarrow U_{q, t}\left(\hat{\hat{\mathfrak{sl}}}_n\right) \widehat{\otimes} U_{q, t}\left(\hat{\hat{\mathfrak{sl}}}_n\right)$ is defined by:

\begin{equation*}
\begin{array}{ll}
\Delta\left(e_i^{+}(z)\right)=\varphi_i^{+}(z) \otimes e_i^{+}(z)+e_i^{+}(z) \otimes 1 & \Delta\left(\varphi_i^{+}(z)\right)=\varphi_i^{+}(z) \otimes \varphi_i^{+}(z) \\
\Delta\left(e_i^{-}(z)\right)=1 \otimes e_i^{-}(z)+e_i^{-}(z) \otimes \varphi_i^{-}(z) & \Delta\left(\varphi_i^{-}(z)\right)=\varphi_i^{-}(z) \otimes \varphi_i^{-}(z).
\end{array}
\end{equation*}

\subsection{Shuffle algebra and slope factorisation}
We now review the construction of the quantum toroidal algebra through the shuffle algebra and the resulting slope factorisation. For details, see \cite{N15}.

Fix the field $\mbb{F}=\mbb{Q}(q,t)$ and consider the following space of symmetric rational functions:
\begin{align*}
\widehat{\text{Sym}}(V):=\bigoplus_{\mbf{k}=(k_1,\cdots,k_{n})\in\mbb{N}^n}\mbb{F}(\cdots,z_{i1},\cdots,z_{ik_i},\cdots)^{\text{Sym}}_{1\leq i\leq n}.
\end{align*}
Here ``Sym'' denotes symmetrisation with respect to the variables $z_{i1},\cdots,z_{ik_i}$ of each colour $i$. We endow this vector space with the shuffle product
\begin{align*}
F*G=\text{Sym}[\frac{F(\cdots,z_{ia},\cdots)G(\cdots,z_{jb},\cdots)}{\mbf{k}!\mbf{l}!}\prod_{1\leq a\leq k_i}^{1\leq i\leq n}\prod_{k_{j}+1\leq b\leq k_{j}+l_{j}}^{1\leq j\leq n}\zeta(\frac{z_{ia}}{z_{jb}})].
\end{align*}

We define the subspace $\mc{S}^{\pm}_{\mbf{k}}\subset\widehat{\text{Sym}}(V)$ by:
\begin{align*}
\mc{S}^{+}:=\{F(\cdots,z_{ia},\cdots)=\frac{r(\cdots,z_{ia},\cdots)}{\prod_{1\leq a\neq b\leq k_i}^{1\leq i\leq n}(qz_{ia}-q^{-1}z_{ib})}\}
\end{align*}
where $r(\cdots,z_{i1},\cdots,z_{ik_i},\cdots)^{1\leq i\leq n}_{1\leq a\leq k_i}$ is any symmetric Laurent polynomial satisfying the wheel conditions
\begin{align*}
r(\cdots,q^{-1},t^{\pm},q,\cdots)=0
\end{align*}
for any three variables of colors $i,\cdots,i\pm1,i$.

The shuffle algebra $\mc{S}^+$ is naturally bigraded by the vector $\mbf{k}\in\mbb{N}^n$ recording the number of variables of each color and by the homogeneous degree $d\in\mbb{Z}$ of the rational function:
\begin{align*}
\mc{S}^+=\bigoplus_{(\mbf{k},d)\in\mbb{N}^n\times\mbb{Z}}\mc{S}^+_{\mbf{k},d}.
\end{align*}

Similarly, we define the negative shuffle algebra $\mc{S}^{-}:=(\mc{S}^{+})^{op}$ by equipping the same vector space with the opposite shuffle product. Its grading is $\mc{S}^-=\bigoplus_{(\mbf{k},d)\in\mbb{N}^n\times\mbb{Z}}\mc{S}^-_{-\mbf{k},d}$. We then enlarge the positive and negative shuffle algebras by adjoining the generators $\{\varphi_{i,d}^{\pm}\}_{1\leq i\leq n}^{d\geq0}$:
\begin{align*}
\mc{S}^{\geq}=\langle\mc{S}^+,\{(\varphi^+_{i,d})^{d\geq0}_{1\leq i\leq n}\}\rangle,\qquad \mc{S}^{\leq}=\langle\mc{S}^-,\{(\varphi^-_{i,d})^{d\geq0}_{1\leq i\leq n}\}\rangle.
\end{align*}

The elements $\varphi^{\pm}_{i,d}$ commute among themselves and satisfy the following relations with $\mc{S}^{\pm}$:
\begin{align*}
\varphi^{+}_i(w)F=F\varphi^+_i(w)\prod_{1\leq a\leq k_j}^{1\leq j\leq n}\frac{\zeta(w/z_{ja})}{\zeta(z_{ja}/w)}
\end{align*}
\begin{align*}
\varphi^{-}_i(w)G=G\varphi^-_i(w)\prod_{1\leq a\leq k_j}^{1\leq j\leq n}\frac{\zeta(z_{ja}/w)}{\zeta(w/z_{ja})}.
\end{align*}

There is a nondegenerate bialgebra pairing:
\begin{align}\label{Drinfeld-pairing}
\langle-,-\rangle:\mc{S}^{\leq}\otimes\mc{S}^{\geq}\rightarrow\mbb{Q}(q,t)
\end{align} 
between the negative and positive shuffle algebras $\mc{S}^{\leq,\geq}$ is given by \cite{N15}:
\begin{align*}
&\langle\varphi^{-}_{i}(z),\varphi^{+}_{j}(w)\rangle=\frac{\zeta_{ij}(w/z)}{\zeta_{ji}(z/w)}\\
&\langle G,F\rangle=\frac{1}{\mbf{k}!}\int_{\lvert z_{ia}\lvert=1}^{\lvert q\lvert<1\lvert p\lvert}\frac{G(\cdots,z_{ia},\cdots)F(\cdots,z_{ia},\cdots)}{\prod_{a\leq k_i,b\leq k_j}^{1\leq i,j\leq n}\zeta_{p}(z_{ia}/z_{jb})}\prod_{1\leq a\leq k_i}^{1\leq i\leq n}\frac{dz_{ia}}{2\pi iz_{ia}}|_{p\mapsto q}
\end{align*}
for $F\in\mc{S}^{+}$, $G\in\mc{S}^{-}$. Here $\zeta_{p}$ is defined as follows:
\begin{align*}
\zeta_{p}(\frac{x_i}{x_j})=\zeta(\frac{x_i}{x_j})\frac{[\frac{x_j}{p^2x_i}]^{\delta^{i}_{j}}}{[\frac{x_j}{q^2x_i}]^{\delta^{i}_{j}}}.
\end{align*}

The pairing is called the Drinfeld pairing, and for the construction one can refer to Chapter one in \cite{N15}. This defines the double shuffle algebra $\mc{S}:=\mc{S}^{\leq}\hat{\otimes}\mc{S}^{\geq}$. Its coproduct $\Delta:\mc{S}\rightarrow\mc{S}\hat{\otimes}\mc{S}$ is given by
\begin{align}\label{Drinfeld-coproduct}
&\Delta(\varphi_{i}^{\pm}(w))=\varphi^{\pm}_{i}(w)\otimes\varphi^{\pm}_{i}(w)\\
&\Delta(F)=\sum_{0\leq\mbf{l}\leq\mbf{k}}\frac{[\prod_{1\leq j\leq n}^{b>l_{j}}\varphi_{j}^+(z_{jb})]F(\cdots,z_{i1},\cdots,z_{il_i}\otimes z_{i,l_i+1},\cdots,z_{ik_i},\cdots)}{\prod_{1\leq i\leq n}^{a\leq l_i}\prod_{1\leq j\leq n}^{b>l_j}\zeta(z_{jb}/z_{ia})},\qquad F\in\mc{S}^+\\
&\Delta(G)=\sum_{0\leq\mbf{l}\leq\mbf{k}}\frac{G(\cdots,z_{i1},\cdots,z_{il_i}\otimes z_{i,l_i+1},\cdots,z_{ik_i},\cdots)[\prod_{1\leq j\leq n}^{b>l_{j}}\varphi_{j}^-(z_{jb})]}{\prod_{1\leq i\leq n}^{a\leq l_i}\prod_{1\leq j\leq n}^{b>l_j}\zeta(z_{ia}/z_{jb})},\qquad G\in\mc{S}^-.
\end{align}

In concrete terms, the right-hand side of the coproduct formula is expanded in the region $\lvert z_{ia}\rvert\ll\lvert z_{jb}\rvert$ for all $a\leq l_i$ and $b>l_i$. Monomials in $\{z_{ia}\}_{a\leq l_i}$ are placed to the left of the tensor sign, while monomials in $\{z_{jb}\}_{b>l_j}$ are placed to its right. The shuffle algebra also has an antipode $S:\mc{S}\rightarrow\mc{S}$, which is an anti-homomorphism for both algebras and coalgebras:
\begin{align*}
&S(\varphi^+_i(z))=(\varphi^+_i(z))^{-1},\qquad S(F)=[\prod_{1\leq a\leq k_i}^{1\leq i\leq n}(-\varphi_i^+(z_{ia}))^{-1}]*F\\
&S(\varphi^-_i(z))=(\varphi^-_i(z))^{-1},\qquad S(G)=G*[\prod_{1\leq a\leq k_i}^{1\leq i\leq n}(-\varphi_i^-(z_{ia}))^{-1}].
\end{align*}

The following theorem has been proved in \cite{N15}:
\begin{thm}
There is an isomorphism of bialgebras
\begin{align*}
Y:U_{q,t}(\hat{\hat{\mf{sl}}}_{n})\rightarrow\mc{S}
\end{align*}
given by:
\begin{align*}
Y(\varphi_{i,d}^{\pm})=\varphi^{\pm}_{i,d},\qquad Y(e_{i,d}^{\pm})=\frac{z_{i1}^d}{[q^{-2}]}.
\end{align*}
\end{thm}

\subsubsection{Slope subalgebra of the shuffle algebra}

We first review the definition of the slope subalgebra of the quantum toroidal algebra $U_{q,t}(\hat{\hat{\mf{sl}}}_n)$ from Section $5$ of \cite{N15}.

For a fixed vector $\mbf{m}\in\mbb{R}^n$, we define the subspace of the positive shuffle algebra $\mc{S}^+$ of slope at most $\mbf{m}$ by
\begin{align*}
\mc{S}^{+}_{\leq\mbf{m}}:=\{F\in\mc{S}^+|\lim_{\xi\rightarrow\infty}\frac{F(z_{i1},\cdots,z_{ia_i},\xi z_{i,a_i+1},\cdots,\xi z_{i,a_i+b_i})}{\xi^{\mbf{m}\cdot\mbf{b}-\frac{(\mbf{a},\mbf{b})}{2}}}<\infty\text{ }\forall\mbf{a},\mbf{b}\in\mbb{N}^n\}.
\end{align*}

Similarly, for the negative half:
\begin{align*}
\mc{S}^{-}_{\leq\mbf{m}}=\{G\in\mc{S}^-|\lim_{\xi\rightarrow0}\frac{G(\xi z_{i1},\cdots,\xi z_{ia_i},z_{i,a_i+1},\cdots,z_{i,a_i+b_i})}{\xi^{-\mbf{m}\cdot\mbf{a}+\frac{(\mbf{a},\mbf{b})}{2}}}<\infty\text{ }\forall\mbf{a},\mbf{b}\in\mbb{N}^n\}
\end{align*}
where $(\mbf{a},\mbf{b}):=\sum_{i\in\mbb{Z}/n\mbb{Z}}2a_ib_i-a_{i}b_{i-1}-a_ib_{i+1}$.

From Section $5.2$ in \cite{N15}, one can see that for the elements $F\in\mc{S}_{\leq\mbf{m}}^{+}\cap\mc{S}^+_{\mbf{k},d}$, $G\in\mc{S}_{\leq\mbf{m}}^{-}\cap\mc{S}^-_{-\mbf{k},d}$, one has the following expansion for their coproducts:
\begin{align*}
\Delta(F)=\Delta_{\mbf{m}}(F)+(\text{anything})\otimes(\text{slope}<\mbf{m}),\qquad\Delta_{\mbf{m}}(F)=\sum_{\mbf{a}+\mbf{b}=\mbf{k}}\lim_{\xi\rightarrow\infty}\frac{\Delta(F)_{\mbf{a}\otimes\mbf{b}}}{\xi^{\mbf{m}\cdot\mbf{b}}}
\end{align*}

\begin{align*}
\Delta(G)=\Delta_{\mbf{m}}(G)+(\text{slope}<\mbf{m})\otimes(\text{anything}),\qquad\Delta_{\mbf{m}}(G)=\sum_{\mbf{a}+\mbf{b}=\mbf{k}}\lim_{\xi\rightarrow\infty}\frac{\Delta(G)_{\mbf{a}\otimes\mbf{b}}}{\xi^{-\mbf{m}\cdot\mbf{a}}}.
\end{align*}

The positive and negative slope subalgebras are therefore defined by
\begin{align*}
\mc{B}^{\pm}_{\mbf{m}}:=\bigoplus_{(\mbf{k},d)\in\mbb{N}^n\times\mbb{Z}}^{\mbf{m}\cdot\mbf{k}=d}\mc{S}^{\pm}_{\leq\mbf{m}}\cap\mc{S}^{\pm}_{\pm\mbf{k},d}.
\end{align*}

The coproduct $\Delta_{\mbf{m}}$ descends to a genuine coproduct on $\mc{B}_{\mbf{m}}^{\geq,\leq}:=\mc{B}_{\mbf{m}}^{\pm}\otimes\mbb{Q}(q,t)[\varphi_i^{\pm1}]$. Taking the Drinfeld double in the shuffle-algebra setting then produces the slope subalgebra $\mc{B}_{\mbf{m}}\subset\mc{S}$.

From Theorem V.4 in \cite{N15}, the slope subalgebra $\mc{B}_{\mbf{m}}$ for the quantum toroidal algebra $U_{q,t}(\widehat{\widehat{\mf{sl}}}_{n})$ is isomorphic to:
\begin{align}\label{rootquantum}
\mc{B}_{\mbf{m}}\cong\bigotimes^{g}_{h=1}U_{q}(\widehat{\mf{gl}}_{l_h}).
\end{align}

In Subsection \ref{section4-3}, we give a simple rule for determining $l_h$ and $g$ as described in Section 5.3 in \cite{N15}. For the present section, we describe the slope subalgebra $\mc{B}_{\mbf{m}}$ through its generators.

Let us define the generators of the slope algebra $\mc{B}_{\mbf{m}}$. For $\mbf{m}\cdot[i;j)\in\mbb{Z}$, we denote the following elements:
\begin{equation}\label{positive-generators}
P_{\pm[i ; j)}^{\pm \mathbf{m}}=\operatorname{Sym}\left[\frac{\prod_{a=i}^{j-1} z_a^{\left\lfloor m_i+\ldots+m_a\right\rfloor-\left\lfloor m_i+\ldots+m_{a-1}\right\rfloor}}{t^{\text {ind } d_{[i, j)}^m q^{i-j}} \prod_{a=i+1}^{j-1}\left(1-\frac{q_2 z_a}{z_{a-1}}\right)} \prod_{i \leq a<b<j} \zeta\left(\frac{z_b}{z_a}\right)\right]\in\mc{B}_{\mbf{m}}^{\pm}
\end{equation}
\begin{equation}\label{negative generators}
Q_{\mp[i ; j)}^{\pm \mathbf{m}}=\operatorname{Sym}\left[\frac{\prod_{a=i}^{j-1} z_a^{\left\lfloor m_i+\ldots+m_{a-1}\right\rfloor-\left\lfloor m_i+\ldots+m_a\right\rfloor}}{t^{-\mathrm{ind}_{[i, j)}^m} \prod_{a=i+1}^{j-1}\left(1-\frac{q_1 z_{a-1}}{z_a}\right)} \prod_{i \leq a<b<j} \zeta\left(\frac{z_a}{z_b}\right)\right]\in\mc{B}_{\mbf{m}}^{\pm}.
\end{equation}

Here $\text{ind}^{\mbf{m}}_{[i;j)}$ is defined as:
\begin{align*}
\text{ind}^{\mbf{m}}_{[i;j)}=\sum_{a=i}^{j-1}(m_i+\cdots+m_a-\lfloor m_i+\cdots+m_{a-1}\lfloor).
\end{align*}

It is proved in \cite{N15} that the positive and negative parts $\mc{B}_{\mbf{m}}^{\pm}$ are generated by $\{P^{\mbf{m}}_{\pm[i;j)}\}_{i\leq j}$ and $\{Q^{\mbf{m}}_{\pm[i;j)}\}_{i\leq j}$. Moreover, the restriction of the bialgebra pairing \eqref{Drinfeld-pairing} to $\mc{B}_{\mbf{m}}^{\leq,\geq}$ equips $\mc{B}_{\mbf{m}}$ with a Hopf algebra structure satisfying the following relations.

For the antipode map $S_{\mbf{m}}$:
\begin{align*}
S_{\mbf{m}}(P^{\mbf{m}}_{[i;j)})=Q^{\mbf{m}}_{[i;j)},\qquad S_{\mbf{m}}(Q^{\mbf{m}}_{-[i;j)})=P^{\mbf{m}}_{-[i;j)}.
\end{align*}

For the coproduct $\Delta_{\mbf{m}}$:
\begin{equation*}
\Delta_{\mathbf{m}}\left(P_{[i ; j)}^{\mathbf{m}}\right)=\sum_{a=i}^j P_{[a ; j)}^{\mathbf{m}} \varphi_{[i ; a)} \otimes P_{[i ; a)}^{\mathbf{m}} \quad \Delta_{\mathbf{m}}\left(Q_{[i ; j)}^{\mathrm{m}}\right)=\sum_{a=i}^j Q_{[i ; a)}^{\mathbf{m}} \varphi_{[a ; j)} \otimes Q_{[a ; j)}^{\mathbf{m}}
\end{equation*}

\begin{equation*}
\Delta_{\mathbf{m}}\left(P_{-[i ; j)}^{\mathbf{m}}\right)=\sum_{a=i}^j P_{-[a ; j)}^{\mathbf{m}} \otimes P_{-[i ; a)}^{\mathbf{m}} \varphi_{-[a ; j)} \quad \Delta_{\mathbf{m}}\left(Q_{-[i ; j)}^{\mathbf{m}}\right)=\sum_{a=i}^j Q_{-[i ; a)}^{\mathbf{m}} \otimes Q_{-[a ; j)}^{\mathbf{m}} \varphi_{-[i ; a)}.
\end{equation*}

Here $\Delta_{\mbf{m}}$ is defined as $(5.27)$ and $(5.28)$ in \cite{N15}. The isomorphism \ref{rootquantum} is given by:
\begin{align*}
e_{[i;j)}=P^{\mbf{m}}_{[i;j)_{h}},\qquad e_{-[i;j)}=Q^{\mbf{m}}_{-[i;j)_h},\qquad\varphi_{k}=\varphi_{[k;v_{\mbf{m}}(k))}.
\end{align*}
Here the number $v_{\mbf{m}}(k)$ is defined after \eqref{oriented-cycles}. The bialgebra pairing determines a universal $R$-matrix $R_{\mbf{m}}\in\mc{B}^{\leq}_{\mbf{m}}\hat{\otimes}\mc{B}_{\mbf{m}}^{\geq}$ for $\mc{B}_{\mbf{m}}$.

It is worth noting that there is no immediate explicit formula for the Heisenberg operators $U_{q}(\hat{\mf{gl}}_1)$ inside $U_{q}(\hat{\mf{gl}}_n)$.
Here we show the explicit formula for the simplest primitive elements of $\mc{B}_{\mbf{m}}$ for $\mbf{m}=\mu\bm{\theta}$. Let $k=(j-i)\mu$ with $\mu$ a rational number. We assume that for $0\leq n\leq j-i-1$, $\mu n\notin\mbb{Z}$, in this case we have
\begin{align*}
&\Delta_{\mbf{m}}(P^{\mbf{m}}_{[i;j)})=P^{\mbf{m}}_{[i;j)}\otimes1+\varphi_{[i;j)}\otimes P^{\mbf{m}}_{[i;j)},\qquad\Delta_{\mbf{m}}(Q^{\mbf{m}}_{[i;j)})=Q^{\mbf{m}}_{[i;j)}\otimes1+\varphi_{[i;j)}\otimes Q^{\mbf{m}}_{[i;j)}\\
&\Delta_{\mbf{m}}(P^{\mbf{m}}_{-[i;j)})=P^{\mbf{m}}_{-[i;j)}\otimes\varphi_{[i;j)}+1\otimes P^{\mbf{m}}_{[i;j)},\qquad\Delta_{\mbf{m}}(Q^{\mbf{m}}_{-[i;j)})=Q^{\mbf{m}}_{-[i;j)}\otimes\varphi_{[i;j)}+1\otimes Q^{\mbf{m}}_{[i;j)}.
\end{align*}

Now take $j=i+l\bm{\theta}$. Since $\varphi_{[i;i+l\bm{\theta})}=\varphi_1^{l}\cdots\varphi_{n}^{l}=c^l$, the elements $P^{\mbf{m}}_{\pm[i;j)}$ are primitive. If $\mu=p/q$ with $\gcd(p,q)=1$, the first condition for primitivity is $ln\leq q$. This provides a family of primitive elements among the generators $P^{\mbf{m}}_{\pm[i;j)}$.

\subsection{Other coproduct structure on $U_{q,t}(\hat{\hat{\mf{sl}}}_{n})$}

In \cite{N15}, Negut proved that there is a factorisation of $U_{q,t}(\hat{\hat{\mf{sl}}}_{n})$ into the product of the slope subalgebras:
\begin{align}\label{slope-factorisation}
U_{q,t}(\hat{\hat{\mf{sl}}}_{n})=\bigotimes_{m\in\mbb{R}}^{\rightarrow}\mc{B}_{s+m\bm{\theta}}^{+}\otimes\mc{B}_{\infty\bm{\theta}}\otimes\bigotimes_{m\in\mbb{R}}^{\rightarrow}\mc{B}_{s+m\bm{\theta}}^{-}
\end{align}
which preserves the Hopf pairing for arbitrary $s\in\mbb{Q}^{n}$ and $\bm{\theta}\in(\mbb{Q}^+)^{n}$. Moreover, using the same argument in \cite{N15}, we can also show that we have the following generalisation of the factorization \ref{slope-factorisation}:
\begin{align*}
U_{q,t}^{+}(\hat{\hat{\mf{sl}}}_n)\cong\bigotimes_{m\in\mbb{R}}^{\rightarrow}\mc{B}_{s+\bm{\mu}(m)}^{+}
\end{align*}
and here $\bm{\mu}:\mbb{R}\rightarrow\mbb{R}^n$ is a continuous function such that $\mu(t_1)>\mu(t_2)$ if $t_1>t_2$.

It is also shown in \cite{N15} that this isomorphism follows from the factorisation of the pieces of slope at most $\mbf{m}$:
\begin{align}\label{slope-m-factorisation}
\mc{S}^{+}_{\leq\mbf{m}}\cong\bigotimes_{m\in\mbb{R}_{\geq0}}^{\rightarrow}\mc{B}_{\mbf{m}+\bm{\mu}(m)}^{+},\qquad\bm{\mu}(0)=0.
\end{align}

Each slope subalgebra $\mc{B}_{\mbf{m}}$ admits the isomorphism given in \eqref{rootquantum}. The algebra $\mathcal{B}_{\infty \cdot \boldsymbol{\theta}}=\left\langle\varphi_{i, d}^{ \pm}\right\rangle_{1 \leq i \leq n}^{d \in \mathbb{N}_0}$ is the Cartan subalgebra. This factorisation yields the following Khoroshkin--Tolstoy-type factorisation of the universal $R$-matrix of $U_{q,t}(\hat{\hat{\mf{sl}}}_{r})$:
\begin{align}\label{wall-crossing-factorization-for-universal-R-matrices}
R^{\infty}=\prod_{m\in\mbb{R}}^{\rightarrow}(\hat{R}_{\bm{\mu}(m)})R_{\infty}=(\prod_{m\in\mbb{R}}^{\rightarrow}\prod_{h=1}^{g^{(m)}}\hat{R}_{U_{q}(\hat{\mf{gl}}_{h}^{(m)})})\cdot R_{\text{Heisenberg}}^{\otimes n}.
\end{align}

The meaning of the infinite product in \eqref{wall-crossing-factorization-for-universal-R-matrices} is that if we express $R^{\infty}=\sum_{\alpha}F_{\alpha}\otimes G_{\alpha}$, and each $F_{\alpha}$ and $G_{\alpha}$ is a finite linear combination of the finite products of elements in the slope subalgebras $\mc{B}_{\mbf{m}}$.

In the above formula, each $R_{\mbf{m}}$ factorises as $R_{\mbf{m}}=\hat{R}_{\mbf{m}}\hat{K}$, where $\hat{K}$ is the diagonal part, and $\hat{R}_{\mbf{m}}$ is called the \textbf{reduced part of the universal $R$-matrix}. Moreover, the universal $R$-matrix $R^{\infty}$ is independent of the choice of the continuous function $\bm{\mu}$, because it depends only on the bialgebra pairing \eqref{Drinfeld-pairing} of the quantum toroidal algebra.

For an arbitrary $\mbf{m}\in\mbb{R}^n$, we can define a new coproduct on $U_{q,t}(\hat{\hat{\mf{sl}}}_n)$ by the following:
\begin{align}\label{twisted-coproduct}
\Delta_{(\mbf{m})}(a)=[\prod^{\rightarrow}_{r\in\mbb{Q}_{>0}\cup\{\infty\}}\hat{R}_{\mbf{m}+\bm{\mu}(r)}^+]^{-1}\cdot\Delta(a)\cdot[\prod^{\rightarrow}_{r\in\mbb{Q}_{>0}\cup\{\infty\}}\hat{R}_{\mbf{m}+\bm{\mu}(r)}^+].
\end{align}
with $\bm{\mu}:\mbb{R}^n\rightarrow\mbb{R}$ an increasing monotone path such that $\bm{\mu}(0)=0$, 
and here the meaning of the infinite product $\prod^{\rightarrow}_{r\in\mbb{Q}_{>0}\cup\{\infty\}}\hat{R}_{\mbf{m}+\bm{\mu}(r)}^+$ should be understood the same as the one in \eqref{wall-crossing-factorization-for-universal-R-matrices}, i.e. it is the sum of the linear combination of the finite product of elements in the slope subalgebras $\mc{B}_{\mbf{m}+\bm{\mu}(r)}$ with $r\geq0$.

With respect to this coproduct, the quantum toroidal algebra admits the slope factorisation
\begin{align*}
U_{q,t}(\hat{\hat{\mf{sl}}}_n)\cong\bigotimes_{\mu\in\mbb{Q}^-}^{\leftarrow}\mc{B}_{\mbf{m}+\bm{\mu}(r)}\otimes\bigotimes_{\mu\in\mbb{Q}^+\sqcup\{\infty\}}^{\rightarrow}\mc{B}_{\mbf{m}+\bm{\mu}(r)}.
\end{align*}

\begin{prop}
The coproduct $\Delta_{(\mbf{m})}$ is independent of the choice of the increasing continuous function $\bm{\mu}$.
\end{prop}
\begin{proof}
It suffices to prove that the infinite product $\prod^{\rightarrow}_{r\in\mbb{Q}_{>0}\cup\{\infty\}}\hat{R}_{\mbf{m}+\bm{\mu}(r)}^+$ is independent of $\bm{\mu}$. This follows from the construction below.

By \ref{slope-m-factorisation}, the restriction of the Drinfeld pairing \eqref{Drinfeld-pairing} of the quantum toroidal algebra pairs $\mc{S}_{\leq\mbf{m}}^{+}$ with $\mc{S}_{\leq\mbf{m}}^{-}$ and $\mbb{Q}(q,t)[\varphi_{i,0}^{\pm1}]$. Note that the operator
\begin{align}\label{partial-universal-m-R-matrix}
\prod^{\rightarrow}_{r\in\mbb{Q}_{>0}\cup\{\infty\}}\hat{R}_{\mbf{m}+\bm{\mu}(r)}^+\cdot K.
\end{align}
is written as $(\sum_{\alpha}F_{\alpha}\otimes G_{\alpha})\cdot K$ with $F_{\alpha}, G_{\alpha}$ the corresponding orthogonal basis in $S_{\leq\mbf{m}}^{\pm}$ respectively. Here $K$ is the universal $R$-matrix induced by $\mbb{Q}(q,t)[\varphi_{i,0}^{\pm}]$ via identifying $\varphi_{i,0}=\exp(\hbar b_i)$. By the standard result of the universal $R$-matrix (For example, Exercise V.14 in \cite{N15}), the operator \eqref{partial-universal-m-R-matrix} is independent of the choice of the orthogonal basis in $S_{\leq\mbf{m}}^{\pm}$ respectively.

The universal $R$-matrix of $\mc{S}_{\mbf{m}}$ depends only on the Drinfeld pairing, whose restriction to each slope subalgebra $\mc{B}_{\mbf{m}+\bm{\mu}(m)}^{\pm}$ is again a Drinfeld pairing. Different paths in \ref{slope-m-factorisation} therefore produce only different choices of basis. Since the universal $R$-matrix is basis-independent, it is independent of the function $\bm{\mu}$.

\end{proof}

Now we prove the following important proposition:
\begin{prop}\label{twist-of-the-coproduct}
$\Delta_{(\mbf{m})}(a)=\Delta_{\mbf{m}}(a)$ when $a\in\mc{B}_{\mbf{m}}$.
\end{prop}
\begin{proof}
The factorisation $U_{q,t}^{\pm}(\hat{\hat{\mf{sl}}}_n)\cong\bigotimes_{t\in\mbb{R}}^{\rightarrow}\mc{B}_{\bm{\mu}(t)}^{\pm}$ means that a basis of $U_{q,t}^{\pm}(\hat{\hat{\mf{sl}}}_n)$ may be represented by finite products of basis elements from the subalgebras $\mc{B}_{\bm{\mu}(t)}^{\pm}$. Accordingly, the universal $R$-matrix of $U_{q,t}(\hat{\hat{\mf{sl}}}_n)$ with respect to the coproduct $\Delta_{(\mbf{m})}$ factorises as (See Theorem 2.3.4 in \cite{Maj95}):
\begin{align*}
R^{\mbf{m}}=\prod^{\leftarrow}_{t\in\mbb{R}_{>0}\cup\{\infty\}}\hat{R}_{\mbf{m}+\mu(t)}^{-}\prod^{\rightarrow}_{t\in\mbb{R}_{\leq0}}\hat{R}_{\mbf{m}+\mu(t)}^+
\end{align*}
Here $\hat{R}_{\mbf{m}}^{-}:=\hat{R}_{\mbf{m}}^{21}$. The factorisation is understood so that each term of $R^{\mbf{m}}$ contains only finitely many non-unit factors from the $R$-matrix $\hat{R}_{\mbf{m}+\mu(t)}^{\pm}$.

For each $\mc{B}^{\pm}_{\mbf{m}+\mu(t)}$, choose dual bases $\{F_{(\mbf{m}+\mu(t))}^{m}\}$ and $\{G_{(\mbf{m}+\mu(t))}^{m}\}$. We adopt the convention that $F_{(\mbf{m}+\mu(t))}^{m}\in\mc{B}_{\mbf{m}+\mu(t)}^{-}$ when $t>0$ or $t=\infty$, and $F_{(\mbf{m}+\mu(t))}^{m}\in\mc{B}_{\mbf{m}+\mu(t)}^{+}$ when $t\leq0$. Here $m\in S_{\mbf{m}+\mu(t)}$, where $S_{\mbf{m}+\mu(t)}$ indexes a basis of $\mc{B}^{\pm}_{\mbf{m}+\mu(t)}$. In terms of these bases, the reduced universal $R$-matrix of $\mc{B}_{\mbf{m}+\mu(t)}$ is
\begin{align*}
\hat{R}_{\mbf{m}+\mu(t)}=
\begin{cases}
\sum_{m\in S_{\mbf{m}+\mu(t)}}F_{(\mbf{m}+\mu(t)}^{m}\otimes G_{(\mbf{m}+\mu(t))}^{m}&t\leq0\\
\sum_{m\in S_{\mbf{m}+\mu(t)}}G_{(\mbf{m}+\mu(t)}^{m}\otimes F_{(\mbf{m}+\mu(t))}^{m}&t>0\text{ or }t=\infty
\end{cases}.
\end{align*}

Here we denote $m=0$ for the basis $1\otimes 1$.

The factorisation implies that the sets $\{\prod_{t\in\mbb{R}_{>0}\sqcup\infty\sqcup\mbb{R}_{\leq0}}F_{(\mbf{m}+\bm{\mu}(t))}^{m_t}\}$ and $\{\prod_{t\in\mbb{R}_{>0}\sqcup\infty\sqcup\mbb{R}_{\leq0}}G_{(\mbf{m}+\bm{\mu}(t))}^{m_t}\}$ form bases of $U_{q,t}^{\pm}(\hat{\hat{\mf{sl}}}_n)$. Here the product notation $\prod$ only means that we only take finitely many elements in the set $\mbb{R}_{>0}\sqcup\infty\sqcup\mbb{R}_{\leq0}$.  Since $\mbf{m}$ is the final element in the ordering, we write these bases as $\{\prod_{t\in\mbb{R}_{>0}\sqcup\infty\sqcup\mbb{R}_{<0}}F_{(\mbf{m}+\bm{\mu}(t))}^{m_t}F_{(\mbf{m})}^{m}\}$ and $\prod_\{\prod_{t\in\mbb{R}_{>0}\sqcup\infty\sqcup\mbb{R}_{\leq0}}G_{(\mbf{m}+\bm{\mu}(t))}^{m_t}G_{(\mbf{m})}^{m}\}$.

Now we denote $I:=\mbb{R}_{>0}\sqcup\{\infty\}\sqcup\mbb{R}_{<0}$. We also denote $u:=\mbf{m}\in I$. In this case, we denote $R_{i}:=\hat{R}_{\mbf{m}+\mu(t)}$. Thus $R^{\mbf{m}}$ may be written as
\begin{align*}
R^{\mbf{m}}=\prod_{i\in I}^{\rightarrow}R_{i}=\prod_{i\in I}[\sum_{m_i\in S_i}F^{m_i}_{(i)}\otimes G^{m_i}_{(i)}]=\sum_{\substack{m_i\in S_i}}\prod_{i\in I}F_{(i)}^{m_i}\otimes \prod_{i\in I}G_{(i)}^{m_i}.
\end{align*}

To prove the statement, we need to use the following formula:
\begin{align*}
(\Delta_{(\mbf{m})}\otimes1)R^{\mbf{m}}=R_{13}^{\mbf{m}}R_{23}^{\mbf{m}},\qquad (1\otimes\Delta_{(\mbf{m})})R^{\mbf{m}}=R_{12}^{\mbf{m}}R_{13}^{\mbf{m}}.
\end{align*}

Now we write the product $R_{13}R_{23}$ in the explicit form:
\begin{equation*}
\begin{aligned}
R_{13}^{\mbf{m}}R_{23}^{\mbf{m}}=&\sum_{\substack{m_i,n_i\in S_i\\i\in I}}(\prod_{i\in I}F_{(i)}^{m_i}F_{(u)}^{m_u}\otimes 1\otimes\prod_{i\in I}G_{(i)}^{m_i}G_{(u)}^{m_u})(1\otimes\prod_{i\in I}F_{(i)}^{n_i}F_{(u)}^{n_u}\otimes\prod_{i\in I}G_{(i)}^{n_i}G_{(u)}^{n_u})\\
=&\sum_{\substack{m_i,n_i\in S_i\\i\in I}}\prod_{i\in I}F_{(i)}^{m_i}F_{(u)}^{m_u}\otimes\prod_{i\in I}F_{(i)}^{n_i}F_{(u)}^{n_u}\otimes\prod_{i\in I}G_{(i)}^{m_i}G_{(u)}^{m_u}\prod_{i\in I}G_{(i)}^{n_i}G_{(u)}^{n_u})\\
=&\sum_{\substack{m_i,n_i,l_i\in S_i\\i\in I}}\prod_{i\in I}F_{(i)}^{m_i}F_{(u)}^{m_u}\otimes\prod_{i\in I}F_{(i)}^{n_i}F_{(u)}^{n_u}\otimes\prod_{i\in I}G_{(i)}^{m_i}G_{(u)}^{m_u}\prod_{i\in I}G_{(i)}^{n_i}G_{(u)}^{n_u})\\
=&\sum_{\substack{m_i,n_i,l_i\in S_i\\i\in I}}a_{\mbf{m}',\mbf{n}}^{\mbf{l}}\prod_{i\in I}F_{(i)}^{m_i}F_{(u)}^{m_u}\otimes\prod_{i\in I}F_{(i)}^{n_i}F_{(u)}^{n_u}\otimes\prod_{i\in I}G_{(i)}^{l_i}G_{(u)}^{l_u}.
\end{aligned}
\end{equation*}

Here $\mbf{l}:=\{l_i\}_{i\in I},\mbf{n}:=\{n_i\}_{i\in I},\mbf{m}':=\{m_i\}_{i\in I}$, and 
\begin{align*}
a_{\mbf{m},\mbf{n}}^{\mbf{l}}=\langle\prod_{i\in I}F_{(i)}^{l_i}F_{(u)}^{l_u},\prod_{i\in I}G_{(i)}^{m_i}G_{(u)}^{m_u}\prod_{i\in I}G_{(i)}^{n_i}G_{(u)}^{n_u})\rangle
\end{align*}
is well-defined since $\mc{B}_{\mbf{m}}^{+}$ is $\mbb{N}^n$-graded finite-dimensional and we only allow finitely many terms in the product for the expression of $R^{\mbf{m}}$.

The coproduct $(\Delta_{(\mbf{m})}\otimes1)R^{\mbf{m}}$ may now be written as
\begin{equation*}
(\Delta_{(\mbf{m})}\otimes1)R^{\mbf{m}}=\sum_{\substack{l_i\in S_i\\i\in I}}\Delta_{(\mbf{m})}(\prod_{i\in I}F_{(i)}^{l_i}F_{(u)}^{l_u})\otimes\prod_{i\in I}G_{(i)}^{l_i}G_{(u)}^{l_u}.
\end{equation*}

This gives the general coproduct formula:
\begin{align*}
\Delta_{(\mbf{m})}(\prod_{i\in I}F_{(i)}^{l_i}F_{(u)}^{l_u})=\sum_{\substack{m_i,n_i\in S_i\\i\in I}}a_{\mbf{m},\mbf{n}}^{\mbf{l}}\prod_{i\in I}F_{(i)}^{m_i}F_{(u)}^{m_u}\otimes\prod_{i\in I}F_{(i)}^{n_i}F_{(u)}^{n_u}.
\end{align*}

We next compare the coproduct of $G_{(u)}^{m_u}$. Observe that
\begin{equation*}
\begin{aligned}
&\langle F_{(u)}^{l_u},\prod_{i\in I}G_{(i)}^{m_i}G_{(u)}^{m_u}\prod_{i\in I}G_{(i)}^{n_i}G_{(u)}^{n_u})\rangle\\
=&\sum_{r_i\in S_i,i\in I}a^{\mbf{r}}_{\mbf{m},\mbf{n}}\langle F_{(u)}^{l_u},\prod_{i\in I}G_{(i)}^{r_i}G_{(u)}^{r_u}\rangle\\
=&\sum_{r_i\in S_i,i\in I}a^{\mbf{r}}_{\mbf{m},\mbf{n}}\prod_{i\in I}\langle1,G_{(i)}^{r_i}\rangle\langle F^{l_u}_{(u)},G_{(u)}^{r_u}\rangle\\
=&\sum_{r_i\in S_i,i\in I}a^{\mbf{r}}_{\mbf{m},\mbf{n}}\delta_{l_u,r_u}\prod_{i\in I}\epsilon(G_{(i)}^{r_i}).
\end{aligned}
\end{equation*}

The product is nonzero only when every $G_{(i)}^{m_i}$ is either $1$ or of the form $\prod_{i=1}^{n}\varphi_i^{k_i}$. Since the elements of $\prod_{i=1}^{n}\varphi_i^{k_i}$ still lie in $\mc{B}_{i}$, this forces $G_{(i)}^{m_i}=G_{(i)}^{n_i}=1$. Hence the coproduct formula contains only elements of $\mc{B}_{\mbf{m}}$ and takes the form
\begin{align*}
\Delta_{\mbf{m}}(F_{(u)}^{l_u})=\sum_{m_u,n_u\in S_u}a^{l_u}_{m_u,n_u}F^{m_u}_{(u)}\otimes F^{n_u}_{(u)}
\end{align*}
which coincides with the coproduct $\Delta_{\mbf{m}}$ of $\mc{B}_{\mbf{m}}$.

Thus $\mbf{m}$ is the maximal element in this ordering, which completes the proof.
\end{proof}

A useful consequence for the universal $R$-matrix is the following.
\begin{prop}\label{finiteness-of-wall-set}
Fix a curve $\mu:\mbb{R}\rightarrow\mbb{R}^n$ such that $\mu(t_1)>\mu(t_2)$ for $t_1>t_2$ and consider a finite length segment $\mu^{-1}([t_1,t_2])$ on this curve. Then among operators $R_{\mu(t)}^+$, only finitely many reduced $R$-matrices $R_{\mu(t)}^+$ act nontrivially on
\[
K_T(M(\mbf{v}_1,\mbf{w}_1))\otimes K_T(M(\mbf{v}_2,\mbf{w}_2)).
\]
\end{prop}
\begin{proof}
This follows from the shuffle-algebra representation of each $\mc{B}_{\mbf{m}}$ and from the finite generation of
\[
K_T(M(\mbf{v}_1,\mbf{w}_1))\otimes K_T(M(\mbf{v}_2,\mbf{w}_2))
\]
as a $K_G(pt)$-module. The reduced part of the universal $R$-matrix has the form
\begin{align*}
R_{\mbf{m}}^+=1\otimes 1+\sum_{\alpha}e_{-[i;i+1)_{\mbf{m}}}\otimes e_{[i;i+1)_{\mbf{m}}}+\cdots.
\end{align*}
The operator $e_{[i;i+1)_{\mbf{m}}}$ maps $K_T(M(\mbf{v}_1,\mbf{w}_1))$ to $K_T(M(\mbf{v}_1-[i;i+1)_{\mbf{m}},\mbf{w}_1))$. As the common denominator of the components of $\mbf{m}$ increases, the positive components of $[i;i+1)_{\mbf{m}}$ eventually exceed the corresponding components of $\mbf{v}_1$. For such $\mbf{m}$, the reduced part of $R_{\mbf{m}}$ acts trivially.

Now we choose a point $\mbf{m}\in\mu^{-1}([t_1,t_2])$ such that for $\mbf{n}\in\mbb{N}^I$ less than or equal to $\mbf{v}_1$ we have $\mbf{m}\cdot\mbf{n}\in\mbb{Z}$. Now we consider the set of hyperplanes defined by $H_{\mbf{n},k}=\{\mbf{x}\in\mbb{R}^I|\mbf{n}\cdot\mbf{x}=k\in\mbb{Z}\}$. Now by definition, we can see that the intersection of $H_{\mbf{n},k}$ and the curve $\mu^{-1}([t_1,t_2])$ contains only finitely many points. Moreover, it is easy to see that for the fixed $\mbf{n}\in\mbb{N}^I$, only finitely many $k\in\mbb{Z}$ satisfy the property that $H_{\mbf{n},k}$ has nonempty intersection with $\mu^{-1}([t_1,t_2])$. Therefore combining the above properties, we have finished the proof.

\end{proof}

This result motivates the following definition of the \textbf{wall set} of $M(\mbf{v},\mbf{w})$:
\begin{align}\label{wall-set}
\text{Walls}(M(\mbf{v},\mbf{w}))=\{\mbf{m}\in\mbb{R}^{n}|\text{Reduced part of }R_{\mbf{m}}\text{ acts on } K_{T}(M(\mbf{v},\mbf{w}))\otimes K_{T}(M(\mbf{v},\mbf{w}))\text{non-trivially}\}.
\end{align}
By the preceding proposition, $\text{Walls}(M(\mbf{v},\mbf{w}))\cap\mu^{-1}([t_1,t_2])$ is finite.

For our situation, we choose the following coproduct:
\begin{align*}
\Delta_{(\mbf{s})}(a)=\prod^{\rightarrow}_{r\in\mbb{Q}_{\geq0}}(R_{\mbf{s}+\mu(r)}^+)^{-1}\Delta_{\infty}(a)\prod^{\leftarrow}_{r\in\mbb{Q}_{\geq0}}R_{\mbf{s}+\mu(r)}^+,\qquad\mu(0)=0.
\end{align*}
Here $\Delta_{\infty}:=R_{\infty}\Delta^{op}R_{\infty}^{-1}$. The corresponding universal $R$-matrix for $\Delta_{(\mbf{s})}$ is
\begin{align*}
R^{\mbf{s}}=\prod^{\leftarrow}_{r\in\mbb{Q}_{<0}}R_{\mbf{s}+\mu(r)}^{-}\cdot R_{\infty}\cdot\prod^{\leftarrow}_{r\in\mbb{Q}_{\geq0}}R_{\mbf{s}+\mu(r)}^+.
\end{align*}
Here $R_{\mbf{m}}^{-}:=(R_{\mbf{m}})_{21}$. Given two slopes $s$ and $s'$, and let $\Gamma$ be the path from $s$ to $s'$ the strictly increasing monotone curve $\mu$. Define
\begin{align*}
T^+_{s,s'}=\prod^{\leftarrow}_{\mbf{m}\in\Gamma}R_{\mbf{m}}^{+},\qquad T^-_{s,s'}=\prod^{\rightarrow}_{\mbf{m}\in\Gamma}(R_{\mbf{m}}^{+})^{-1}.
\end{align*}

Comparing the two formulas gives the following lemma.

\begin{lem}\label{coproduct-difference:lemma}
The coproduct $\Delta_{(\mbf{s})}$ and $\Delta_{(\mbf{s}')}$ are related by the following formula
\begin{align*}
\Delta_{(\mbf{s})}(a)=T_{s,s'}^{-}\Delta_{(\mbf{s}')}(a)T_{s,s'}^{+}.
\end{align*}

\end{lem}

The lemma has the following useful consequence.
\begin{cor}
Let $\mc{L}\in\text{Pic}(M(\mbf{v},\mbf{w}))$, then we have:
\begin{align*}
\Delta_{\infty}(\mc{L})=R_{\mbf{m}_{m-1}}^{+}\cdots R_{\mbf{m}_0}^{+}\Delta_{s}(\mc{L})
\end{align*}
in $\text{Hom}(K(\mbf{v},\mbf{w})\otimes K(\mbf{v},\mbf{w}),K(\mbf{w})\otimes K(\mbf{w}))$. Here $\Delta_{\infty}$ is the standard coproduct on $U_{q,t}(\hat{\hat{\mf{sl}}}_{n})$. $\mbf{m}_0,\cdots,\mbf{m}_{m-1}$ are the slope points between $s$ and $s+\mc{L}$.
\end{cor}
\begin{proof}
Note that $\mc{L}$ acts on $K_{T}(M(\mbf{v},\mbf{w}))$ as the Cartan element of $U_{q,t}(\hat{\hat{\mf{sl}}}_{n})$, and thus we have
\begin{align*}
\Delta_{\infty}(\mc{L})=\mc{L}\otimes\mc{L}.
\end{align*}

Applying the preceding formula gives
\begin{align*}
\Delta_{s}(\mc{L})=(R_{\mbf{m}_0}^+)^{-1}\cdots (R_{\infty}^{+})^{-1}\Delta_{\infty}(\mc{L})(R_{\infty}^{+}\cdots R_{\mbf{m}_0}^{+}).
\end{align*}

From Proposition \ref{translationpic} we have
\begin{align*}
\Delta_{\infty}(\mc{L})R_{\mbf{m}_k}^+\Delta_{\infty}(\mc{L})^{-1}=R_{\mbf{m}_k+\mc{L}}^{+}=R_{\mbf{m}_{k+m}}^+.
\end{align*}

This proves the result.
\end{proof}


\subsection{Geometric action on $K_{T}(M(\mbf{w}))$}
We first fix the notation for the affine type $A$ Nakajima quiver varieties.

Given a finite quiver $Q=(I,E)$, consider the following quiver representation space
\begin{align*}
\text{Rep}(\mbf{v},\mbf{w}):=\bigoplus_{h\in E}\text{Hom}(V_{i(h)},V_{o(h)})\oplus\bigoplus_{i\in I}\text{Hom}(V_i,W_i).
\end{align*}

Here $\mbf{v}=(\dim V_1,\ldots,\dim V_I)$ is the dimension vector of the vertex spaces, and $\mbf{w}=(\dim W_1,\ldots,\dim W_I)$ is the dimension vector of the framing spaces.

Denote $G_{\mbf{v}}:=\prod_{i\in I}GL(V_i)$ and $G_{\mbf{w}}:=\prod_{i\in I}GL(W_i)$. The natural action of $G_{\mbf{v}}\times G_{\mbf{w}}$ on $\text{Rep}(\mbf{v},\mbf{w})$ induces a Hamiltonian action on $T^*\text{Rep}(\mbf{v},\mbf{w})$ with respect to the standard symplectic form $\omega$. Let
\begin{align*}
\mu:T^*\text{Rep}(\mbf{v},\mbf{w})\rightarrow\mf{g}_{\mbf{v}}^*.
\end{align*}
For the generic stability condition $\theta\in\mbb{Z}^I$, we define the \textbf{Nakajima quiver variety} as the GIT quotient:
\begin{align}\label{Quiver-variety}
M_{\theta,0}(\mbf{v},\mbf{w}):=\mu^{-1}(0)\sslash_{\theta}G_{\mbf{v}}.
\end{align}

The Nakajima variety also has a natural action of $\mbb{C}_{q}^{\times}$ which scales the cotangent fibre with a character $q^{-1}$.
We abbreviate the Nakajima quiver variety in \ref{Quiver-variety} by $M(\mbf{v},\mbf{w})$.

We focus on affine type $A$ quivers $Q=(I,E)$, for which the representation space is
\begin{align*}
\text{Rep}(\mbf{v},\mbf{w})=\bigoplus_{i\in\mbb{Z}/n\mbb{Z}}\text{Hom}(V_i,V_{i+1})\oplus\bigoplus_{i\in I}\text{Hom}(V_i,W_i).
\end{align*}

The stability condition is chosen as $\theta=(1,\cdots,1)$.

The action of $\mbb{C}^*_{q}\times\mbb{C}^*_t\times G_{\mbf{w}}$ on $M(\mbf{v},\mbf{w})$ is given by:
\begin{align*}
(q,t,U_i)\cdot(X_e,Y_e,A_i,B_i)_{e\in E,i\in I}=(qtX_e,qt^{-1}Y_e,qA_iU_i,qU_i^{-1}B_i).
\end{align*}

Let $T_{\mbf{w}}\subset G_{\mbf{w}}$ denote as the maximal torus of $G_{\mbf{w}}$. The fixed-point set of the torus action $\mbb{C}^*_q\times\mbb{C}^*_t\times T_{\mbf{w}}$ on $M_{\theta,0}(\mbf{v},\mbf{w})$ is indexed by the $|\mbf{w}|$-partitions $\bm{\lambda}=(\lambda_1,\cdots,\lambda_{\mbf{w}})$. If we choose the cocharacter $\sigma:\mbb{C}^*\rightarrow G_{\mbf{w}}$ such that $\mbf{w}=u_1\mbf{w}_1+\cdots+u_k\mbf{w}_k$, we have:
\begin{align*}
M(\mbf{v},\mbf{w})^{\sigma}=\bigsqcup_{\mbf{v}_1+\cdots+\mbf{v}_k=\mbf{v}}M(\mbf{v}_1,\mbf{w}_1)\times\cdots\times M(\mbf{v}_k,\mbf{w}_k).
\end{align*}

We denote:
\begin{align*}
K(\mbf{w}):=K_{T}(M(\mbf{w}))_{loc}:=\bigoplus_{\mbf{w}}K_{\mbb{C}_q^*\times\mbb{C}_t^*\times T_{\mbf{w}}}(M(\mbf{v},\mbf{w}))_{loc}
\end{align*}
as the localised equivariant $K$-theory for $M(\mbf{v},\mbf{w})$. 

The $\mbb{C}_q^*\times\mbb{C}_t^*\times T$-fixed points for $M(\mbf{v},\mbf{w})$ are isolated, and the set of such fixed points is bijective to the following set of multipartitions:
\begin{align*}
\bm{\lambda}=(\lambda_1,\cdots,\lambda_{|\mbf{w}|})
\end{align*}
such that $|\bm{\lambda}|=|\mbf{v}|$, and such that the set of the boxes
\begin{align}
\{\square\in\bm{\lambda}|c_{\square}\equiv i(\text{mod} n)\}
\end{align}
has size $v_{i}$. Here if $\square\in\lambda_{j}^{(k)}$ with $\lambda_{j}=(\lambda_{j}^{(k)})$ a partition in the multipartition $\bm{\lambda}$, the content function for $\square\in\lambda_{j}^{(k)}$ is defined as $c_{\square}=x-y+k$ with $(x,y)$ the coordinates of $\square$.

We denote $I_{\bm{\lambda}}$ as the skyscraper sheaf corresponding to the fixed point $\bm{\lambda}$. We define the fixed-point basis by:
\begin{align*}
|\bm{\lambda}\rangle=\frac{I_{\bm{\lambda}}}{[T_{\bm{\lambda}}M(\mbf{v},\mbf{w})]}
\end{align*}

It has a natural basis $|\bm{\lambda}\rangle$ indexed by the corresponding partitions $\bm{\lambda}$, and we call it the fixed-point basis of $K_T(M(\mbf{w}))_{loc}$. 

Since $K_T(M(\mbf{v},\mbf{w}))$ is generated by tautological classes and $K_{T}(pt)$, if we denote $x_{ij}$ with $1\leq i\leq n$ and $1\leq j\leq v_{i}$ as the Chern roots. Given $f\in K_{T}(M(\mbf{v},\mbf{w}))$ the colored-symmetric Laurent polynomial in Chern roots (or Laurent polynomial of tautological classes), we define $\overline{f}$ as:
\begin{align}\label{fixed-point-expression}
\overline{f}:=\sum_{\bm{\lambda}}f|_{\bm{\lambda}}|\bm{\lambda}\rangle
\end{align}

Here we briefly review the geometric action of $U_{q,t}(\hat{\hat{\mf{sl}}}_{n})$ on $K_{T}(M(\mbf{w}))$. The construction is based on Nakajima's simple correspondence. One can refer to \cite{Nak01}\cite{N15} for the detailed construction.

Consider a pair of vectors $(\mbf{v}_+,\mbf{v}_-)$ such that $\mbf{v}_{+}=\mbf{v}_{-}+\mbf{e}_{i}$. We denote $Z(\mbf{v}_+,\mbf{v}_-,\mbf{w})$ as the closed subscheme in $M(\mbf{v}_+,\mbf{w})\times M(\mbf{v}_-,\mbf{w})$
parametrising the following short exact sequence of quiver representations
\begin{equation*}
\begin{tikzcd}
0\arrow[r]&K_i\arrow[d,shift left]\arrow[r]&V^+_i\arrow[r]\arrow[d,shift left]&V^-_i\arrow[r]\arrow[d,shift left]&0\\
0\arrow[r]&0\arrow[u,shift left]\arrow[r]&W_i\arrow[r]\arrow[u,shift left]&W_i\arrow[r]\arrow[u,shift left]&0
\end{tikzcd}
\end{equation*}
with $K_i=\mbb{C}^{\delta_{i0}}$. The scheme $Z(\mbf{v}_+,\mbf{v}_{-},\mbf{w})$ is smooth with a tautological line bundle:
\begin{align*}
\mc{L}|_{V^+\rightarrow V^-}=\text{Ker}(V_{i}^+\rightarrow V_{i}^{-}).
\end{align*}
There are natural projection maps:
\begin{equation*}
\begin{tikzcd}
&Z(\mbf{v}_+,\mbf{v}_{-},\mbf{w})\arrow[ld,"\pi_{+}"]\arrow[rd,"\pi_{-}"]&\\
M(\mbf{v}_+,\mbf{w})&&M(\mbf{v}_-,\mbf{w})
\end{tikzcd}.
\end{equation*}

Using $(\pi_+,\pi_-)$, we construct the operator
\begin{align*}
e_{i,d}^{\pm}:K_{T}(M(\mbf{v}^{\mp},\mbf{w}))\rightarrow K_{T}(M(\mbf{v}^{\pm},\mbf{w})),\qquad e_{i,d}^{\pm}(\alpha)=\pi_{\pm*}(\text{Cone}(d\pi_{\pm})\mc{L}^d\pi_{\mp}^{*}(\alpha)).
\end{align*}

Here $\text{Cone}(d\pi_{\pm})$ is defined in \cite{N15}. Summing over all $\mbf{v}$ gives an operator $e_{i,d}^{\pm}:K_T(M(\mbf{w}))\rightarrow K_T(M(\mbf{w}))$. The operators $\varphi_{i,d}^{\pm}$ act by multiplication by tautological classes; explicitly,
\begin{align}
\varphi_{i}^{\pm}(z)=\overline{\frac{\zeta{(\frac{z}{X})}}{\zeta(\frac{X}{z})}}\prod^{u_{j}\equiv i}_{1\leq j\leq\mbf{w}}\frac{[\frac{u_j}{\hbar z}]}{[\frac{z}{\hbar u_{j}}]},\qquad\overline{\frac{\zeta{(\frac{z}{X})}}{\zeta(\frac{X}{z})}}:=\prod_{j\in I}\prod_{1\leq a\leq v_j}\overline{\frac{\zeta_{ij}(\frac{z}{x_{ja}})}{\zeta_{ji}(\frac{x_{ja}}{z})}}
\end{align}
and here the overline is defined in \eqref{fixed-point-expression}. In particular, the element $\varphi_{i,0}$ acts on $K_{T}(M(\mbf{v},\mbf{w}))$ as $q^{\alpha_{i}^T(\mbf{w}-C\mbf{v})}$.

The above construction gives the following well-known result (\cite{N15}\cite{Nak01}):
\begin{thm}
For all $\mbf{w}\in\mbb{N}^r$, the operator $e_{i,d}^{\pm}$ and $\varphi^{\pm}_{i,d}$ give rise to an action of $U_{q,t}(\hat{\hat{\mf{sl}}}_{n})$ on $K_{T}(M(\mbf{w}))$.
\end{thm}

In the shuffle-algebra presentation, the action has the following explicit formula of the quantum toroidal algebra $U_{q,t}(\hat{\hat{\mf{sl}}}_{n})$ on $K_{T}(M(\mbf{w}))$:

Given $F\in\mc{S}_{\mbf{k}}^+$, we have
\begin{align*}
\langle\bm{\lambda}|F|\bm{\mu}\rangle=F(\chi_{\bm{\lambda}\backslash\bm{\mu}})\prod_{\blacksquare\in\bm{\lambda}\backslash\bm{\mu}}[\prod_{\square\in\bm{\mu}}\zeta(\frac{\chi_{\blacksquare}}{\chi_{\square}})\prod_{i=1}^{\mbf{w}}[\frac{u_i}{q\chi_{\blacksquare}}]]
\end{align*}
where $\bm{\lambda}\backslash\bm{\mu}$ is the skewed partition given by a pair of partitions $\bm{\mu}\subset\bm{\lambda}$. 

Similarly, for $G\in\mc{S}_{-\mbf{k}}^{-}$, we have
\begin{align*}
\langle\bm{\mu}|G|\bm{\lambda}\rangle=G(\chi_{\bm{\lambda}\backslash\bm{\mu}})\prod_{\blacksquare\in\bm{\lambda}\backslash\bm{\mu}}[\prod_{\square\in\bm{\lambda}}\zeta(\frac{\chi_{\square}}{\chi_{\blacksquare}})\prod_{i=1}^{\mbf{w}}[\frac{\chi_{\blacksquare}}{qu_i}]]^{-1}.
\end{align*}

\subsection{Khoroshkin-Tolstoy construction of $R$-matrix}
The Khoroshkin-Tolstoy construction \cite{KT92} of the universal $R$-matrix is, roughly speaking, a direct generalisation of the method in \cite{KT91} by constructing the Cartan-Weyl basis with respect to the root system of the corresponding quantum group. The Cartan-Weyl basis can be constructed via the $q$-Weyl group method in \cite{KR90}.

We recommend the reference \cite{BGKNV10} as a good reference for an in-depth and comprehensive treatment of the universal $R$-matrix for the quantum affine algebra $U_{q}(\hat{\mf{g}})$.

In the following subsection, we use the notation of \cite{BGKNV10}.

The formula for the universal $R$-matrix of $U_{q}(\hat{\mf{g}})$ is given by:
\begin{align}\label{KT-facto}
R=R_{<\delta}R_{\cong\delta}R_{>\delta}K.
\end{align}

Here:
\begin{equation}\label{KT-facto2}
\begin{aligned}
&R_{<\delta}=\prod_{m\geq0}^{\rightarrow}\prod_{\gamma\in\Delta(A)}^{\leftarrow}\exp_{q^{-(\gamma,\gamma)}}((q-q^{-1})s_{m,\gamma}^{-1}e_{\gamma+m\delta}\otimes f_{\gamma+m\delta})\\
&R_{\cong\delta}=\exp((q-q^{-1})\sum_{m\in\mbb{Z}_+}\sum_{i,j=1}^{r}u_{m,ij}e_{m\delta,\alpha_{i}}\otimes f_{m\delta,\alpha_{j}})\\
&R_{>\delta}=\prod_{m\geq0}^{\rightarrow}\prod_{\gamma\in\Delta(A)}^{\rightarrow}\exp_{q^{-(\gamma,\gamma)}}((q-q^{-1})s_{m,\delta-\gamma}^{-1}e_{(\delta-\gamma)+m\delta}\otimes f_{(\delta-\gamma)+m\delta})
\end{aligned}
\end{equation}

Here $\Delta(A)$ is the set of positive roots for $\mf{g}$, and $\delta=\alpha_0+\cdots+\alpha_{n}$ is the imaginary root. $e_{\alpha}$ and $f_{\alpha}$ are the Cartan-Weyl basis constructed in \cite{KT92}. Here $K$ is the universal $R$-matrix induced by $\mbb{Q}(q,t)[\varphi_{i,0}^{\pm}]$ via identifying $\varphi_{i,0}=\exp(\hbar b_i)$

\subsection{$R$-matrix presentation in the geometric module}
In this subsection, we prove the Laurent polynomiality of the $R$-matrix in each module $K_{T\times\mbb{C}^{\times}_{q}}(M(\mbf{w}))$.

\begin{prop}\label{polynomiality-of-R-matrices}
Let $K_{T}(pt):=\mbb{C}[u_{i}^{\pm1}]_{i\in I}$, and let the equivariant parameter of $M(\mbf{w}_1)$ and $M(\mbf{w}_2)$ be $u_1$ and $u_2$. Then the matrix coefficients of the representation of $R_{\mbf{m}}\in\mc{B}_{\mbf{m}}\hat{\otimes}\mc{B}_{\mbf{m}}$ in $\text{End}(K_{T\times\mbb{C}^{\times}_{q}}(M(\mbf{w}_1))_{loc}\otimes K_{T\times\mbb{C}^{\times}_{q}}(M(\mbf{w}_2)))_{loc}$ are a Laurent monomial in $u_1/u_2$, with coefficients in $\mbb{Q}(q,t)$. 
\end{prop}
\begin{proof}
Now we fix two quiver varieties $M(\mbf{w}_1)$ and $M(\mbf{w}_2)$ and the equivariant parameter torus on $M(\mbf{w}_1)$ and $M(\mbf{w}_2)$ to be $u_1$ and $u_2$. For the given representation:
\begin{align*}
(\pi_1\otimes \pi_2):\mc{B}_{\mbf{m}}\hat{\otimes}\mc{B}_{\mbf{m}}\rightarrow \text{End}(K_{T\times\mbb{C}^{\times}_{q}}(M(\mbf{w}_1))_{loc}\otimes K_{T\times\mbb{C}^{\times}_{q}}(M(\mbf{w}_2))_{loc})
\end{align*}

Using the KT factorisation \ref{KT-facto}, we obtain the isomorphism of $\mc{B}_{\mbf{m}}$ in \ref{rootquantum}.

Now we fix a vector $v_1\otimes v_2\in K_{T\times\mbb{C}^{\times}_{q}}(M(\mbf{v}_1,\mbf{w}_1))_{loc}\otimes K_{T\times\mbb{C}^{\times}_{q}}(M(\mbf{v}_2,\mbf{w}_2))_{loc}$. Because the universal $R$-matrix $R_{\mbf{m}}$ lies in $\mc{B}_{\mbf{m}}^{\geq}\hat{\otimes}\mc{B}_{\mbf{m}}^{\leq}$, elements of $\mc{B}_{\mbf{m}}^{\leq}$ decrease the dimension vector of the factor $M(\mbf{v}_2,\mbf{w}_2)$. Since every graded piece $\mc{B}_{\mbf{m},-\mbf{k}}^{\leq}$ is finite-dimensional, only finitely many terms of $R_{\mbf{m}}$ act nontrivially on $v_1\otimes v_2$.

Consequently, to show that $(\pi_1\otimes\pi_2)(R_{\mbf{m}})$ is a Laurent polynomial in $u_1/u_2$, it suffices to show that $(\pi_{1}\otimes\pi_{2})(P^{\mbf{m}}_{[i;j)}\otimes Q^{\mbf{m}}_{-[i;j)})$ is a Laurent polynomial in $u_1/u_2$. This follows from the fact that $P^{\mbf{m}}_{[i;j)}|\bm{\lambda}_1\rangle$ is a monomial in $u_1$ of degree $\mbf{m}\cdot[i;j)$ and $Q^{\mbf{m}}_{-[i;j)}|\bm{\lambda}_2\rangle$ is a monomial in $u_2$ of degree $-\mbf{m}\cdot[i;j)$.

\end{proof}

Using this result, we can prove the rationality of the operator $\prod_{\mu\in\mbb{Q}_{>0}}R_{\mbf{m}+\mu\bm{\theta}}^{+}(u)$ valued in $\text{End}(K(\mbf{w}_1)\otimes K(\mbf{w}_2))$ with $u=u_1/u_2$:
\begin{prop}\label{rationality-of-R-matrix}
The matrix coefficients of the operator $\prod_{\mu\in\mbb{Q}_{>0}}R_{\mbf{m}+\mu\bm{\theta}}^{+}(u)$ valued in $\text{End}(K(\mbf{w}_1)\otimes K(\mbf{w}_2))$ are rational functions of $u$. 
\end{prop}
\begin{proof}
We will prove this proposition in Appendix \ref{appendix}.
\end{proof}

\section{\textbf{Quantum difference equation}}

\subsection{Background: Stable quasimaps to Nakajima varieties}

A quasimap $f:C\rightarrow M(\mbf{v},\mbf{w})$ with $C$ the irreducible component isomorphic to $\mbb{P}^1$ consists of the following data:
\begin{itemize}
	\item A collection of vector bundles $\mc{V}_{i}$, $i\in I$, on $C$, with $\operatorname{rk}(\mc{V}_i)=v_i$. These bundles are induced by the corresponding bundles $V_i$ on $M(\mbf{v},\mbf{w})$.
	\item A collection of trivial vector bundles $\mc{W}_{i}$ of rank $w_{i}$.

	\item A section
	\begin{align}
	f\in H^0(C,\mc{M}\oplus\mc{M}^*\otimes q^{-1})
	\end{align}
	satisfying the moment-map condition $\mu=0$, where
	\begin{align}
	\mc{M}=(\bigoplus_{h\in E}\text{Hom}(\mc{V}_{i(h)},\mc{V}_{o(h)})\oplus\bigoplus_{i\in I}\text{Hom}(\mc{V}_i,\mc{W}_i)).
	\end{align}
\end{itemize}
The degree of a quasimap is $d=(\text{deg}(\mc{V}_{i}))_{i=1}^I\in\text{Pic}(M(\mbf{v},\mbf{w}))\cong\mbb{Z}^{I}$.

Let $QM^d(X)$ denote the moduli space of stable quasimaps of degree $d$, and let $QM^d(X)_{\text{nonsing }p}\subset QM^d(X)$ be the open subspace parametrising quasimaps that are nonsingular at $p$. This subspace carries a natural evaluation morphism $ev_p:QM^d(X)_{\text{nonsing }p}\rightarrow X$. The moduli space of \textbf{relative quasimaps}, denoted $QM^d(X)_{\text{rel }p}$, compactifies this evaluation morphism:
\begin{equation*}
\begin{tikzcd}
&QM^d_{\text{rel }p}(X)\arrow[dr,"\hat{ev}_{p}"]&\\
QM^d_{\text{nonsing }p}(X)\arrow[ur]\arrow[rr]&&X
\end{tikzcd}.
\end{equation*}

This moduli space is equipped with the natural symmetrised virtual structure sheaves $\hat{\mc{O}}_{vir}$. 

Consider the moduli space $QM^{d}_{\text{rel }p_1,\text{nonsing }p_2}(X)$ of quasimaps with relative conditions at $p_1\in C$ and nonsingular at $p_2\in C$. Thus we can define the evaluation map:
\begin{align*}
ev=\hat{ev}_{p_1}\times ev_{p_2}:QM^{d}_{\text{rel }p_1,\text{nonsing }p_2}(X)\rightarrow X\times X.
\end{align*}

This moduli space is equipped with an action of $G\times\mbb{C}_{p}^{\times}$ with $G=\prod_{i}GL_{w_i}$ coming from $X$ and $\mbb{C}_{q}^{\times}$ coming from $C$. 

Note that this map $ev$ is not proper, but it becomes proper on the subset of fixed points $QM^{d}_{\text{rel }p_1,\text{nonsing }p_2}(X)^{\mbb{C}_{p}^{*}}$.

Thus we can define the \textbf{Capping operator}:
\begin{align}\label{capping-difference}
\bm{\Psi}(z)=\sum_{d\in\mbb{Z}^n}z^d\text{ev}_*(QM^{d}_{\text{rel }p_1,\text{nonsing }p_2}(X),\hat{\mc{O}}_{vir})\in K_{G\times\mbb{C}_{p}^{\times}}(X)^{\otimes 2}_{loc}\otimes\mbb{Q}[[z]].
\end{align}

Here $z\in\text{Pic}(X)\otimes\mbb{C}^{*}$ called the \textbf{K\"ahler parameter}. It is proved in \cite{AO21} that the capping operator is the matrix of the fundamental solution of a system of $q$-difference equations:
\begin{align*}
\mbf{J}(zp^{\mc{L}})\mc{L}=\mbf{M}_{\mc{L}}(z)\mbf{J}(z).
\end{align*}

In this paper, we focus on the operator $\mbf{M}_{\mc{L}}(z)$, which is called the \textbf{quantum difference operator}.

\subsection{Quantum difference equation in the algebraic setting}
In this paper, we construct the algebraic quantum difference equation using the quantum toroidal algebra. This construction may be viewed as an algebraic counterpart of the Okounkov--Smirnov quantum difference equation for affine type $A$ Nakajima quiver varieties.

The construction parallels that of \cite{OS22} and may be summarised as follows.

Denote the translation operator with respect to $\mc{L}\in\mbb{Z}^n$ by $T_{\mc{L}}$:
\begin{align*}
T_{\mc{L}}f(z)=f(zp^{\mc{L}}).
\end{align*}

We define the difference operator:
\begin{align*}
A_{\mc{L}}=T_{\mc{L}}^{-1}\mbf{M}_{\mc{L}}(z).
\end{align*}
These operators satisfy the commutativity relation
\begin{align*}
[A_{\mc{L}},A_{\mc{L}'}]=0.
\end{align*}

In this section, we prove the algebraic analogue of this result.

We now denote $n:=\#I$, and we consider the K\"ahler parameter $z=(z_i)_{i\in I}\in\text{Pic}(X)\otimes\mbb{C}^{\times}\cong(\mbb{C}^{\times})^{n}$.

First we define the operator $z_{(k)}\in\text{End}(K_{T}(F))$ with $F=M(\mbf{v}_1,\mbf{w}_1)\times\cdots\times M(\mbf{v}_k,\mbf{w}_k)\subset X^A$ to be diagonal in the basis supported on fixed-point set:
\begin{align}\label{kahler-parametre-multiplication-operator}
z_{(k)}(\gamma)=z_1^{v_{k,1}}\cdots z_{n}^{v_{k,n}}\gamma,\qquad\gamma\in K_T(F).
\end{align}

We define $q^{\Omega}\in\text{End}(K_{G}(F))$ as:
\begin{align*}
q^{\Omega}(\gamma)=q^{\text{codim}(F)}(\gamma),\qquad\text{codim}(F)=\text{dim}(M(\sum_{i}\mbf{v}_i,\sum_{i}\mbf{w}_i))-\sum_{i}\text{dim}(M(\mbf{v}_i,\mbf{w}_i)).
\end{align*}

To each slope subalgebra $\mc{B}_{\mbf{m}}$, we associate fusion operators $J_{\mbf{m}}^{\pm}(z)$ as follows.

\begin{prop}\label{ABRR-equation}
There exist unique solutions $J_{\mbf{m}}^{+}(z)$ and $J_{\mbf{m}}^{-}(z)$ of the following ABRR equations, with $J_{\mbf{m}}^{+}(z)$ strictly upper triangular and $J_{\mbf{m}}^{-}(z)$ strictly lower triangular:
\begin{align}\label{ABRR-eqn}
J_{\mbf{m}}^{+}(z)z_{(1)}^{-1}q^{\Omega}R_{\mbf{m}}^{+}=z_{(1)}^{-1}q^{\Omega}J_{\mbf{m}}^{+}(z),\qquad q^{\Omega}R_{\mbf{m}}^{-}z_{(1)}^{-1}J_{\mbf{m}}^{-}(z)=J_{\mbf{m}}^{-}(z)q^{\Omega}z_{(1)}^{-1}
\end{align}
where $R_{\mbf{m}}^{-}:=(R_{\mbf{m}}^+)_{21}$ is the transposition of $R_{\mbf{m}}^+$. Moreover, $J_{\mbf{m}}^{\pm}(z)$ are elements in a completion $\mc{B}_{\mbf{m}}\hat{\otimes}\mc{B}_{\mbf{m}}$ satisfying:
\begin{align*}
S_{\mbf{m}}\otimes S_{\mbf{m}}((J_{\mbf{m}}^{+}(z))_{21})=J_{\mbf{m}}^{-}(z).
\end{align*}
\end{prop}
\begin{proof}
For details, see Proposition $7$ of \cite{OS22}. The argument can be applied since we can replace the wall subalgebra $U_q(\mf{g}_w)$ by the slope subalgebra $\mc{B}_{\mbf{m}}$, and replace the Maulik-Okounkov quantum affine algebra $U_{q,t}^{MO}(\hat{\mf{g}}_Q)$ by the quantum toroidal algebras $U_{q,t}(\hat{\hat{\mf{sl}}}_n)$. 
\end{proof}

The ABRR equation originates in \cite{ABRR97}, where it is used to construct solutions of the quantum dynamical Yang--Baxter equation.

An immediate corollary is the following.
\begin{lem}
\begin{align*}
\lim_{z\rightarrow\infty}J_{\mbf{m}}^+(z^{\theta})=1,\qquad\lim_{z\rightarrow0}J_{\mbf{m}}^+(z^{\theta})=R_{\mbf{m}}^{+}
\end{align*}
where $\theta$ is the parameter appearing in the root factorisation of the quantum toroidal algebra, and $\theta\in(\mbb{Q}^{+})^n$.
\end{lem}
\begin{proof}
The argument is essentially the same as that in \cite{OS22}; for completeness, we give a direct computation.

Let us decompose $J_{\mbf{m}}^{+}$ and $R_{\mbf{m}}^{+}$ according to the shuffle grading of the slope subalgebra:
\begin{align*}
J_{\mbf{m}}^+(z)=1+\sum_{\mbf{n}>0}J_{\mbf{m}|\mbf{n}}(z),\qquad (R_{\mbf{m}}^+)^{-1}=1+\sum_{\mbf{n}>0}R_{\mbf{m}|\mbf{n}}^{+}.
\end{align*}

The ABRR equation gives the recursion
\begin{align}\label{recursive-formula-for-fusion-operator}
J_{\mbf{m}|\mbf{n}}(z)=\frac{1}{z^{\mbf{n}}q^{k}-1}\sum_{\substack{\mbf{n}_1+\mbf{n}_2=\mbf{n}\\\mbf{n}_1<\mbf{n}}}J_{\mbf{m}|\mbf{n}_1}(z)R_{\mbf{m}|\mbf{n}_2}^{+}.
\end{align}

Set $(z_1,\ldots,z_n)=(z^{\theta_1},\ldots,z^{\theta_n})$. As $z\rightarrow\infty$, the recursion implies that $J_{\mbf{m}|\mbf{n}}\rightarrow0$. 

As $z\rightarrow0$, we obtain
\begin{align*}
J_{\mbf{m}|\mbf{n}}(0^{\theta})=\sum_{k}\sum_{\mbf{n}_1+\cdots+\mbf{n}_k=\mbf{n}}(-1)^{k}R_{\mbf{m}|\mbf{n}_1}^{+}R_{\mbf{m}|\mbf{n}_2}^{+}\cdots R_{\mbf{m}|\mbf{n}_k}^{+}.
\end{align*}

Thus we have:
\begin{equation*}
\begin{aligned}
J_{\mbf{m}}(0^{\theta})=&1+\sum_{k}\sum_{\mbf{n}_1,\cdots,\mbf{n}_k}(-1)^kR_{\mbf{m}|\mbf{n}_1}^{+}R_{\mbf{m}|\mbf{n}_2}^{+}\cdots R_{\mbf{m}|\mbf{n}_k}^{+}\\
=&(1+\sum_{\mbf{n}>0}R_{\mbf{m}|\mbf{n}}^{+})^{-1}\\
=&R_{\mbf{m}}^{+}.
\end{aligned}
\end{equation*}

\end{proof}

Let us define the shift-automorphism map $T_u:\mc{B}_{\mbf{m}}\rightarrow\mc{B}_{\mbf{m}}$ by
\begin{align*}
T_u(F)=u^{\pm\mbf{m}\cdot\mbf{k}}F,\qquad F\in\mc{B}_{\mbf{m},\pm\mbf{k}}^{\pm}
\end{align*}
and in this way, we define $R_{\mbf{m}}^{\pm}(u):=(D_{u}\otimes\text{Id})R_{\mbf{m}}^{\pm}$. Let $\tau_{\mbf{m}}=\mbf{s}\mc{L}_{\mbf{m}}$ with $q^s=p$. Then:
\begin{align*}
q_{(1)}^{\tau_{\mbf{m}}}R^{+}_{\mbf{m}}(u)q_{(1)}^{-\tau_{\mbf{m}}}=R_{\mbf{m}}^{+}(up).
\end{align*}

By the previous proposition we have
\begin{align*}
q_{(1)}^{\tau_{\mbf{m}}}J_{\mbf{m}}^{+}(u,z)q_{(1)}^{-\tau_{\mbf{m}}}=J_{\mbf{m}}^{+}(up,z).
\end{align*}
Thus we can rewrite the above equation as:
\begin{align*}
J_{\mbf{m}}^{+}(u,zp^{-\mbf{m}})z_{(1)}^{-1}q^{\Omega}R_{\mbf{m}}^{+}(uq)=z_{(1)}^{-1}q^{\Omega}J_{\mbf{m}}^{+}(up,zp^{-\mbf{m}})
\end{align*}

We define the fusion operator
\begin{align*}
\mbf{J}_{\mbf{m}}^{\pm}(z)=J_{\mbf{m}}^{\pm}(zp^{-\mbf{m}}).
\end{align*}

The preceding relation yields the following proposition.
\begin{prop}\label{wall-KZ-eqn}
There exist unique operators $\mbf{J}_{\mbf{m}}^{+}(z)\in\mc{B}_{\mbf{m}}^{-}\hat{\otimes}\mc{B}_{\mbf{m}}^{+}(z)$ and $\mbf{J}_{\mbf{m}}^{-}(z)\in\mc{B}_{\mbf{m}}^{+}\hat{\otimes}\mc{B}_{\mbf{m}}^{-}(z)$ satisfying the following shifted ABRR equations, where $\mbf{J}_{\mbf{m}}^{+}(z)$ is strictly upper triangular and $\mbf{J}_{\mbf{m}}^{-}(z)$ is strictly lower triangular:
\begin{align*}
&\mbf{J}_{\mbf{m}}^{+}(z)z_{(1)}^{-1}T_{u}q^{\Omega}R_{\mbf{m}}^{+}=z_{(1)}^{-1}T_{u}q^{\Omega}\mbf{J}_{\mbf{m}}^{+}(z)\\
&R_{\mbf{m}}^{-}z_{(1)}^{-1}T_{u}q^{\Omega}\mbf{J}_{\mbf{m}}^{-}(z)=\mbf{J}_{\mbf{m}}^{-}(z)z_{(1)}^{-1}T_{u}q^{\Omega}.
\end{align*}
\end{prop}
In \cite{OS22}, these equations are called the \textbf{wall Knizhnik--Zamolodchikov equations}.

We now define the monodromy operator by
\begin{align}\label{defn-of-quantum-difference-operator}
\mbf{B}_{\mbf{m}}(z)=m((1\otimes S_{\mbf{m}})(\mbf{J}_{\mbf{m}}^{-}(z)^{-1}))|_{z\mapsto zq^{\kappa}}.
\end{align}

Here $\kappa=\frac{C\mbf{v}-\mbf{w}}{2}$, and $m$ denotes the multiplication map. We call this dynamical operator the \textbf{monodromy operator}.\footnote{This terminology reflects the fact that, as shown in \cite{OS22,S21}, these dynamical operators represent the $q$-monodromy of the quantum difference equation. We will elaborate on this interpretation in the sequel.}

If $\mbf{m}$ is not a wall for $M(\mbf{v},\mbf{w})$, $R_{\mbf{m}}^{\pm}$ gives the trivial action. Proposition \ref{wall-KZ-eqn} then implies that $J_{\mbf{m}}^{\pm}$ is also trivial. It therefore suffices to define $\mbf{B}_{\mbf{m}}(z)$ for points $\mbf{m}$ in the wall set, where $R_{\mbf{m}}^{\pm}$ acts nontrivially.

Let $\mc{L}\in\operatorname{Pic}(M(\mbf{v},\mbf{w}))$ be a line bundle. Fix a slope $s\in H^2(M(\mbf{v},\mbf{w}),\mbb{R})$ and choose a path from $s$ to $s-\mc{L}$. By Proposition \ref{finiteness-of-wall-set}, the path crosses only finitely many slope points, say $\mbf{m}_1,\ldots,\mbf{m}_m$, in that order. To this choice of slope, line bundle, and path, we associate the following operator on $K_T(M(\mbf{v},\mbf{w}))$:
\begin{align}\label{defnofqdeoperator}
\mbf{B}_{\mc{L}}^{s}(z):=\mc{L}\prod^{\leftarrow}_{\mbf{m}\in(s,s-\mc{L}]}\mbf{B}_{\mbf{m}}(z)=\mc{L}\mbf{B}_{\mbf{m}_m}(z)\cdots\mbf{B}_{\mbf{m}_1}(z).
\end{align}

We define the $p$-difference operators:
\begin{align*}
\mc{A}^{s}_{\mc{L}}=T_{\mc{L}}^{-1}\mbf{B}^s_{\mc{L}}(z).
\end{align*}

The algebraic quantum difference equation for the affine type $A$ Nakajima quiver variety is defined as:
\begin{align*}
\Psi(zp^{\mc{L}})=\mbf{B}^s_{\mc{L}}(z)\Psi(z)\in\text{End}(K_{T}(M(\mbf{v},\mbf{w}))_{loc})
\end{align*}
or equivalently
\begin{align*}
\mc{A}^s_{\mc{L}}\Psi(z)=\Psi(z).
\end{align*}

\subsection{Monodromy operator for the slope subalgebra $\mc{B}_{\mbf{m}}$}
In this subsection, we compute the monodromy operator $\mbf{B}_{\mbf{m}}(z)\in\mc{B}_{\mbf{m}}$ associated with the slope subalgebra $\mc{B}_{\mbf{m}}$, for $\mbf{m}\in\mbb{Q}^n$.

By \cite{N15}, the slope subalgebra $\mc{B}_{\mbf{m}}$ for the quantum toroidal algebra $U_{q,t}(\hat{\hat{\mf{sl}}}_{n})$ is isomorphic to:
\begin{align*}
\mc{B}_{\mbf{m}}\cong\bigotimes^{g}_{h=1}U_{q}(\widehat{\mf{gl}}_{l_h}).
\end{align*}

The isomorphism is constructed as follows. It is proved in \cite{N15} that, as a shuffle algebra, $\mc{B}_{\mbf{m}}$ is the Drinfeld double of $\mc{B}^{\pm}_{\mbf{m}}$:
\begin{align*}
\mc{B}_{\mbf{m}}:=\mc{B}_{\mbf{m}}^{+}\otimes\mbb{C}[\varphi_{[i;j)}^{\pm1}]\otimes\mc{B}_{\mbf{m}}^{-}/(\text{relations}).
\end{align*}

The generators for the slope subalgebra can be written as follows: Recall the elements in \ref{positive-generators} and \ref{negative generators}, for $\mbf{m}\cdot[i;j)\in\mbb{Z}$:
\begin{equation}\label{Pformula}
P_{\pm[i ; j)}^{\pm \mbf{m}}=\operatorname{Sym}\left[\frac{\prod_{a=i}^{j-1} z_a^{\left\lfloor m_i+\ldots+m_a\right\rfloor-\left\lfloor m_i+\ldots+m_{a-1}\right\rfloor}}{t^{\text {ind } d_{[i, j)}^m q^{i-j}} \prod_{a=i+1}^{j-1}\left(1-\frac{q_2 z_a}{z_{a-1}}\right)} \prod_{i \leq a<b<j} \zeta\left(\frac{z_b}{z_a}\right)\right]\in\mc{B}_{\mbf{m}}^{\pm}
\end{equation}
\begin{equation}\label{Qformula}
Q_{\mp[i ; j)}^{\pm \mbf{m}}=\operatorname{Sym}\left[\frac{\prod_{a=i}^{j-1} z_a^{\left\lfloor m_i+\ldots+m_{a-1}\right\rfloor-\left\lfloor m_i+\ldots+m_a\right\rfloor}}{t^{-\mathrm{ind}_{[i, j)}^m} \prod_{a=i+1}^{j-1}\left(1-\frac{q_1 z_{a-1}}{z_a}\right)} \prod_{i \leq a<b<j} \zeta\left(\frac{z_a}{z_b}\right)\right]\in\mc{B}_{\mbf{m}}^{\pm}
\end{equation}
and here:
\begin{align*}
\text{ind}^{\mbf{m}}_{[i;j)}=\sum_{a=i}^{j-1}(m_i+\cdots+m_a-\lfloor m_i+\cdots+m_{a-1}\lfloor).
\end{align*}

A direct calculation shows that the antipode map $S_{\mbf{m}}:\mc{B}_{\mbf{m}}\rightarrow\mc{B}_{\mbf{m}}$ satisfies the relations
\begin{align*}
S_{\mbf{m}}(P^{\mbf{m}}_{[i;j)})=Q^{\mbf{m}}_{[i;j)},\qquad S_{\mbf{m}}(Q^{\mbf{m}}_{-[i;j)})=P^{\mbf{m}}_{-[i;j)}.
\end{align*}

Now the slope subalgebra $\mc{B}_{\mbf{m}}$ is generated by $Q_{\mp[i ; j)}^{\pm \mbf{m}}$ and $P_{\pm[i ; j)}^{\pm \mbf{m}}$ and $\varphi_{k}$. The isomorphism \ref{rootquantum} is given by:
\begin{align*}
e_{[i;j)}=P^{\mbf{m}}_{[i;j)_{h}},\qquad e_{-[i;j)}=Q^{\mbf{m}}_{-[i;j)_h},\qquad\varphi_{k}=\varphi_{[k;v_{\mbf{m}}(k))}.
\end{align*}

Under this isomorphism, the generators act by the morphisms
\begin{align*}
&e_{[i;j)}:K_{T}(M(\mbf{v},\mbf{w}))\rightarrow K_{T}(M(\mbf{v}+[i;j)_{h},\mbf{w}))\\
&e_{-[i;j)}:K_{T}(M(\mbf{v},\mbf{w}))\rightarrow K_{T}(M(\mbf{v}-[i;j)_{h},\mbf{w})).
\end{align*}

In particular:
\begin{align*}
&e_{[i;i+1)}:K_{T}(M(\mbf{v},\mbf{w}))\rightarrow K_{T}(M(\mbf{v}+[i;i+1)_{h},\mbf{w}))\\
&e_{-[i;i+1)}:K_{T}(M(\mbf{v},\mbf{w}))\rightarrow K_{T}(M(\mbf{v}-[i;i+1)_{h},\mbf{w})).
\end{align*}

To compute the monodromy operator $\mbf{B}_{\mbf{m}}(z)$, we need to solve the ABRR equation \ref{ABRR-eqn}. Now we compute the formal solution for the monodromy operator $\mc{B}_{\mbf{m}}$:
\begin{align*}
&J_{\mbf{m}}^{-}(z)=\prod_{k=0}^{\substack{\rightarrow\\\infty}}\text{Ad}_{(q^{\Omega}z_{(1)}^{-1})^k}\text{Ad}_{q^{\Omega}}R_{\mbf{m}}^{-},\qquad R_{\mbf{m}}^{-}=(R_{\mbf{m}})_{21}\\
&\mbf{B}_{\mbf{m}}(z)=m((1\otimes S_{\mbf{m}})(J_{\mbf{m}}^{-}(zp^{-\mbf{m}}))|_{z\rightarrow zq^{\kappa}}.
\end{align*}
and here $S_{\mbf{m}}$ is the antipode map on the slope subalgebra $\mc{B}_{\mbf{m}}$, $m$ stands for the multiplication map $m:\mc{B}_{\mbf{m}}\otimes\mc{B}_{\mbf{m}}\rightarrow\mc{B}_{\mbf{m}}$.

By the isomorphism \ref{quantum gln decomposition} the universal $R$-matrix is written as:
\begin{align*}
R_{\mbf{m}}^-=\prod_{h=1}^{g}(R_{\hat{\mf{sl}}_{l_h}}^-R_{\hat{\mf{gl}}_1}^-).
\end{align*}

Using the KT factorisation of $R_{\mbf{m}}^-$, the computation reduces to the following factors:
\begin{align}\label{A1}
\exp_{q^{-2}}((q-q^{-1})e_{\gamma+m\delta}\otimes f_{\gamma+m\delta})
\end{align}
\begin{align}\label{A2}
\exp((q-q^{-1})\sum_{m\in\mbb{Z}_{+}}\sum_{i,j=1}^{l_h}u_{m,ij}e_{m\delta,\alpha_i}\otimes f_{m\delta,\alpha_i})
\end{align}
\begin{align}\label{A3}
\exp_{q^{-2}}((q-q^{-1})e_{(\delta-\gamma)+m\delta}\otimes f_{(\delta-\gamma)+m\delta})
\end{align}
\begin{align}\label{A4}
R_{\mf{gl}_1}^-=\exp(\sum_{k=1}^{\infty}n_kp_k\otimes p_{-k}).
\end{align}
and here $u_{m,ij}$ are the matrix elements of the matrix $u_{m}$ inverse to the matrix $t_m$ with the matrix elements
\begin{align}
t_{m,ij}=(-1)^{m(1-\delta_{ij})}m^{-1}\frac{q^{ma_{ij}}-q^{-ma_{ij}}}{q-q^{-1}}
\end{align}
and $(a_{ij})$ is the corresponding Cartan matrix for $\hat{\mf{sl}}_{l_h}$.

Now we choose specific $U_{q}(\hat{\mf{gl}}_{l_h})$, which corresponds to the cycle $(i_1,i_2,\cdots,i_{l_h})$ as defined in \eqref{oriented-cycles}. Set $\delta=\sum_{k=1}^{l_h}[i_k;i_{k+1})$. A direct computation shows that:

\begin{equation*}
\begin{aligned}
&\text{Ad}_{(q^{\Omega}z_{(1)}^{-1})^k}\exp_{q^{-2}}((q-q^{-1})e_{\gamma+m\delta}\otimes f_{\gamma+m\delta})\\
=&\exp_{q^{-2}}((q-q^{-1})z^{-k(\mbf{v}_{\gamma}+m\bm{\delta})}q^{2k(\mbf{v}_{\gamma}+m\bm{\delta})^TC(\mbf{v}_{\gamma}+m\bm{\delta})}(\varphi^{k(\mbf{v}_{\gamma}+m\bm{\delta})}\otimes\varphi^{-k(\mbf{v}_{\gamma}+m\bm{\delta})})e_{\gamma+m\delta}\otimes f_{\gamma+m\delta})
\end{aligned}
\end{equation*}

\begin{equation*}
\begin{aligned}
&\text{Ad}_{(q^{\Omega}z_{(1)}^{-1})^k}\exp((q-q^{-1})\sum_{m\in\mbb{Z}_{+}}\sum_{i,j=1}^{l_h}u_{m,ij}e_{m\delta,\alpha_i}\otimes f_{m\delta,\alpha_i})\\
=&\exp((q-q^{-1})\sum_{m\in\mbb{Z}_{+}}\sum_{i,j=1}^{l_h}z^{-km\bm{\delta}}q^{2km^2\bm{\delta}^TC\bm{\delta}}(\varphi^{km\bm{\delta}}\otimes\varphi^{-km\bm{\delta}})u_{m,ij}e_{m\delta,\alpha_i}\otimes f_{m\delta,\alpha_i})
\end{aligned}
\end{equation*}

\begin{equation*}
\begin{aligned}
&\text{Ad}_{(q^{\Omega}z_{(1)}^{-1})^k}\exp_{q^{-2}}((q-q^{-1})f_{(\delta-\gamma)+m\delta}\otimes e_{(\delta-\gamma)+m\delta})\\
=&\exp_{q^{-2}}((q-q^{-1})z^{-k(-\mbf{v}_{\gamma}+(m+1)\bm{\delta})}q^{2k(-\mbf{v}_{\gamma}+(m+1)\bm{\delta})^TC(-\mbf{v}_{\gamma}+(m+1)\bm{\delta})}\\
&(\varphi^{k(-\mbf{v}_{\gamma}+(m+1)\bm{\delta})}\otimes\varphi^{-k(-\mbf{v}_{\gamma}+(m+1)\bm{\delta})})e_{(\delta-\gamma)+m\delta}\otimes f_{(\delta-\gamma)+m\delta})
\end{aligned}
\end{equation*}

\begin{equation*}
\begin{aligned}
&\text{Ad}_{(q^{\Omega}z_{(1)}^{-1})^k}\exp(\sum_{k=1}^{\infty}n_kp_k\otimes p_{-k})\\
=&\exp(\sum_{k=1}^{\infty}z^{-km\bm{\delta}}(\varphi^{km\bm{\delta}}\otimes\varphi^{-km\bm{\delta}})n_kp_k\otimes p_{-k}).
\end{aligned}
\end{equation*}

Here the operator $\varphi^{\bm{\alpha}}\in\text{End}(K_{T}(M(\mbf{v},\mbf{w}))$ acts on $K_{T}(M(\mbf{v},\mbf{w})$ as:
\begin{align*}
\varphi^{\bm{\alpha}}v=q^{2\alpha^T(\mbf{w}-C\mbf{v})}v.
\end{align*}

The shift $z\rightarrow zp^{-\mbf{m}}$ has the indicated effect, while $z\rightarrow zq^{\kappa}$ becomes $z\rightarrow zq^{\frac{1}{2}(\mbf{v}^TC\mbf{v}-\mbf{w}^T\mbf{w})}$. The antipode maps the $(e_{\pm[i;j)})$-basis to the $(f_{\mp[i;j)})$-basis.

Combining these computations yields the following theorem.
\begin{thm}\label{universal-formula-for-difference-operator}
On the representation space $\text{End}(K_T(M(\mbf{v},\mbf{w})))$, the monodromy operator $\mbf{B}_{\mbf{m}}(z)$ defined in \ref{defn-of-quantum-difference-operator} has the following formal expansion:
\begin{equation*}
\begin{aligned}
&\mbf{B}_{\mbf{m}}(z)\\
=&\mu(\prod_{h=1}^{g}(\textbf{Heisenberg algebra part})\prod_{k=0}^{\substack{\rightarrow\\\infty}}\\
&\prod_{\substack{\gamma\in\Delta(A)\\m\geq0}}^{\leftarrow}(\exp_{q^{2}}(-(q-q^{-1})z^{k(-\mbf{v}_{\gamma}+(m+1)\bm{\delta}_h)}p^{-k\mbf{m}\cdot(-\mbf{v}_{\gamma}+(m+1)\bm{\delta}_h)}q^{-\frac{k}{2}(\mbf{v}^TC\mbf{v}-\mbf{w}^T\mbf{w})\bm{\theta}\cdot(-\mbf{v}_{\gamma}+(m+1)\bm{\delta}_h))}(\\
&q^{2k(-\mbf{v}_{\gamma}+(m+1)\bm{\delta})^TC(-\mbf{v}_{\gamma}+(m+1)\bm{\delta})}\varphi^{(k+1)(-\mbf{v}_{\gamma}+(m+1)\bm{\delta})}e_{(\delta-\gamma)+m\delta}\otimes f_{(\delta-\gamma)+m\delta}'\varphi^{-(k+1)(-\mbf{v}_{\gamma}+(m+1)\bm{\delta})})\\
&\times\exp(-(q-q^{-1})\sum_{m\in\mbb{Z}_{+}}\sum_{i,j=1}^{l_h}z^{km\bm{\delta}_h}p^{-km\mbf{m}\cdot\bm{\delta}_{h}}u_{m,ij}q^{2km^2\bm{\delta}^TC\bm{\delta}}\varphi^{(k+1)m\bm{\delta}}e_{m\delta,\alpha_i}\otimes f_{m\delta,\alpha_i}'\varphi^{-(k+1)m\bm{\delta}})\\
&\times\prod_{\substack{\gamma\in\Delta(A)\\m\geq0}}^{\rightarrow}\exp_{q^{2}}(-(q-q^{-1})z^{k(\mbf{v}_{\gamma}+m\bm{\delta}_h)}p^{-k\mbf{m}\cdot(\mbf{v}_{\gamma}+m\bm{\delta}_h)}q^{-\frac{k}{2}(\mbf{v}^TC\mbf{v}-\mbf{w}^T\mbf{w})\mbf{\theta}\cdot(\mbf{v}_{\gamma}+m\bm{\delta}_h))}\\
&q^{2k(\mbf{v}_{\gamma}+m\bm{\delta})^TC(\mbf{v}_{\gamma}+m\bm{\delta})}(\varphi^{(k+1)(\mbf{v}_{\gamma}+m\bm{\delta})}e_{\gamma+m\delta}\otimes f_{\gamma+m\delta}'\varphi^{-(k+1)(\mbf{v}_{\gamma}+m\bm{\delta})}))).
\end{aligned}
\end{equation*}

Here $\mu$ is the multiplication map. $e_v'=S(e_v)$ and $f_v'=S(f_v)$ are the images under the antipode map. $\bm{\theta}=(1,\cdots,1)$. Here $e_{\alpha}$ and $f_{\alpha}$ stand for the generators in $\mc{B}_{\mbf{m}}$ without using $\mbf{m}$. The \textbf{Heisenberg algebra part} denotes a term of the form \ref{instanton-moduli-difference-operator}; it can be constructed recursively from the following formula:
\begin{align*}
S(P_{\pm[i;j)}^{\pm\mbf{m}})=Q^{\pm\mbf{m}}_{\pm[i;j)},\qquad S(Q_{\pm[i;j)}^{\pm\mbf{m}})=P^{\pm\mbf{m}}_{\pm[i;j)}
\end{align*}
with $P^{\pm\mbf{m}}_{\pm[i;j)}$ and $Q^{\pm\mbf{m}}_{\pm[i;j)}$ defined in \ref{Pformula} and \ref{Qformula}.

\end{thm}

The following statement gives a useful convergence property of the monodromy operator $\mbf{B}_{\mbf{m}}(z)$. Although it may be useful in analysing convergence of solutions to the difference equation, it will not be used in this paper.
\begin{prop}
The monodromy operator $\mbf{B}_{\mbf{m}}(z)$ is convergent in each $K_{T}(M(\mbf{v},\mbf{w}))$ of $\mc{B}_{\mbf{m}}(z)$ for $\{z_i\}$ sufficiently small.
\end{prop}
\begin{proof}
Since the fusion operator is the rational function over $z$ by the formula \eqref{recursive-formula-for-fusion-operator}, the definition of the monodromy operator in \eqref{defn-of-quantum-difference-operator} implies  that the above expansion is the expansion of a rational function over $z$. Hence it is convergent for small $z$.
\end{proof}

\subsection{Some properties of $\mbf{B}_{\mbf{m}}(z)$ and $J_{\mbf{m}}^{\pm}(z)$}

Here we show some essential properties of monodromy operators $\mbf{B}_{\mbf{m}}(z)$ and fusion operators $J_{\mbf{m}}^{\pm}(z)$. They are similar to the properties of the monodromy operators and fusion operators defined in \cite{OS22}. It is noteworthy that we have successfully adapted and replicated similar formulations within the language of the quantum toroidal algebra by requiring some modification in the construction.

We state the following theorem without proof in our current framework. For the proof, see Theorem 5 and Theorem 6 in \cite{OS22}, where the same principles apply within the quantum toroidal algebra setting.

\begin{thm}
The operators $J^{\pm}_{\mbf{m}}(z)$ satisfy the dynamical cocycle conditions on $\text{End}(K_{G}(M(\mbf{v}_1,\mbf{w}_1)\times M(\mbf{v}_2,\mbf{w}_2)\times M(\mbf{v}_3,\mbf{w}_3)))$:
\begin{equation*}
\begin{aligned}
&J^{-}_{\mbf{m}}(z)^{12,3}J^{-}_{\mbf{m}}(zq^{\kappa_{(3)}})^{12}=J^{-}_{\mbf{m}}(z)^{1,23}J^{-}_{\mbf{m}}(zq^{-\kappa_{(1)}})^{23}\\
&J^{+}_{\mbf{m}}(zq^{\kappa_{(3)}})^{12}J^{+}_{\mbf{m}}(z)^{12,3}=J^{+}_{\mbf{m}}(zq^{-\kappa_{(1)}})^{23}J^{+}_{\mbf{m}}(z)^{1,23}
\end{aligned}
\end{equation*}
in $\mc{B}_{\mbf{m}}\hat{\otimes}\mc{B}_{\mbf{m}}\hat{\otimes}\mc{B}_{\mbf{m}}$. Here $X^{12,3}:=(\Delta_{\mbf{m}}\otimes\text{Id})X, X^{12}=X\otimes\text{Id}, X^{1,23}=(\text{Id}\otimes\Delta_{\mbf{m}})X$.
$\kappa(\mbf{v},\mbf{w})=(C\mbf{v}-\mbf{w})/2$ is the dynamical shift operator, which is the product of $\varphi_{i,0}$ in $U_{q,t}(\hat{\hat{\mf{sl}}}_{n})$. Moreover, we set
\begin{align*}
\tilde{B}_{\mbf{m}}(z)=\mbf{m}(1\otimes S_{\mbf{m}}(J_{\mbf{m}}^{-}(z)^{-1}))|_{z\mapsto zq^{\kappa}},\qquad B_{\mbf{m}}(z)=\mbf{m}_{21}(S_{\mbf{m}}^{-1}\otimes 1(J_{\mbf{m}}^{-}(z)^{-1}))|_{z\mapsto zq^{-\kappa}}.
\end{align*}

We have the following coproduct formula:
\begin{align*}
\Delta_{\mbf{m}}\tilde{B}_{\mbf{m}}(z)=J_{\mbf{m}}^{-}(z)(\tilde{B}_{\mbf{m}}(zq^{\kappa_{(2)}})\otimes\tilde{B}_{\mbf{m}}(zq^{-\kappa_{(1)}}))J_{\mbf{m}}^{+}(z)
\end{align*}
\begin{align}\label{coproduct-for-monodromy}
\Delta_{\mbf{m}} B_{\mbf{m}}(z)=J_{\mbf{m}}^{-}(z)(B_{\mbf{m}}(zq^{\kappa_{(2)}})\otimes B_{\mbf{m}}(zq^{-\kappa_{(1)}}))J_{\mbf{m}}^{+}(z).
\end{align}
\end{thm}

Now with the shift $z\mapsto zp^{-\mbf{m}}$, we have the coproduct formula for $\mbf{B}_{\mbf{m}}(z)$:
\begin{align}\label{coproduct-for-monodromy-operator}
\Delta_{\mbf{m}}(\mbf{B}_{\mbf{m}}(z))=\mbf{J}_{\mbf{m}}^{-}(z)(\mbf{B}_{\mbf{m}}(zq^{\kappa_{(2)}})\otimes\mbf{B}_{\mbf{m}}(zq^{-\kappa_{(1)}}))\mbf{J}_{\mbf{m}}^{+}(z).
\end{align}

Now given a slope subalgebra $\mc{B}_{\mbf{m}}\subset U_{q,t}(\hat{\hat{\mf{sl}}}_{n})$ and its corresponding universal $R$-matrix $R_{\mbf{m}}^+$ and $R_{\mbf{m}}^-$, we set the qKZ operator for each slope subalgebra $\mc{B}_{\mbf{m}}$:
\begin{align*}
\mc{R}_{\mbf{m}}^+=T_{u}z_{(1)}^{-1}q^{\Omega}R_{\mbf{m}}^+,\qquad\mc{R}_{\mbf{m}}^{-}=R_{\mbf{m}}^{-}q^{\Omega}T_{u}z_{(1)}^{-1}.
\end{align*}

It can be checked that it commutes with the monodromy operator:
\begin{prop}
\begin{align*}
\mc{R}_{\mbf{m}}^{-}\Delta_{\mbf{m}}(\mbf{B}_{\mbf{m}}(z))=\Delta_{\mbf{m}}(\mbf{B}_{\mbf{m}}(z))\mc{R}_{\mbf{m}}^{+}.
\end{align*}
\end{prop}
\begin{proof}
From the ABRR equation and the above coproduct formula, we have the formula, for details see Proposition $11$ in \cite{OS22}.
\end{proof}

\begin{prop}\label{translationpic}
For $\mc{L}\in\text{Pic}(X)$, where $X=M(\mbf{v},\mbf{w})$ is an affine type $A$ quiver variety, we have the following translation formula:
\begin{align}\label{translationformula}
\mc{L}\mbf{B}_{\mbf{m}}(zp^{-\mc{L}})=\mbf{B}_{\mbf{m}+\mc{L}}(z)\mc{L}\in\text{End}(K_{T}(M(\mbf{v},\mbf{w}))_{loc}).
\end{align}
\end{prop}
\begin{proof}

For two slope points $\mbf{m}$ and $\mbf{m}'$, the corresponding slope subalgebra $\mc{B}_{\mbf{m}}$ and $\mc{B}_{\mbf{m}'}$ are isomorphic if $\mbf{m}'-\mbf{m}\in\mbb{Z}^{|I|}$.

Let $X=M(\mbf{v},\mbf{w})$ and $A\subset T$ such that $X^A=\bigsqcup_{\mbf{v}_1+\mbf{v}_2=\mbf{v}}M(\mbf{v}_1,\mbf{w}_1)\times M(\mbf{v}_2,\mbf{w}_2)$. 

Let $\mc{L}\in\text{Pic}(X)\cong\mbb{Z}^n$. Comparing the actions of $\mc{L}\mbf{B}_{\mbf{m}}(zp^{-\mc{L}})\mc{L}^{-1}$ and $\mbf{B}_{\mbf{m}+\mc{L}}(z)$ on $K_T(X)$ shows that $\mc{L}\mbf{B}_{\mbf{m}}(zp^{-\mc{L}})\mc{L}^{-1}\in\mc{B}_{\mbf{m}+\mc{L}}$. It remains to prove that $U(\mc{L})J_{\mbf{m}}^{-}(zp^{-\mc{L}})U(\mc{L})^{-1}$ equals $J_{\mbf{m}+\mc{L}}^{-}(z)$. Both operators solve the same ABRR equation for $R_{\mbf{m}+\mc{L}}^{-}$ in Proposition \ref{ABRR-eqn}. Here $U(\mc{L})\in\text{End}(K_T(X^A))$ is the block-diagonal operator defined by
\begin{align*}
U(\mc{L})|_{F_i\times F_i}=\mc{L}|_{F_i}.
\end{align*}

We need to show that for $R_{\mbf{m}}$ acting on $K_{T}(X^A)$, $R_{\mbf{m}}$ has the following translation formula:
\begin{align}\label{translation-of-wall-R-matrix}
U(\mc{L})R_{\mbf{m}}^{\pm}U(\mc{L})^{-1}=R_{\mbf{m}+\mc{L}}^{\pm}.
\end{align}

Using the expression of the KT factorisation of $R_{\mbf{m}}$, it is equivalent to prove the following identity:
\begin{align*}
\mc{L}P^{\mbf{m}}_{\pm[i;j)}\mc{L}^{-1}=P^{\mbf{m}+\mc{L}}_{\pm[i;j)},\qquad \mc{L}Q^{\mbf{m}}_{\pm[i;j)}\mc{L}^{-1}=Q^{\mbf{m}+\mc{L}}_{\pm[i;j)}.
\end{align*}

We give the proof only for $P^{\mbf{m}}_{\pm[i;j)}$. The case of $Q^{\mbf{m}}_{\pm[i;j)}$ is analogous.
We examine this via the matrix element of $\mc{L}P^{\mbf{m}}_{[i;j)}\mc{L}^{-1}$ in $K(\mbf{w})$. Recall that the matrix element of $P^{\mbf{m}}_{[i;j)}$ is given by \cite{N15}:
\begin{equation}\label{Pieri-for-P}
\langle\bm{\lambda}|P^{\mbf{m}}_{[i;j)}|\bm{\mu}\rangle=P^{\mbf{m}}_{[i;j)}(\bm{\lambda}\backslash\bm{\mu})\prod_{\blacksquare\in\bm{\lambda}\backslash\bm{\mu}}[\prod_{\square\in\bm{\mu}}\zeta(\frac{\chi_{\blacksquare}}{\chi_{\square}})\prod_{k=1}^{\mbf{w}}[\frac{u_k}{q\chi_{\blacksquare}}]].
\end{equation}

Here $\bm{\lambda}\backslash\bm{\mu}$ is the skew partition with respect to the pair $\bm{\mu}\subset\bm{\lambda}$. $P^{\mbf{m}}_{[i;j)}(\bm{\lambda}\backslash\bm{\mu})$ is given by:
\begin{equation}\label{Pieri-for-P-2}
P_{[i ; j)}^{\mathbf{m}}(\boldsymbol{\lambda} \backslash \boldsymbol{\mu})=\sum_{\mathrm{ASYT}^{+}} \prod_{a=i}^{j-1} \chi_a^{m_a} \cdot \frac{\prod_{a=i}^{j-1}\left(\chi_a t^{-a}\right)^{s_a}}{q^{-1} \prod_{a=i+1}^{j-1}\left(\frac{1}{q}-\frac{\chi_a}{t_{\chi_{a-1}}}\right)} \prod_{i \leq a<b<j} \zeta\left(\frac{\chi_b}{\chi_a}\right)
\end{equation}
for $\bm{\mu}\subset\bm{\lambda}$ and $0$ otherwise. We also use the notation
\begin{equation*}
s_a=\left\lfloor m_i+\ldots+m_a\right\rfloor-\left\lfloor m_i+\ldots+m_{a-1}\right\rfloor-m_a
\end{equation*}
and here $\text{ASYT}^{+}$ is the set of bijections $\psi:\{\bm{\lambda}\backslash\bm{\mu}\}\rightarrow\{i,\cdots,j-1\}$ such that $\psi(\square)=\text{color of }\square$.

Now since $\mc{L}$ is the product of the tautological bundle, for simplicity, we assume that $\mc{L}=\mc{L}_{l}$ is the generator corresponding to the vector space $V_{l}$ in the quiver. In this case we have:
\begin{equation*}
\langle\bm{\lambda}|\mc{L}P^{\mbf{m}}_{\pm[i;j)}\mc{L}^{-1}|\bm{\mu}\rangle=\prod_{\square\in\bm{\lambda}\backslash\bm{\mu}}^{c_{\square}=l}\chi_{\square}P^{\mbf{m}}_{[i;j)}(\bm{\lambda}\backslash\bm{\mu})\prod_{\blacksquare\in\bm{\lambda}\backslash\bm{\mu}}[\prod_{\square\in\bm{\mu}}\zeta(\frac{\chi_{\blacksquare}}{\chi_{\square}})\prod_{k=1}^{\mbf{w}}[\frac{u_k}{q\chi_{\blacksquare}}]].
\end{equation*}

Since the expression after $P^{\mbf{m}}_{[i;j)}(\bm{\lambda}\backslash\bm{\mu})$ only depends on $\bm{\lambda}$ and $\bm{\mu}$, we only need to compute $\prod_{\square\in\bm{\lambda}\backslash\bm{\mu}}\chi_{\square}P^{\mbf{m}}_{[i;j)}(\bm{\lambda}\backslash\bm{\mu})$. Note that:
\begin{equation*}
\begin{aligned}
&\prod_{\square\in\bm{\lambda}\backslash\bm{\mu}}\chi_{\square}P^{\mbf{m}}_{[i;j)}(\bm{\lambda}\backslash\bm{\mu})\\
=&\prod_{\square\in\bm{\lambda}\backslash\bm{\mu}}^{c_{\square}=l}\chi_{\square}\sum_{\mathrm{ASYT}^{+}} \prod_{a=i}^{j-1} \chi_a^{m_a} \cdot \frac{\prod_{a=i}^{j-1}\left(\chi_a t^{-a}\right)^{s_a}}{q^{-1} \prod_{a=i+1}^{j-1}\left(\frac{1}{q}-\frac{\chi_a}{t_{\chi_{a-1}}}\right)} \prod_{i \leq a<b<j} \zeta\left(\frac{\chi_b}{\chi_a}\right).
\end{aligned}
\end{equation*}

We can observe that in the expression, $\chi_{\square}$ contains all the boxes in $\bm{\lambda}\backslash\bm{\mu}$ such that they are of colour $l$, while the expression $\prod_{a=i}^{j-1}\chi_{a}^{m_{a}}$ contains the product $\prod_{i\leq a\leq j-1}^{c_{\psi^{-1}(a)}=l}\chi_{\psi^{-1}(a)}^{m_a}=\prod_{\square\in\bm{\lambda}\backslash\bm{\mu}}^{c_{\square}=l}\chi_{\square}^{m_{\psi(\square)}}$, thus from this expression we have:
\begin{equation*}
\begin{aligned}
&\prod_{\square\in\bm{\lambda}\backslash\bm{\mu}}^{c_{\square}=l}\chi_{\square}\sum_{\mathrm{ASYT}^{+}} \prod_{a=i}^{j-1} \chi_a^{m_a} \cdot \frac{\prod_{a=i}^{j-1}\left(\chi_a t^{-a}\right)^{s_a}}{q^{-1} \prod_{a=i+1}^{j-1}\left(\frac{1}{q}-\frac{\chi_a}{t_{\chi_{a-1}}}\right)} \prod_{i \leq a<b<j} \zeta\left(\frac{\chi_b}{\chi_a}\right)\\
=&\sum_{\mathrm{ASYT}^{+}} \prod_{a=i}^{j-1} \chi_a^{m_a+\delta_{ai}} \cdot \frac{\prod_{a=i}^{j-1}\left(\chi_a t^{-a}\right)^{s_a}}{q^{-1} \prod_{a=i+1}^{j-1}\left(\frac{1}{q}-\frac{\chi_a}{t_{\chi_{a-1}}}\right)} \prod_{i \leq a<b<j} \zeta\left(\frac{\chi_b}{\chi_a}\right).
\end{aligned}
\end{equation*}

Note also that applying the integral translation $\mbf{m}\mapsto\mbf{m}+\mc{L}_{i}$ does not change $s_{a}$, thus we have proved that:
\begin{align*}
\mc{L}P^{\mbf{m}}_{\pm[i;j)}\mc{L}^{-1}=P^{\mbf{m}+\mc{L}}_{\pm[i;j)}.
\end{align*}

The same argument applies to $Q^{\mbf{m}}_{\pm[i;j)}$.

In conclusion, we see that $U(\mc{L})J_{\mbf{m}}^{-}(z)U(\mc{L})^{-1}$ and $J_{\mbf{m}+\mc{L}}^{-}(z)$ solve the same ABRR equation. By the uniqueness of the solution of the ABRR equation, we have:
\begin{align*}
U(\mc{L})J_{\mbf{m}}^{-}(z)U(\mc{L})^{-1}=J_{\mbf{m}+\mc{L}}^{-}(z).
\end{align*}

This completes the proof.
\end{proof}

\subsection{Main result for the quantum difference operator $\mbf{B}^{s}_{\mc{L}}(z)$}
In this subsection, we prove that the quantum difference operator $\mbf{B}^{s}_{\mc{L}}(z)$ is generically independent of the choice of slope.

Now we fix a universal $R$-matrix $R^s$ for $U_{q,t}(\hat{\hat{\mf{sl}}}_{n})$ with respect to the factorisation
\begin{align*}
U_{q,t}(\hat{\hat{\mf{sl}}}_{n})=\bigotimes_{m\in\mbb{Q}}^{\rightarrow}\mc{B}_{s+m\bm{\theta}}^{+}\otimes\mc{B}_{\infty\bm{\theta}}\otimes\bigotimes_{m\in\mbb{Q}}^{\rightarrow}\mc{B}_{s+m\bm{\theta}}^{-}.
\end{align*}

Let $\mc{R}^s(u)$ be the image of $R^s$ acting on $K_{T}(M(\mbf{w}'))\otimes K_{T}(M(\mbf{w}''))$. Define the qKZ operator with the slope $s$ by
\begin{align*}
\mc{A}_s^{qKZ}:=z_{(1)}T_{u}^{-1}\mc{R}^s(u).
\end{align*}

We define the quantum difference connection as:
$$\mc{A}_{\mc{L}}^s=T_{\mc{L}}^{-1}\mbf{B}_{\mc{L}}^{s}(z)$$.

\begin{prop}\label{coproduct-for-difference}
On $K_{T}(M(\mbf{v}_1,\mbf{w}_1))\otimes K_{T}(M(\mbf{v}_2,\mbf{w}_2))$ we have
\begin{align*}
\Delta_{s}(\mc{A}^{s}_{\mc{L}})=W_{\mbf{m}_0}(z)W_{\mbf{m}_1}(z)\cdots W_{\mbf{m}_{m-1}}(z)\Delta_{\infty}(\mc{L})T_{\mc{L}}^{-1}
\end{align*}
with $\mbf{m}_0,\cdots,\mbf{m}_{m-1}$ the ordered set of slope points on a strictly monotone path from $s$ to $s+\mc{L}$ . $W_{\mbf{m}}(z)=\Delta_{\mbf{m}}(\mbf{B}_{\mbf{m}}(z))(R_{\mbf{m}}^{+})^{-1}$ and $\Delta_{\infty}$ is the original coproduct \eqref{Drinfeld-coproduct} of $U_{q,t}(\hat{\hat{\mf{sl}}}_{n})$.
\end{prop}
\begin{proof}
By definition we know that:
\begin{align*}
\mc{A}^s_{\mc{L}}=T_{\mc{L}}^{-1}\mc{L}\mbf{B}_{\mbf{m}_{-m}}(z)\cdots\mbf{B}_{\mbf{m}_{-2}}(z)\mbf{B}_{\mbf{m}_{-1}}(z).
\end{align*}

Here $\mbf{m}_{-1},\cdots,\mbf{m}_{-m}$ are the ordered set of slope points between $s$ and $s-\mc{L}$. By the translation symmetry in Proposition \ref{translationpic} we know that:
\begin{align*}
\mc{A}^s_{\mc{L}}=\mbf{B}_{\mbf{m}_0}(z)\mbf{B}_{\mbf{m}_1}(z)\cdots\mbf{B}_{\mbf{m}_{m-1}}(z)\mc{L}T_{\mc{L}}^{-1}.
\end{align*}
and here $\mbf{m}_{k+m}=\mbf{m}_{k}+\mc{L}$.

Now use the coproduct formula:
\begin{align*}
\Delta_{s}(\mc{A}^s_{\mc{L}})=\Delta_{s}(\mbf{B}_{\mbf{m}_0}(z)\mbf{B}_{\mbf{m}_1}(z)\cdots\mbf{B}_{\mbf{m}_{m-1}}(z)\mc{L})T_{\mc{L}}^{-1}.
\end{align*}

Note that the coproducts $\Delta_{s}$ and $\Delta_{\mbf{m}_k}$ differ as follows by Lemma \ref{coproduct-difference:lemma}:
\begin{align*}
\Delta_{s}(\mbf{B}_{\mbf{m}_k}(z))=(T_{s,\mbf{m}_{k}}^{-})\Delta_{\mbf{m}_{k}}(\mbf{B}_{\mbf{m}_k}(z))T_{s,\mbf{m}_k}^{+}.
\end{align*}

Writing this relation explicitly gives
\begin{align*}
\Delta_{s}(\mbf{B}_{\mbf{m}_k}(z))=(R_{\mbf{m}_0}^{+})^{-1}\cdots (R_{\mbf{m}_{k-1}}^{+})^{-1}\Delta_{\mbf{m}_k}(\mbf{B}_{\mbf{m}_k}(z))R_{\mbf{m}_{k-1}}^{+}\cdots R_{\mbf{m}_0}^{+}.
\end{align*}

Thus we have:
\begin{align*}
\Delta_{s}(\mc{A}^s_{\mc{L}})=&\Delta_{\mbf{m}_0}(\mbf{B}_{\mbf{m}_0}(z))(R_{\mbf{m}_0}^{+})^{-1}\cdots\Delta_{\mbf{m}_{m-1}}(\mbf{B}_{\mbf{m}_{m-1}}(z))(R_{\mbf{m}_{m-1}}^+)^{-1}R_{\mbf{m}_{m-1}}^{+}\cdots R_{\mbf{m}_0}^{+}\Delta_{s}(\mc{L})T_{\mc{L}}^{-1}\\
=&\Delta_{\mbf{m}_0}(\mbf{B}_{\mbf{m}_0}(z))(R_{\mbf{m}_0}^{+})^{-1}\cdots\Delta_{\mbf{m}_{m-1}}(\mbf{B}_{\mbf{m}_{m-1}}(z))(R_{\mbf{m}_{m-1}}^+)^{-1}\Delta_{\infty}(\mc{L})T_{\mc{L}}^{-1}.
\end{align*}
This completes the proof.
\end{proof}

\begin{thm}
On the space $\text{End}(K_{T}(M(\mbf{w}_1))\otimes K_{T}(M(\mbf{w}_2)))$, for two slopes $s$ and $s'$ separated by a single slope point $\mbf{m}$, we have:
\begin{align*}
W_{\mbf{m}}(z)^{-1}\mc{A}_s^{qKZ}W_{\mbf{m}}(z)=\mc{A}_{s'}^{qKZ},\qquad W_{\mbf{m}}(z)^{-1}\Delta_{s}(\mc{A}^s_{\mc{L}})W_{\mbf{m}}(z)=\Delta_{s'}(\mc{A}^{s'}_{\mc{L}}).
\end{align*}
\end{thm}
\begin{proof}
By the computation:
\begin{equation*}
\begin{aligned}
\mc{A}_s^{qKZ}W_{\mbf{m}}(z)=&z_{(1)}T_{u}^{-1}\mc{R}^s(u)\Delta_{\mbf{m}}(\mbf{B}_{\mbf{m}}(z))(R_{\mbf{m}}^+)^{-1}\\
=&z_{(1)}T_{u}^{-1}\Delta_{\mbf{m}}^{op}(\mbf{B}_{\mbf{m}}(z))\mc{R}^s(u)(R_{\mbf{m}}^+)^{-1}\\
=&z_{(1)}T_{u}^{-1}(R_{\mbf{m}}^{-})^{-1}(R_{\mbf{m}}^{-})\Delta_{\mbf{m}}^{op}(\mbf{B}_{\mbf{m}}(z))\mc{R}^{s}(u)(R_{\mbf{m}}^+)^{-1}\\
=&z_{(1)}T_{u}^{-1}(R_{\mbf{m}}^{-})^{-1}\Delta_{\mbf{m}}(\mbf{B}_{\mbf{m}}(z))(R_{\mbf{m}}^{-})\mc{R}^s(u)(R_{\mbf{m}}^{+})^{-1}\\
=&\Delta_{\mbf{m}}(\mbf{B}_{\mbf{m}}(z))(R_{\mbf{m}}^{+})^{-1}z_{(1)}T_{u}^{-1}\mc{R}^{s'}(u)=W_{\mbf{m}}(z)\mc{A}_{s'}^{qKZ}.
\end{aligned}
\end{equation*}

Now for the second equality, the proof is analogous and is omitted.
\end{proof}

\begin{thm}\label{qKZ-commute}
For arbitrary line bundles $\mc{L}$, $\mc{L}'\in\text{Pic}(M(\mbf{v},\mbf{w}))$ and a slope $s$, the qKZ operators commute with the $q$-difference operators
\begin{align*}
\Delta_{s}(\mc{A}^s_{\mc{L}})\mc{A}_s^{qKZ}=\mc{A}_s^{qKZ}\Delta_{s}(\mc{A}^s_{\mc{L}})
\end{align*}
\end{thm}
\begin{proof}
\begin{equation*}
\begin{aligned}
\mc{A}_s^{qKZ}\Delta_s(\mc{A}^s_{\mc{L}})=&\mc{A}_s^{qKZ}W_{\mbf{m}_0}(z)W_{\mbf{m}_1}(z)\cdots W_{\mbf{m}_{m-1}}(z)\Delta_{\infty}(\mc{L})T_{\mc{L}}^{-1}\\
=&W_{\mbf{m}_0}(z)W_{\mbf{m}_1}(z)\cdots W_{\mbf{m}_{m-1}}(z)\mc{A}_{s+\mc{L}}^{qKZ}\Delta_{\infty}(\mc{L})T_{\mc{L}}^{-1}\\
=&\Delta_{s}(\mc{A}^{s}_{\mc{L}})\mc{A}_s^{qKZ}.
\end{aligned}
\end{equation*}
\end{proof}

We now prove the path independence of the quantum difference operator.
\begin{thm}\label{path-independence-theorem}
The quantum difference operator $\mbf{B}^s_{\mc{L}}(z)$ is independent of the choice of the path from $s$ to $s-\mc{L}$ provided that the path $\mu:\mbb{R}\rightarrow\mbb{R}^n$ satisfies $\mu(t_1)>\mu(t_2)$ for $t_1>t_2$.
\end{thm}
\begin{proof}
Let $\mbf{B}^s_{\mc{L}}(z)'$ and $\mbf{B}^s_{\mc{L}}(z)$ be the operators associated with two different paths, acting on $K_{T}(M(\mbf{v},\mbf{w}))_{loc}$. Let $u$ denote the relevant torus-equivariant parameter, and set
\begin{align*}
D^s_{\mc{L}}(z)=\mbf{B}^s_{\mc{L}}(z)^{-1}\mbf{B}^s_{\mc{L}}(z)'\in\text{End}(K_{T}(M(\mbf{v},\mbf{w}))_{loc}).
\end{align*}

Since $\mc{A}^{s}_{\mc{L}}$ commutes with $\mc{A}_s^{qKZ}$, we have $[\Delta_{s}(D^s_{\mc{L}})(u,z),\mc{A}_s^{qKZ}]=0$, and this means that:
\begin{align}\label{conjugate-equation}
\Delta_{s}(D^s_{\mc{L}})(up,z)=z_{(1)}\mc{R}^s(u)\Delta_{s}(D^s_{\mc{L}})(u,z)(z_{(1)}\mc{R}^s(u))^{-1}=z_{(1)}\Delta_s^{op}(D^s_{\mc{L}})(u,z)z_{(1)}^{-1}
\end{align}
and here the second equality comes from the fact that $\mc{R}^{s}(u)\Delta_s(D^s_{\mc{L}})(u,z)(\mc{R}^s(u))^{-1}=\Delta_s^{op}(D^s_{\mc{L}})(u,z)$. Doing this equation again and we obtain the equation:
\begin{equation*}
\Delta_s(D^s_{\mc{L}})(up^2,z)=\Delta_s(D^s_{\mc{L}})(u,z)
\end{equation*}
and since $\Delta_s(D^s_{\mc{L}})(u,z)$ is a Laurent polynomial in $u$ by Proposition \ref{coproduct-for-difference} and formula \eqref{coproduct-for-monodromy-operator}, we have that $\Delta_s(D^s_{\mc{L}})(u,z)=\Delta_s(D^s_{\mc{L}}(z))$ is independent of $u$.

If we restrict the operator $\Delta_{s}(D^s_{\mc{L}})(z)$ to the lowest component $K_{T\times\mbb{C}^{\times}_{q}}(M(0,\mbf{w}_1))_{loc}\otimes K_{T\times\mbb{C}^{\times}_{q}}(M(\mbf{v}_2+\mbf{v}_1,\mbf{w}_2))_{loc}$, we have the following claim:

\textbf{Claim.}
For $s\in\mbb{R}^n$ with all components non-negative, when restricted to the component $\text{End}(K_{T\times\mbb{C}^{\times}_{q}}(M(0,\mbf{w}_1))_{loc}\otimes K_{T\times\mbb{C}^{\times}_{q}}(M(\mbf{v}_2+\mbf{v}_1,\mbf{w}_2))_{loc})$, the operator $\Delta_s(D^s_{\mc{L}})(z)$ can be written as:
\begin{align}\label{restricted-coproduct-D}
\Delta_{s}(D^s_{\mc{L}}(z))=D^s_{\mc{L}}(zq^{\kappa_{(2)}})\otimes D^s_{\mc{L}}(zq^{-\kappa_{(1)}})=1\otimes D^s_{\mc{L}}(zq^{-\kappa_{(1)}}).
\end{align}

The proof of the claim follows from Lemma \ref{vanishing-theorem:lem}. For $s\in\mbb{R}^n$ with all components non-negative. We have that:
\begin{align*}
\Delta_s(\mbf{B}_{\mc{L}}^s(z)')(u)=\Delta_s(\mbf{B}_{\mc{L}}^s(z))(u)\Delta_s(D^s_{\mc{L}}(z))
\end{align*}
When restricted to the component $\text{End}(K_{T\times\mbb{C}^{\times}_{q}}(M(0,\mbf{w}_1))_{loc}\otimes K_{T\times\mbb{C}^{\times}_{q}}(M(\mbf{v}_2+\mbf{v}_1,\mbf{w}_2))_{loc})$, we will have:
\begin{equation*}
\begin{aligned}
\lim_{u\rightarrow0}\Delta_s(\mbf{B}_{\mc{L}}^s(z)')(u)_{00}=&\lim_{u\rightarrow0}\Delta_s(\mbf{B}_{\mc{L}}^s(z))(u)_{00}\Delta_s(D^s_{\mc{L}}(z))_{00}+\lim_{u\rightarrow0}\sum_{\mbf{v}\neq0}\Delta_s(\mbf{B}_{\mc{L}}^s(z))(u)_{\mbf{v}0}\Delta_s(D^s_{\mc{L}}(z))_{0\mbf{v}}\\
=&\lim_{u\rightarrow0}\Delta_s(\mbf{B}_{\mc{L}}^s(z))(u)_{00}\Delta_s(D^s_{\mc{L}}(z))_{00}
\end{aligned}
\end{equation*}
the second term on the right hand side of the first equality is zero by Lemma \ref{vanishing-theorem:lem}. Thus again using Lemma \ref{vanishing-theorem:lem}. One obtain the desired result.

Now since $\Delta_s(D)(z)$ only contains the diagonal part, we have
\begin{align*}
\Delta_s(D^s_{\mc{L}})(z)=\Delta_s^{op}(D^s_{\mc{L}})(z)
\end{align*}
and thus we have
\begin{align*}
1\otimes D^s_{\mc{L}}(zq^{-\kappa_{(1)}})=1\otimes D^s_{\mc{L}}(zq^{\kappa_{(2)}})
\end{align*}
and thus we conclude that $D^s_{\mc{L}}(z)$ is $q$-periodic in $z$. Because $D^s_{\mc{L}}(z)$ is rational in $z$, it must be constant. The definition of $D^s_{\mc{L}}(z)$ therefore implies that this constant operator is the identity.

For the general $s\in\mbb{R}^n$, choose an integral vector $\mbf{n}\in\mbb{Z}^n$ such that all components of $s+\mbf{n}\in\mbb{Z}^n$ are non-negative. Using \ref{translationformula} we have that:
\begin{align*}
\mbf{B}_{\mc{L}}^s(z)=\mc{L}_{\mbf{n}}^{-1}\mbf{B}_{\mc{L}}^{s+\mbf{n}}(zp^{\mbf{n}})\mc{L}_{\mbf{n}}
\end{align*}
and thus we have that:
\begin{align*}
D^s_{\mc{L}}(z)=\mbf{B}_{\mc{L}}^{s}(z)^{-1}\mbf{B}_{\mc{L}}^{s}(z)'=\mc{L}_{\mbf{n}}^{-1}\mbf{B}_{\mc{L}}^{s+\mbf{n}}(zp^{\mbf{n}})^{-1}\mbf{B}_{\mc{L}}^{s+\mbf{n}}(zp^{\mbf{n}})'\mc{L}_{\mbf{n}}=\mc{L}_{\mbf{n}}^{-1}D_{\mc{L}}^{s+\mbf{n}}(zp^{\mbf{n}})\mc{L}_{\mbf{n}}=\text{Id}.
\end{align*}
and the last equality follows from the fact that $D_{\mc{L}}^{s+\mbf{n}}$ is an identity operator.

\end{proof}

To prove the claim in the proof of the theorem, we denote by $\Delta_s(\mbf{B}^s_{\mc{L}}(z))(u)_{\mbf{v}_1,\mbf{v}_1'}$ the component of the operator $\Delta_s(\mbf{B}^s_{\mc{L}}(z))(u)$ sending $K_{T\times\mbb{C}_q^*}(M(\mbf{v}_1,\mbf{w}_1))\otimes K_{T\times\mbb{C}_q^*}(M(\mbf{v}_2,\mbf{w}_2))$ to $K_{T\times\mbb{C}_q^*}(M(\mbf{v}_1',\mbf{w}_1'))\otimes K_{T\times\mbb{C}_q^*}(M(\mbf{v}_2',\mbf{w}_2'))$. We show the following lemma:
\begin{lem}\label{vanishing-theorem:lem}
For $s\in\mbb{R}^n$ such that all components of $s$ are non-negative, we have that
\begin{equation*}
\begin{aligned}
&\lim_{u\rightarrow0}\Delta_s(\mbf{B}^s_{\mc{L}}(z))(u)_{\mbf{v}_1,\mbf{v}_1'}=\mbf{B}^s_{\mc{L}}(zq^{\kappa_{(2)}})\otimes \mbf{B}^s_{\mc{L}}(zq^{-\kappa_{(2)}}),\qquad\mbf{v}_1=\mbf{v}_1'=0\\
&\lim_{u\rightarrow0}\Delta_s(\mbf{B}_{\mc{L}}^{s}(z))(u)_{\mbf{v}_1,\mbf{v}_1'}=0,\qquad\mbf{v}_1\neq0
\end{aligned}
\end{equation*}
\end{lem}
\begin{proof}
The argument is similar to the proof of Proposition 15 in \cite{OS22}. By Proposition \ref{coproduct-for-difference}:
\begin{align*}
\Delta_s(\mbf{B}^s_{\mc{L}}(z))=\Delta_{\mbf{m}_0}(\mbf{B}_{\mbf{m}_0}(z))(R_{\mbf{m}_0}^{+})^{-1}\cdots\Delta_{\mbf{m}_{m-1}}(\mbf{B}_{\mbf{m}_{m-1}}(z))(R_{\mbf{m}_{m-1}}^{+})^{-1}\Delta_{\infty}(\mc{L})
\end{align*}
and for each item:
\begin{align*}
\Delta_{\mbf{m}}(\mbf{B}_{\mbf{m}}(z))(R_{\mbf{m}}^+)^{-1}=\mbf{J}_{\mbf{m}}^-(z)(\mbf{B}_{\mbf{m}}(zq^{\kappa_{(2)}})\otimes\mbf{B}_{\mbf{m}}(zq^{-\kappa_{(1)}}))\mbf{J}_{\mbf{m}}^+(z)(R_{\mbf{m}}^{+})^{-1}
\end{align*}
and its $u$-weight part $\mbf{v}_1$ to $\mbf{v}_1+\mbf{n}$ is given by $\langle\mbf{m},\mbf{n}\rangle$. Therefore the $u$-weight of $\Delta_s(\mbf{B}^{s}_{\mc{L}}(z))(u)_{\mbf{v}_1,\mbf{v}_1'}$ is given by:
\begin{align*}
\langle\mc{L},\mbf{v}_1\rangle+\langle\mbf{m}_0,\mbf{n}_0\rangle+\cdots+\langle\mbf{m}_{m-1},\mbf{n}_{m-1}\rangle
\end{align*}
such that
\begin{align*}
\mbf{n}_0+\cdots+\mbf{n}_{m-1}=\mbf{v}_1'-\mbf{v}_1
\end{align*}
and also:
\begin{equation*}
\begin{aligned}
&\langle\mc{L},\mbf{v}_1\rangle+\langle\mbf{m}_0,\mbf{n}_0\rangle+\cdots+\langle\mbf{m}_{m-1},\mbf{n}_{m-1}\rangle\\
=&\langle\mc{L}+\mc{L}_0-\mbf{m}_{m-1},\mbf{v}_1\rangle+\langle\mbf{m}_{m-1}-\mbf{m}_{m-2},\mbf{v}_1+\mbf{n}_{m-1}\rangle+\cdots+\langle\mbf{m}_0-\mc{L}_0,\mbf{v}_1+\mbf{n}_{m-1}+\cdots+\mbf{n}_1\rangle+\\
&+\langle\mc{L}_0,\mbf{v}_1'-\mbf{v}_1\rangle
\end{aligned}
\end{equation*}
and here $\mc{L}_0$ is the closest integral vector to the slope point $s$. Also note that we always require that $\mbf{v}_1+\mbf{n}_1+\cdots+\mbf{n}_k\geq0$ for $1\leq k\leq m-1$. Otherwise the operator $\Delta_s(\mbf{B}^{s}_{\mc{L}}(z))(u)_{\mbf{v}_1,\mbf{v}_1'}$ will vanish.

Let us assume that $s\in\mbb{R}^n$ is all non-negative: It means that $s_i\geq0$ for every component $1\leq i\leq n$. In the case when $\mbf{v}_1=\mbf{v}_1'=0$, the $u$-weight of $\Delta_s(\mbf{B}^s_{\mc{L}}(z))(u)_{\mbf{v}_1,\mbf{v}_1'}$ is reduced to the following:
\begin{align*}
\langle\mbf{m}_{m-1}-\mbf{m}_{m-2},\mbf{n}_{m-1}\rangle+\cdots+\langle\mbf{m}_0-\mc{L}_0,\mbf{n}_{m-1}+\cdots+\mbf{n}_1\rangle
\end{align*}
and now since $\mbf{n}_1+\cdots+\mbf{n}_k\geq0$ for $1\leq k\leq m-1$ and $\mbf{v}_1-\mbf{v}_1=0$, it also forces that $\mbf{n}_i=0$ for $i=0,\cdots,m-1$.

When $\mbf{v}_1\neq0$, we will have that:
\begin{align*}
\lim_{u\rightarrow0}\Delta_s(\mbf{B}_{\mc{L}}^{s}(z))(u)_{\mbf{v}_1,\mbf{v}_1'}=0,\qquad\mbf{v}_1\neq0
\end{align*}

\end{proof}

\begin{thm}\label{commutativity-theorem}
For arbitrary line bundles $\mc{L}$, $\mc{L}'\in\text{Pic}(M(\mbf{v},\mbf{w}))$ and a slope $s\in\mbb{R}^n$, the corresponding $q$-difference operators commute:
\begin{align*}
\mc{A}^{s}_{\mc{L}}\mc{A}^{s}_{\mc{L}'}=\mc{A}^s_{\mc{L}'}\mc{A}^{s}_{\mc{L}}.
\end{align*}
\end{thm}
\begin{proof}
This follows from Proposition \ref{translationformula}:
\[
\mc{A}^{s}_{\mc{L}}\mc{A}^{s}_{\mc{L}'}=\mc{A}^{s}_{\mc{L}+\mc{L}'}=\mc{A}^s_{\mc{L}'}\mc{A}^{s}_{\mc{L}}.
\]
\end{proof}

To conclude, let $M(\mbf{v},\mbf{w})$ be the Nakajima quiver variety of the affine type $A$, let $\mc{L}_i$, $i=1,\cdots,n$ be the generators of $\text{Pic}(M(\mbf{v},\mbf{w}))$. The main result of the paper is the following. 

\begin{thm}\label{holonomic-theorem}
For an affine type $A$ quiver variety $M(\mbf{v},\mbf{w})$ and a line bundle $\mc{L}\in\text{Pic}(M(\mbf{v},\mbf{w}))$, there is a holonomic system of operators $\mbf{B}_{\mc{L}}^s(z)\in\text{End}(K_{T}(M(\mbf{v},\mbf{w}))[[z^{\mc{L}}]]_{\mc{L}\in\text{Pic}(M(\mbf{v},\mbf{w}))})$ defined in \ref{defnofqdeoperator}; explicitly,
\begin{align*}
T_{p,\mc{L}}^{-1}\mbf{B}^s_{\mc{L}}(z)T_{p,\mc{L}'}^{-1}\mbf{B}^s_{\mc{L}'}(z)=T_{p,\mc{L}'}^{-1}\mbf{B}^s_{\mc{L}'}(z)T_{p,\mc{L}}^{-1}\mbf{B}^s_{\mc{L}}(z).
\end{align*}
\end{thm}
\begin{proof}
This is a direct corollary of Theorem \ref{commutativity-theorem}.
\end{proof}

\subsection{Connection to the Okounkov-Smirnov quantum difference equation}
In \cite{OS22}, Okounkov and Smirnov constructed the quantum difference operator $\mbf{M}_{\mc{L}}(z)$ for Nakajima quiver varieties. It is shown in \cite{OS22} that the quantum difference operator may be written as
\begin{align*}
\mbf{M}_{\mc{L}}(z)=\text{Const}_s\cdot\mc{L}\cdot\prod_{w\in(s,s-\mc{L}]}\mbf{B}_{w}^{MO}(z)
\end{align*}
with $s\in\mbb{R}^n$ infinitesimally close to $0\in\mbb{R}^n$. Here $\mbf{B}_{w}^{MO}(z)$ are the monodromy operators defined as in \eqref{defn-of-quantum-difference-operator} on the wall subalgebra $U_q^{MO}(\mf{g}_w)$ defined in \cite{OS22}. 

In \cite{Z25}, we have proved that the quantum toroidal algebra $U_{q,t}(\hat{\hat{\mf{sl}}}_n)$ is isomorphic to the geometric quantum loop group defined in \cite{OS22}. Moreover, by Theorem 6.9 in \cite{Z25}, we have the following isomorphism of $\mbb{Q}(q,t)$-quasitriangular Hopf algebra between the wall subalgebra $U_{q}^{MO}(\mf{g}_w)$ and the slope subalgebra $\mc{B}_{\mbf{m}}$ of the generic slope $\mbf{m}\in\mbb{R}^n$ corresponding to the wall $w$:
\begin{equation*}
(\mc{B}_{\mbf{m}},R_{\mbf{m}}^+,\Delta_{\mbf{m}},S_{\mbf{m}},\epsilon,\eta)\cong(U_{q}^{MO}(\mf{g}_w),R_{\mbf{m}}^{+,MO},\Delta_{\mbf{m}}^{MO},S_{\mbf{m}}^{MO},\epsilon,\eta).
\end{equation*}

By the identification of the universal $R$-matrix $R_{\mbf{m}}^+=R_{\mbf{m}}^{+,MO}$, and since the upper-triangular solution with block-diagonal part as the identity operator for the ABRR equation defined in \ref{ABRR-equation} is unique, the identification of the universal $R$-matrix will give the same monodromy operator:
\begin{align*}
\mbf{B}_{w}^{MO}(z)=\mbf{B}_{\mbf{m}}(z).
\end{align*}
for the generic $\mbf{m}\in\mbb{R}^n$.

Hence:
\begin{equation*}
\mbf{M}_{\mc{L}}(z)=\text{Const}_s\cdot\mbf{B}_{\mc{L}}^s(z).
\end{equation*}

Thus we have the following result:
\begin{thm}\label{geometric-qde=algebraic-qde:label}
The algebraic quantum difference operator for the affine type $A$ Nakajima quiver variety of the slope $s$ infinitesimally close to $0\in\mbb{R}^n$ coincides with the Okounkov-Smirnov geometric quantum difference operator up to a constant operator.
\end{thm}

Therefore we can see that we have established a $\mbb{R}^n$-family of $q$-difference connections $\{\mc{A}^{s}_{\mc{L}}\}_{\mc{L}\in\text{Pic}(M(\mbf{v},\mbf{w}))}$ over $(\mbb{C}^*)^n$. When $s\in\mbb{R}^n$ is infinitesimally close to $0\in\mbb{R}^n$, the $q$-difference connection is the Okounkov-Smirnov geometric quantum difference equations. Combining the above formula, one can write down an algebraic formulation of the Okounkov-Smirnov geometric quantum difference equation of the Nakajima quiver variety of the affine type $A$ using the quantum toroidal algebra $U_{q,t}(\hat{\hat{\mf{sl}}}_n)$.

\section{\textbf{Examples}}
\subsection{Cotangent bundle of the Grassmannians}
We first realise finite type $A$ quiver varieties as special cases of affine type $A$ quiver varieties by setting the vector spaces at selected vertices equal to zero. For example,
\begin{equation*}
\begin{tikzcd}
	{V_1} & {V_2} & {V_3} & {V_4} \\
	{W_1} & {W_2} & {W_3} & {W_4}
	\arrow[curve={height=-6pt}, from=1-1, to=1-2]
	\arrow[curve={height=-6pt}, from=1-2, to=1-1]
	\arrow[curve={height=-6pt}, from=1-2, to=1-3]
	\arrow[curve={height=-6pt}, from=1-3, to=1-2]
	\arrow[curve={height=-6pt}, from=1-3, to=1-4]
	\arrow[curve={height=-6pt}, from=1-4, to=1-3]
	\arrow[curve={height=6pt}, from=1-1, to=2-1]
	\arrow[curve={height=6pt}, from=2-1, to=1-1]
	\arrow[curve={height=6pt}, from=1-2, to=2-2]
	\arrow[curve={height=6pt}, from=2-2, to=1-2]
	\arrow[curve={height=6pt}, from=1-3, to=2-3]
	\arrow[curve={height=6pt}, from=2-3, to=1-3]
	\arrow[curve={height=6pt}, from=1-4, to=2-4]
	\arrow[curve={height=6pt}, from=2-4, to=1-4]
\end{tikzcd}\hookrightarrow
\begin{tikzcd}
	&& 0 \\
	&& 0 \\
	{W_1} & {V_1} && 0 & 0 \\
	{W_2} & {V_2} && {V_4} & {W_4} \\
	&& {V_3} \\
	&& {W_3}
	\arrow[curve={height=-6pt}, from=3-2, to=3-1]
	\arrow[curve={height=-6pt}, from=3-1, to=3-2]
	\arrow[curve={height=-6pt}, from=4-2, to=4-1]
	\arrow[curve={height=-6pt}, from=4-1, to=4-2]
	\arrow[curve={height=-6pt}, from=5-3, to=6-3]
	\arrow[curve={height=-6pt}, from=6-3, to=5-3]
	\arrow[curve={height=-6pt}, from=4-4, to=4-5]
	\arrow[curve={height=-6pt}, from=4-5, to=4-4]
	\arrow[curve={height=-6pt}, from=3-4, to=3-5]
	\arrow[curve={height=-6pt}, from=3-5, to=3-4]
	\arrow[curve={height=-6pt}, from=1-3, to=2-3]
	\arrow[curve={height=-6pt}, from=2-3, to=1-3]
	\arrow[curve={height=6pt}, from=3-2, to=2-3]
	\arrow[curve={height=6pt}, from=2-3, to=3-2]
	\arrow[curve={height=6pt}, from=2-3, to=3-4]
	\arrow[curve={height=6pt}, from=3-4, to=2-3]
	\arrow[curve={height=6pt}, from=3-4, to=4-4]
	\arrow[curve={height=6pt}, from=4-4, to=3-4]
	\arrow[curve={height=6pt}, from=3-2, to=4-2]
	\arrow[curve={height=6pt}, from=4-2, to=3-2]
	\arrow[curve={height=6pt}, from=4-2, to=5-3]
	\arrow[curve={height=6pt}, from=5-3, to=4-2]
	\arrow[curve={height=6pt}, from=5-3, to=4-4]
	\arrow[curve={height=6pt}, from=4-4, to=5-3]
\end{tikzcd}
\end{equation*}

In this case, the monodromy operator reduces to
\begin{equation*}
\begin{aligned}
&\mbf{B}_{\mbf{m}}(z)=\prod_{k=0}^{\substack{\rightarrow\\\infty}}\prod_{\gamma\in\Delta(A)}^{\rightarrow}\exp_{q^{2}}(-(q-q^{-1})z^{-k(\mbf{v}_{\gamma}+m\bm{\delta}_h)}p^{-k\mbf{m}\cdot(\mbf{v}_{\gamma})}q^{-\frac{k}{2}(\mbf{v}^TC\mbf{v}-\mbf{w}^T\mbf{w})\mbf{\theta}\cdot(\mbf{v}_{\gamma}))}\\
&q^{2k(\mbf{v}_{\gamma})^TC(\mbf{v}_{\gamma})}(\varphi^{(k+1)(\mbf{v}_{\gamma})}f_{\gamma}\otimes e_{\gamma}'\varphi^{-(k+1)(\mbf{v}_{\gamma})}))).
\end{aligned}
\end{equation*}

We next consider cotangent bundles of Grassmannians.
\begin{equation*}
T^*Gr(k,n)=
\begin{tikzcd}
	{\mathbb{C}^k} & {\mathbb{C}^n}
	\arrow[curve={height=-6pt}, from=1-1, to=1-2]
	\arrow[curve={height=-6pt}, from=1-2, to=1-1]
\end{tikzcd}
\end{equation*}

A direct calculation shows that for generic $\mbf{m}\notin\mbb{Z}$, the monodromy operator $\mbf{B}_{\mbf{m}}(z)$ acts as zero on $K_{T}(T^*Gr(k,n))$. The quantum difference operator is therefore
\begin{align*}
\mbf{M}_{\mc{O}(1)}(z)=\mc{O}(1)\mbf{B}_{-e_1}(z)
\end{align*}
Moreover, $\mbf{B}_{-e_1}(z)$ is
\begin{equation*}
\begin{aligned}
\mbf{B}_{-e_1}(z)=&:\prod_{k=0}^{\infty}\exp_{q^2}(-(q-q^{-1})z^{-k}p^{-k}q^{-2k}f\otimes e):\\
=&\sum_{n=0}^{\infty}\frac{(q-q^{-1})^n}{[n]_{q^2}!}\frac{(-1)^n}{\prod_{k=1}^{n}(1-z^{-1}p^{-1}q^{2k})}f^ne^n.
\end{aligned}
\end{equation*}

This agrees with the formula in \cite{ABRR97,EV02}\cite{OS22}.

\subsection{Instanton moduli space}
This example appears in \cite{OS22}; we review it here for completeness.

The instanton moduli space is constructed from the Jordan quiver $\text{Rep}(n,r)=\text{Hom}(\mbb{C}^r,\mbb{C}^n)\oplus\text{Hom}(\mbb{C}^n,\mbb{C}^n)$, whose diagram is
\begin{equation*}
\begin{tikzcd}
	{\mathbb{C}^n} && {\mathbb{C}^r}
	\arrow[from=1-1, to=1-1, loop, in=145, out=215, distance=10mm]
	\arrow[from=1-1, to=1-1, loop, in=140, out=220, distance=15mm]
	\arrow[curve={height=6pt}, from=1-1, to=1-3]
	\arrow[curve={height=6pt}, from=1-3, to=1-1]
\end{tikzcd}
\end{equation*}

When $r=1$, the corresponding variety is isomorphic to the Hilbert scheme $\text{Hilb}_{m}(\mbb{C}^2)$ of points over $\mbb{C}^2$. As a vector space, the equivariant $K$-theory of Hilbert schemes can be identified with the ring of symmetric functions in infinitely many variables.
\begin{align*}
\oplus_{m=0}^{\infty}K_{T}(\text{Hilb}_{m}(\mbb{C}^2))_{loc}=F(u_1):=\mbb{Q}(q_1,q_2)[x_1,x_2,\cdots]^{Sym}.
\end{align*}

The fixed-point basis corresponds to the Macdonald polynomials $P_z$, which form an orthogonal basis of the space of symmetric functions over $\mbb{Q}(q,t)$.

Similarly, for the instanton moduli space $M(n,r)$ of rank $r$ and instanton number $n$, choose the whole framed torus $T=(\mbb{C}^{*})^{r}$ such that $\mbf{w}=u_1+\cdots+u_{r}$, and we have:
\begin{align*}
\oplus_{n=0}^{\infty}K_{T\times\mbb{C}^*_{q_1}\times\mbb{C}^*_{q_2}}(M(n,r))_{loc}=F(u_1)\otimes\cdots\otimes F(u_r).
\end{align*}

For convenience, we first recall the elliptic Hall algebra description of the quantum toroidal $U_{q,t}(\hat{\hat{\mf{gl}}}_1)$. 

\subsubsection{Quantum toroidal $U_{q,t}(\hat{\hat{\mf{gl}}}_1)$}
Set
\begin{align*}
n_{k}=\frac{(q_1^{k/2}-q_{1}^{-k/2})(q_2^{k/2}-q_{2}^{-k/2})(q^{-k/2}-q^{k/2})}{k},\qquad q=q_1q_2.
\end{align*}

For simplicity, we can use the elliptic Hall algebra description for the quantum toroidal algebra $U_{q}(\hat{\hat{\mf{gl}}}_{1})$:
\begin{align*}
U_{q}(\hat{\hat{\mf{gl}}}_{1})=\mbb{Q}(q_1,q_2)\langle K_{\mbf{a}},e_{\mbf{a}}\rangle/\{\text{relation in 1.1 of \cite{SV13} or 7.1 of \cite{OS22}}\}.
\end{align*}

Let $\mbf{a}=(a_1,a_2)\in\mbb{Z}^2\backslash\{(0,0)\}$. For $w\in\mbb{Q}\cup\{\infty\}$ we denote by $d(w)$ and $n(w)$ the denominator and numerator of $w$. We set $d(\infty)=0$ and $n(\infty)=1$. It follows directly that 
\begin{align*}
\alpha^{w}_{k}=e_{(d(w)k,n(w)k)},\qquad k\in\mbb{Z}\backslash\{0\}.
\end{align*}

generate a Heisenberg subalgebra of $\mc{B}_{w}\subset U_{q}(\hat{\hat{\mf{gl}}}_{1})$ with the following relations:
\begin{align*}
[\alpha_{-k}^{w},\alpha_{k}^{w}]=\frac{{q}^{-kd(w)}-q^{kd(w)}}{n_k}.
\end{align*}

Consequently, we obtain the slope factorisation
\begin{align*}
U_{q,t}(\hat{\hat{\mf{gl}}}_1)\cong\bigotimes^{\rightarrow}_{w\in\mbb{Q}\sqcup\{\infty\}}\mc{B}_{w}.
\end{align*}
where each $\mc{B}_{w}$ is isomorphic to $U_{q}(\hat{\mf{gl}}_1)$ for $w\in\mbb{Q}$. The algebra $U_{q_1,q_2}(\hat{\hat{\mf{gl}}}_{1})$ acts on $\oplus_{m}K_{T}(\text{Hilb}_{m}(\mbb{C}^2))_{loc}$\cite{SV13}. The central elements act in this representation by:
\begin{align*}
K_{(1,0)}=q_{1}^{-1/2}q_{2}^{-1/2},K_{(0,1)}=1.
\end{align*}

The vertical generators commute:
\begin{align*}
[e_{(0,m)},e_{(0,n)}]=0.
\end{align*}

The algebra $\mc{B}_{\mbf{0}}$ corresponds to the horizontal Heisenberg subalgebra
\begin{align*}
[e_{(m,0)},e_{(n,0)}]=-\frac{m}{(q_{1}^{m/2}-q_{1}^{-m/2})(q_2^{m/2}-q_2^{-m/2})}\delta_{m+n,0}.
\end{align*}

The vertical subalgebra $\{e_{(m,0)}\}$ is commutative, and it acts diagonally in the fixed-point basis, i.e., the basis of Macdonald polynomials:
\begin{align*}
e_{(0,l)}(P_{\lambda})=u_{1}^{-l}\text{sign}(l)(\frac{1}{1-q_1^{-l}}\sum_{i=1}^{\infty}q_{1}^{-lz_i}q_{2}^{-l(i-1)})P_{\lambda}.
\end{align*}

Hence the reduced part of the universal $R$-matrix of $U_{q}(\hat{\hat{\mf{gl}}}_{1})$ can be written as
\begin{align*}
R=\prod_{a/b\in\mbb{Q}\cup\{\infty\}}^{\rightarrow}\exp(-\sum_{k=1}^{\infty}n_{k}\alpha_{-k}^{a/b}\otimes\alpha_{k}^{a/b}).
\end{align*}

We now derive the quantum difference operator $\mbf{M}_{\mc{L}}(z)$.

We first derive the solution of the ABRR equation. We choose a subtorus $A$ such that it splits the framing by $r=r_1u_1+r_2u_2$ so that
\begin{align*}
K_{G}(M(r)^A)_{loc}=F^{\otimes r_1}(u_1)\otimes F^{\otimes r_2}(u_2).
\end{align*}

Let $F=M(m_1,r_1)\times M(m_2,r_2)$ be a component of $M(m,r)^A$. As $\text{dim}(M(m,r))=2mr$, we obtain that the corresponding eigenvalue of $\Omega$ equals:
\begin{align*}
\Omega=\frac{\text{codim}(F)}{4}=\frac{m_1r_1+m_2r_2}{2}
\end{align*}
the ABRR equation for a wall $w\in\mbb{Q}$ takes the form:
\begin{align*}
q^{\Omega}R_{w}^{-}z_{(1)}^{-1}J_{w}^{-}(z)=J_{w}^{-}(z)q^{\Omega}z_{(1)}^{-1}.
\end{align*}

For the derivation, see Section $7$ of \cite{OS22}. The result is
\begin{align*}
J_{w}^{-}(z)=\exp(-\sum_{k=1}^{\infty}\frac{n_{k}K_{(1,0)}^{kd(w)}\otimes K_{(1,0)}^{-kd(w)}}{1-z^{-kd(w)}K_{(1,0)}^{kd(w)}\otimes K_{(1,0)}^{-kd(w)}}\alpha_{-k}^{w}\otimes\alpha_{k}^w).
\end{align*}

This gives the monodromy operators
\begin{align*}
\mbf{B}_{w}(z)=\exp(-\sum_{k=1}^{\infty}\frac{n_kK_{(1,0)}^{kd(w)}}{1-z^{-kd(w)}p^{kn(w)}K_{(1,0)}^{kd(w)}}\alpha_{-k}^w\alpha_{k}^w).
\end{align*}

Let $\mc{L}=\mc{O}(1)$ be the generator of the Picard group. Then for every wall $w\in(s,s-\mc{O}(1))$ satisfies $-1\leq w<0$,
the definition of the wall set gives that the wall set is:
\begin{align*}
\text{Walls}(M(n,r)):=\{\frac{a}{b}|a\in\mbb{Z},|b|\leq n\}\cap[-1,0).
\end{align*}

This coincides with the one given by the computation of the $K$-theoretic stable envelope in \cite{S20}.

The quantum difference operator is therefore
\begin{align}\label{instanton-moduli-difference-operator}
\mbf{M}_{\mc{O}(1)}(z)=\mc{O}(1)\prod^{\leftarrow}_{\substack{w\in\mbb{Q}\\-1\leq w<0}}:\exp(-\sum_{k=1}^{\infty}\frac{n_kq^{-\frac{krd(w)}{2}}}{1-z^{-kd(w)}p^{kn(w)}q^{-krd(w)}}\alpha_{-k}^w\alpha_{k}^w):.
\end{align}

\textbf{Example.} For $M(2,5)$, the quantum difference operator is
\begin{equation*}
\begin{aligned}
\mbf{M}_{\mc{O}(1)}(z)=&\mc{O}(1)\mbf{B}_{-1}(z)\mbf{B}_{-4/5}(z)\mbf{B}_{-3/4}(z)\mbf{B}_{-2/3}(z)\mbf{B}_{-3/5}(z)\mbf{B}_{-1/2}(z)\times\\
&\times\mbf{B}_{-2/5}(z)\mbf{B}_{-1/3}(z)\mbf{B}_{-1/4}(z)\mbf{B}_{-1/5}(z)
\end{aligned}
\end{equation*}

The first few terms are
\begin{align*}
&\mbf{B}_{-1/5}(z)=1-\frac{n_1q^{-5}}{1-z^{-5}pq^{-5}}\alpha^{1/5}_{-1}\alpha^{1/5}_{1}\\
&\mbf{B}_{-1/4}(z)=1-\frac{n_1q^{-4}}{1-z^{-4}pq^{-4}}\alpha^{1/4}_{-1}\alpha^{1/4}_{1}\\
&\mbf{B}_{-1/3}(z)=1-\frac{n_1q^{-3}}{1-z^{-3}pq^{-3}}\alpha^{1/3}_{-1}\alpha^{1/3}_{1}\\
&\mbf{B}_{-2/5}(z)=1-\frac{n_1q^{-5}}{1-z^{-5}p^2q^{-5}}\alpha^{1/5}_{-1}\alpha^{1/5}_{1}\\
&\mbf{B}_{-1/2}(z)=1-\frac{n_1q^{-2}}{1-z^{-2}pq^{-2}}\alpha^{1/2}_{-1}\alpha^{1/2}_{1}-\frac{n_2q^{-4}}{1-z^{4}p^2q^{-4}}\alpha^{1/2}_{-2}\alpha^{1/2}_{2}+\frac{n_1^2q^{-4}}{2(1-z^{-2}pq^{-2})^2}\alpha^{1/2}_{-1}\alpha^{1/2}_{-1}\alpha^{1/2}_{1}\alpha^{1/2}_{1}
\end{align*}
\begin{equation*}
\begin{aligned}
&\mbf{B}_{-1}(z)=:\exp(-\frac{n_1q^{-1}}{1-z^{-1}pq^{-1}}\alpha^{-1}_{-1}\alpha^{-1}_{1}-\frac{n_2q^{-2}}{1-z^{-2}p^2q^{-2}}\alpha^{-1}_{-2}\alpha^{-1}_{2}-\frac{n_3q^{-3}}{1-z^{-3}p^3q^{-3}}\alpha^{-1}_{-3}\alpha^{-1}_{3}\\
&-\frac{n_4q^{-4}}{1-z^{-4}p^4q^{-4}}\alpha^{-1}_{-4}\alpha^{-1}_{4}-\frac{n_5q^{-5}}{1-z^{-5}p^5q^{-5}}\alpha^{-1}_{-5}\alpha^{-1}_{5}):.
\end{aligned}
\end{equation*}

The formula for the operator $\alpha_{-k}^{w}\alpha^{w}_{k}$ acting on the Fock space can be written as follows:
\begin{equation*}
\begin{aligned}
X^{-}X^{+}|z\rangle=&\sum_{\substack{z\supset\mu\subset\nu\\|z|=|\nu|}}R(\nu\backslash\mu)R(z\backslash\mu)q^{\frac{1}{2}|z\backslash\mu|}\prod_{\blacksquare\in z\backslash\mu}\frac{\prod_{\square\text{ i.c. of }\mu}(1-\frac{\chi_{\square}}{q\chi_{\blacksquare}})}{\prod_{\square\text{ o.c. of }\mu}(1-\frac{\chi_{\square}}{q\chi_{\blacksquare}})}\\
&\times\prod_{\blacksquare\in\nu\backslash\mu}\frac{\prod_{\square\text{ o.c. of }\nu}(1-\frac{\chi_{\square}}{\chi_{\blacksquare}})}{\prod_{\square\in\text{ i.c. of }\nu}(1-\frac{\chi_{\square}}{\chi_{\blacksquare}})}|\nu\rangle
\end{aligned}
\end{equation*}
with $X^{+}\in\mc{S}^{+}$ and $X^{-}\in\mc{S}^{-}$. And in our case, the operator $\alpha^{a/b}_{k}$ can be written in terms of the shuffle algebra element as:
\begin{align*}
\alpha^{a/b}_{k}=\frac{(q_1-1)^{bk}(1-q_2)^{bk}}{(q_1^{k}-1)(1-q_2^{k})}\text{Sym}[\frac{\prod_{i=1}^{ak}z_{i}^{\lfloor\frac{ib}{a}\rfloor-\lfloor\frac{(i-1)b}{a}\rfloor}}{\prod_{i=1}^{ak-1}(1-\frac{qz_{i+1}}{z_{i}})}\sum_{s=0}^{k-1}q^s\frac{z_{a(k-1)+1}\cdots z_{a(k-s)+1}}{z_{a(k-1)}\cdots z_{a(k-s)}}\prod_{i<j}\zeta(\frac{z_i}{z_j})].
\end{align*}

\subsection{Hilbert scheme $\text{Hilb}_{n}([\mbb{C}^2/\mbb{Z}_{r+1}])$}\label{section4-3}

Consider the equivariant Hilbert scheme $\text{Hilb}_{n}([\mbb{C}^2/\mbb{Z}_{r+1}])$. This is the Nakajima quiver variety of affine type $\hat{A}_{r}$. Its quiver is written as:
\begin{equation*}
\begin{tikzcd}
	& \cdots \\
	{\mathbb{C}^n} && {\mathbb{C}^n} \\
	{\mathbb{C}^n} && {\mathbb{C}^n} \\
	& {\mathbb{C}^n} \\
	& {\mathbb{C}}
	\arrow[curve={height=-6pt}, from=2-1, to=1-2]
	\arrow[curve={height=-6pt}, from=1-2, to=2-1]
	\arrow[curve={height=-6pt}, from=2-1, to=3-1]
	\arrow[curve={height=-6pt}, from=3-1, to=2-1]
	\arrow[curve={height=-6pt}, from=4-2, to=3-1]
	\arrow[curve={height=-6pt}, from=3-1, to=4-2]
	\arrow[curve={height=-6pt}, from=1-2, to=2-3]
	\arrow[curve={height=-6pt}, from=2-3, to=1-2]
	\arrow[curve={height=-6pt}, from=2-3, to=3-3]
	\arrow[curve={height=-6pt}, from=3-3, to=2-3]
	\arrow[curve={height=-6pt}, from=3-3, to=4-2]
	\arrow[curve={height=-6pt}, from=4-2, to=3-3]
	\arrow[curve={height=6pt}, from=4-2, to=5-2]
	\arrow[curve={height=6pt}, from=5-2, to=4-2]
\end{tikzcd}
\end{equation*}
with $r+1$ vertices and $1$ framed vertex. Recall that for the slope subalgebra $\mc{B}_{\mbf{m}}$, we have an isomorphism:
\begin{align*}
\mc{B}_{\mbf{m}}\cong\bigotimes_{h=1}^{g}U_{q}(\hat{\mf{gl}}_{l_h}).
\end{align*}

It has the slope decomposition as:
\begin{align*}
U_{q,t}(\hat{\hat{\mf{sl}}}_{r+1})\cong\bigotimes_{r\in\mbb{Q}}\mc{B}_{\mbf{m}+r\bm{\theta}}.
\end{align*}

Here $\mbf{m}=(m_1,\cdots,m_{r+1})\in\mbb{Q}^{r+1}$ and $\bm{\theta}=(1,\cdots,1)$, and for $\mbf{m}\in\mbb{Q}^{r+1}$, we have the following result.

We simply explain how to construct the positive integers $g,l_1,\cdots,l_g$. Set
\begin{align*}
e_{i}=\underbrace{(0,\cdots,1}_{i\text{ numbers}},\cdots,0).
\end{align*}

The set of vectors $\{e_{i}\}_{1\leq i\leq r+1}$ spans the $\mathbb{Q}$-vector space $\mbb{Q}^{r+1}$.

We define:
\begin{align}\label{cyclic-vectors}
[i;j):=e_i+e_{i+1}+\cdots+e_{j-1}
\end{align}
with $e_{i}:=e_{i\text{ mod }(r+1)}$. We call this vector an \textbf{arc vector}. We say that an arc vector $[i;j)$ is $\mbf{m}$-\textbf{integral} if
\begin{align*}
\mbf{m}\cdot[i;j)\in\mbb{Z}.
\end{align*}

It is straightforward to show that an $\mbf{m}$-integral arc vector exists as above, starting at each vertex $i$. Therefore there exists a well-defined minimal $\mbf{m}$-integral arc vector starting at each vertex $i$, and we will denote it by
\begin{align*}
[i;v_{\mbf{m}}(i)).
\end{align*}

Because the quiver has only $n$ vertices in the quiver with modulo relations, the uniqueness of the minimal $\mbf{m}$-integral arc vector implies that the resulting graph will be a union of oriented cycles:
\begin{align}\label{oriented-cycles}
C_1\sqcup\cdots\sqcup C_g
\end{align}
where $C_{h}=\{i_1,\cdots,i_{h}\}$ such that $i_{e+1}=v_{\mbf{m}}(i_{e})$ and $C_{h}=\{i_1,\cdots,i_{l_h}\}$.

To determine $\mbf{B}_{\mbf{m}}(z)$ in the present case of $\text{Hilb}_{n}([\mbb{C}^2/\mbb{Z}_{r}])$, $\mbf{v}_2=\mbf{w}_2=0$, $\mbf{v}_1=n\bm{\theta}$, $\mbf{w}_1=\mbf{e}_1$, and thus we have $\mbf{v}_2-\mbf{v}_1=-m\bm{\theta}$, $\mbf{w}_2-\mbf{w}_1=-\mbf{e}_1$.
Applying Theorem \ref{universal-formula-for-difference-operator} gives the following result.
\begin{prop}\label{explicitformulaqdoperator}
For the monodromy operator $\mbf{B}_{\mbf{m}}(z)\in\mc{B}_{\mbf{m}}$, where $\mbf{B}_{\mbf{m}}(z)$ is defined in \eqref{defnofqdeoperator} and $\mc{B}_{\mbf{m}}$ is defined as \eqref{rootquantum}, its representation in $\text{End}(K_{T}(\text{Hilb}_{n}([\mbb{C}^2/\mbb{Z}_{r}]))$ is given by:
\begin{equation}\label{qdehilbar}
\begin{aligned}
&\mbf{B}_{\mbf{m}}(z)\\
=&\mbf{m}(\prod_{h=1}^{g}(\exp(-\sum_{k=1}^{\infty}\frac{n_kq^{-\frac{k\lvert\bm{\delta_h}\lvert}{2}}}{1-z^{-k\lvert\bm{\delta_h}\lvert}p^{k\mbf{m}\cdot\bm{\delta}_{h}}q^{-\frac{k\lvert\bm{\delta_h}\lvert}{2}}}\alpha^{\mbf{m},h}_{-k}\otimes\alpha^{\mbf{m},h}_{k})\prod_{k=0}^{\substack{\rightarrow\\\infty}}\\
&\times\prod_{\substack{\gamma\in\Delta(A)\\m\geq0}}^{\leftarrow}(\exp_{q^{2}}(-(q-q^{-1})z^{-k(-\mbf{v}_{\gamma}+(m+1)\bm{\delta}_h)}p^{-k\mbf{m}\cdot(-\mbf{v}_{\gamma}+(m+1)\bm{\delta}_h)}q^{-k(-\mbf{v}_{\gamma}+(m+1)\bm{\delta}_h)^T((\frac{n^2r-1}{2})\bm{\theta}+\mbf{e}_1)-2-2\delta_{1\gamma}}\\
&f_{(\delta-\gamma)+m\delta}\otimes e_{(\delta-\gamma)+m\delta}')\exp(-(q-q^{-1})\sum_{m\in\mbb{Z}_{+}}\sum_{i,j=1}^{l_h}z^{-km\bm{\delta}_h}p^{-km\mbf{m}\cdot\bm{\delta}_{h}}q^{km\bm{\delta}_h^T((\frac{n^2r-1}{2})\bm{\theta}+\mbf{e}_1)-2}u_{m,ij}f_{m\delta,\alpha_i}\otimes e_{m\delta,\alpha_i}')\\
&\times\prod_{\substack{\gamma\in\Delta(A)\\m\geq0}}^{\rightarrow}\exp_{q^{2}}(-(q-q^{-1})z^{-k(\mbf{v}_{\gamma}+m\bm{\delta}_h)}p^{-k\mbf{m}\cdot(\mbf{v}_{\gamma}+m\bm{\delta}_h)}q^{k(\mbf{v}_{\gamma}+m\bm{\delta}_h)^T((\frac{n^2r-1}{2})\bm{\theta}+\mbf{e}_1)-2-2\delta_{1\gamma}}f_{\gamma+m\delta}\otimes e_{\gamma+m\delta}')))
\end{aligned}
\end{equation}
and $e_{v}'=S(e_{v})$, $f_{v}'=S(f_{v})$ are the images under the antipode map. $n_{k}$ is an element in $\mbb{Q}(q,t)$ such that $[\alpha^{\mbf{m},h}_{-k},\alpha^{\mbf{m},h}_{k}]=\frac{1}{n_{k}}$.
\end{prop}

To determine the quantum difference operator, we choose the slope from $0$ to $-\mc{L}$  with $\mc{L}=\mc{L}_1\cdots\mc{L}_{r+1}$. We choose a straight line of the equation $L_t:=(ts,\cdots,ts)$ to connect from $-\bm{\theta}$ to $0$. We see that the quantum difference equation over this line has the form:
\begin{align*}
\mbf{M}_{\mc{O}(1)}(z)=\mc{O}(1)\prod^{\rightarrow}_{w\in\text{Walls}\cap L_t}\mbf{B}_{w}(z).
\end{align*}

The finite set $\text{Walls}\cap L_t$ is described as follows.
\begin{lem}
The finite set $\text{Walls}\cap L_t$ has the following elements:
\begin{align*}
\text{Walls}_{0}:=\text{Walls}\cap L_t=\{(\frac{a}{b},\cdots,\frac{a}{b})|\lvert a\lvert<\lvert b\lvert, b\in[1,n(r+1)]\}.
\end{align*}
which coincides with the one defined by the $K$-theoretic stable envelope.
\end{lem}
\begin{proof}
Note that in the slope subalgebra $\mc{B}_{\mbf{m}}$ with the slope $\mbf{m}=(\frac{a}{b},\cdots,\frac{a}{b})$, the basic generator $e_{[i;i+1)}$ sends $K_{T}(M(n\bm{\theta},e_j))$ to $K_{T}(M(n\bm{\theta}-[i;i+b)-e_{i},e_i))$. For each vertex $v_i$, under the operation $e_{[i;i+1)}$, it can be sent to $v_{i}-[b/(r+1)]-1$. If we set $v_i=n$, to require that $M(n\bm{\theta}-[i;i+b)-e_{i},e_i)$ is non-empty, we set that:
\begin{align*}
n-[\frac{b-1}{r+1}]-1\geq0.
\end{align*}

Therefore, $b\leq n(r+1)$.

\end{proof}

We compute the quantum difference operator associated with $\mc{O}(1)=\mc{L}_1\cdots\mc{L}_{r+1}\in\text{Pic}(\text{Hilb}_{n}([\mbb{C}^2/\mbb{Z}_{r+1}])$, that is, we choose a path from $0$ to $-\mc{O}(1)$. 

Let us choose the straight line from $0$ to $\mc{O}(1)$. For $-\frac{a}{b}\bm{\theta}\in\text{Walls}\cap L_t$, we know that by the result \ref{rootquantum} the corresponding root subalgebra should be:
\begin{align*}
\mc{B}_{-\frac{a}{b}\bm{\theta}}=U_{q}(\hat{\mf{gl}}_{\frac{r+1}{g}})^{\otimes g},\qquad g=\text{gcd}(b,r+1).
\end{align*}

Thus the corresponding quantum difference operator is:
\begin{align*}
\mbf{M}_{\mc{O}(1)}(z)=\mc{O}(1)\prod^{\rightarrow}_{-\frac{a}{b}\bm{\theta}\in\text{Walls}\cap L_t,a\neq0}\mbf{B}_{-\frac{a}{b}\bm{\theta}}(z).
\end{align*}

As a further example, we choose the slope $\mc{L}_{i}\mc{O}(1)$. We also fix $\mc{L}_{i}\mc{O}(1)$, the interval $I_{i}:=(\underbrace{w,\cdots,2w}_{i\text{ numbers }},\cdots w)$ and $-1<w\leq 0$. A direct computation gives
\begin{align*}
\text{Walls}_{i,0}(\text{Hilb}_{n}([\mbb{C}^2/\mbb{Z}_{r+1}]))\cap I_{i}=\{(\underbrace{\frac{a}{b},\cdots,\frac{2a}{b}}_{i\text{ numbers }},\cdots\frac{a}{b})|a\in\mbb{Z},1\leq b\leq n(r+1)+1, -1<\frac{a}{b}\leq0\}.
\end{align*}

The corresponding quantum difference operator is
\begin{align*}
\mbf{B}^s_{\mc{L}_i\mc{O}(1)}(z)=\mc{L}_i\mc{O}(1)\prod^{\rightarrow}_{w\in\text{Walls}_{i,0}}\mbf{B}_{w}(z).
\end{align*}

For each $w\in\text{Walls}_{i,0}$, the corresponding root subalgebra $\mc{B}_w$ can be computed using the formula above. It is worth noting that the corresponding wall structure of $\text{Hilb}_{n}([\mbb{C}^2/\mbb{Z}_{r+1}])$ above can be computed from the elliptic stable envelope formula using the formula in \cite{D21}.

\textbf{Example: The $A_{r-1}$ surfaces $\widehat{\mbb{C}^2/\mbb{Z}_{r}}$.} The monodromy operator $\mbf{B}_{0}(z)$ is
\begin{equation*}
\begin{aligned}
\mbf{B}_{0}(z)=&1-\frac{n_1q^{-r/2}}{1-z^{-r}q^{-kr/2}}\alpha_{-1}^0\alpha_1^0-\sum_{\gamma\in\Delta(A)}\frac{(q-q^{-1})q^{-2-2\delta_{1\gamma}}}{1-z^{\mbf{v}_{\gamma}-\bm{\delta}_h}q^{-(-\mbf{v}_{\gamma}+\bm{\delta}_h)^T((\frac{r-1}{2})\bm{\theta}+\mbf{e}_1)}}f_{(\delta-\gamma)}e_{(\delta-\gamma)}'\\
&-(q-q^{-1})\sum_{i,j=1}^{r}\frac{1}{1-z^{-\bm{\delta}_h}q^{-\bm{\delta}_h^T((\frac{r-1}{2})\bm{\theta}+\mbf{e}_1)-2}}u_{1,ij}f_{\delta,\alpha_i}e'_{\delta,\alpha_i}\\
&-\sum_{\gamma\in\Delta(A)}\frac{(q-q^{-1})q^{-2-2\delta_{1\gamma}}}{1-z^{-\mbf{v}_{\gamma}}q^{-(\mbf{v}_{\gamma})^T((\frac{r-1}{2})\bm{\theta}+\mbf{e}_1)}}f_{\gamma}e_{\gamma}'.
\end{aligned}
\end{equation*}

In \cite{Z26}, we generalise this construction of quantum difference equations to quiver varieties of arbitrary type and to weight spaces of other classes of representations. It would also be interesting to use the construction to investigate further structures on the quantum difference equations.

\section{\textbf{Appendix: Proof of the Proposition \ref{rationality-of-R-matrix}}}\label{appendix}
In this appendix, we prove that $\prod_{\mu\in\mbb{R}_{>0}}R_{\mbf{m}+\mu\bm{\theta}}^+(u)$ has matrix coefficients that are rational functions over $u$ when acting on $K_{T}(M(\mbf{w}_1))\otimes K_{T}(M(\mbf{w}_2))$. Now we choose the fixed-point basis $|\bm{\mu}_i\rangle$, $|\bm{\lambda}_{i}\rangle$ for $i=1,2$. It is equivalent to proving that the matrix coefficients:
\begin{align*}
\langle\bm{\lambda}_1|\otimes\langle\bm{\lambda}_2|\prod_{\mu\in\mbb{Q}_{>0}}R_{\mbf{m}+\mu\bm{\theta}}^+(u)|\bm{\mu}_1\rangle\otimes|\bm{\mu}_2\rangle
\end{align*}
is a rational function in $u$.

Let $\mc{L}$ denote the corresponding line bundle for $\bm{\theta}\in(\mbb{Z}^+)^n$. In the proof of Proposition \ref{translationpic}, we have the following formula:
\begin{align*}
R_{\mbf{m}+\mc{L}}^+(u)=\mc{L}R_{\mbf{m}}^+(u)\mc{L}^{-1}
\end{align*}

Since $R_{\mbf{m}}^+(u)$ is strictly upper- triangular, we can write:
\begin{align*}
R_{\mbf{m}}^+(u)=\text{Id}+U_{\mbf{m}}^+(u)
\end{align*}
with $U_{\mbf{m}}^+(u)$ the upper-triangular part of $R_{\mbf{m}}^+(u)$ with zero diagonal.

We expand in the fixed-point basis as follows:
\begin{equation*}
\begin{aligned}
&\langle\bm{\lambda}_1|\otimes\langle\bm{\lambda}_2|\prod_{\mu\in\mbb{Q}_{>0}}R_{\mbf{m}+\mu\bm{\theta}}^+(u)|\bm{\mu}_1\rangle\otimes|\bm{\mu}_2\rangle=\langle\bm{\lambda}_1|\otimes\langle\bm{\lambda}_2|\prod_{k=1}^{\substack{\infty\\\rightarrow}}\prod_{\mu\in[0,1)}^{\rightarrow}R_{\mbf{m}+k\mu\bm{\theta}}^+(u)|\bm{\mu}_1\rangle\otimes|\bm{\mu}_2\rangle\\
=&\langle\bm{\lambda}_1|\otimes\langle\bm{\lambda}_2|\prod_{k=0}^{\substack{\infty\\\rightarrow}}\text{Ad}_{\mc{L}^k}\prod_{\mu\in[0,1)}^{\rightarrow}R_{\mbf{m}+\mu\bm{\theta}}^+(u)|\bm{\mu}_1\rangle\otimes|\bm{\mu}_2\rangle\\
=&\langle\bm{\lambda}_1|\otimes\langle\bm{\lambda}_2|\prod_{k=0}^{\substack{\infty\\\rightarrow}}\text{Ad}_{\mc{L}}^k\prod_{\mu\in[0,1)}^{\rightarrow}(\text{Id}+U_{\mbf{m}+\mu\bm{\theta}}^+(u))|\bm{\mu}_1\rangle\otimes|\bm{\mu}_2\rangle
\end{aligned}
\end{equation*}
By construction and Proposition \ref{finiteness-of-wall-set}, the infinite product $\prod_{\mu\in[0,1)}^{\rightarrow}(\text{Id}+U_{\mbf{m}+\mu\bm{\theta}}^+)$ reduces to a finite product on each fixed vector in $K_{T}(M(\mbf{v}_1,\mbf{w}_1))\otimes K_{T}(M(\mbf{v}_2,\mbf{w}_2))$ and $K_{T}(M(\mbf{v}_1',\mbf{w}_1))\otimes K_{T}(M(\mbf{v}_2',\mbf{w}_2))$. Without loss of generality, under the fixed-point basis, we assume that:
\begin{align*}
\prod_{\mu\in[0,1)}^{\rightarrow}(\text{Id}+U_{\mbf{m}+\mu\bm{\theta}}^+(u))=(\text{Id}+U_{\mbf{m}+\mu_1\bm{\theta}}^+(u))\cdots(\text{Id}+U_{\mbf{m}+\mu_n\bm{\theta}}(u))
\end{align*}
Thus in this case we have
\begin{equation*}
\begin{aligned}
&\langle\bm{\lambda}_1|\otimes\langle\bm{\lambda}_2|\prod_{k=0}^{\substack{\infty\\\rightarrow}}\text{Ad}_{\mc{L}}^k\prod_{\mu\in[0,1)}^{\rightarrow}(\text{Id}+U_{\mbf{m}+\mu\bm{\theta}}^+(u))|\bm{\mu}_1\rangle\otimes|\bm{\mu}_2\rangle\\
=&\langle\bm{\lambda}_1|\otimes\langle\bm{\lambda}_2|\prod_{k=0}^{\substack{\infty\\\rightarrow}}\text{Ad}_{\mc{L}}^k(\text{Id}+U_{\mbf{m}+\mu_1\bm{\theta}}^+(u))\cdots(\text{Id}+U_{\mbf{m}+\mu_n\bm{\theta}}(u))|\bm{\mu}_1\rangle\otimes|\bm{\mu}_2\rangle\\
=&\langle\bm{\lambda}_1|\otimes\langle\bm{\lambda}_2|\prod_{k=0}^{\substack{\infty\\\rightarrow}}\text{Ad}_{\mc{L}}^k[\text{Id}+\sum_{l=1}^{n}U_l^+(u)|\bm{\mu}_1]\rangle\otimes|\bm{\mu}_2\rangle,\qquad U_l(u):=\sum_{1\leq i_1<i_2<\cdots<i_l\leq n}U^{+}_{\mbf{m}+\mu_{i_1}\bm{\theta}}(u)\cdots U^{+}_{\mbf{m}+\mu_{i_l}\bm{\theta}}(u)\\
=&\langle\bm{\lambda}_1|\otimes\langle\bm{\lambda}_2|\text{Id}+\sum_{m=1}^{\infty}\sum_{k_1,\cdots,k_m=0}^{\infty}\sum_{1\leq l_1,\cdots,l_m\leq n}\text{Ad}_{\mc{L}}^{k_1}U_{l_1}(u)\cdots\text{Ad}_{\mc{L}}^{k_m}U_{l_m}(u)|\bm{\mu}_1\rangle\otimes|\bm{\mu}_2\rangle\\
=&\langle\bm{\lambda}_1|\otimes\langle\bm{\lambda}_2|\text{Id}|\bm{\mu}_1\rangle\otimes|\bm{\mu}_2\rangle+\sum_{m=1}^{\infty}\sum_{k_1,\cdots,k_m=0}^{\infty}\sum_{1\leq l_1,\cdots,l_m\leq n}\sum_{\bm{\lambda}_1,\cdots,\bm{\lambda}_{m-1}}\langle\bm{\lambda}|\text{Ad}_{\mc{L}}^{k_1}U_{l_1}(u)|\bm{\lambda}_1\rangle\langle\bm{\lambda}_1|\text{Ad}_{\mc{L}}^{k_2}U_{l_2}(u)|\bm{\lambda}_2\rangle\cdots\\
&\cdots\langle\bm{\lambda}_{m-1}|\text{Ad}_{\mc{L}}^{k_m}U_{l_m}(u)|\bm{\mu}\rangle\\
=&\delta_{\bm{\lambda},\bm{\mu}}+\sum_{m=1}^{\infty}\sum_{k_1,\cdots,k_m=0}^{\infty}\sum_{1\leq l_1,\cdots,l_m\leq n}\sum_{\substack{\square_i\in\bm{\lambda}_i,i=1}}^{m-1}\sum_{\substack{\square\in\bm{\lambda}\\\blacksquare\in\bm{\mu}}}(\frac{\chi_{\square}}{\chi_{\square_1}})^{k_1}\cdots(\frac{\chi_{\square_{m-1}}}{\chi_{\blacksquare}})^{k_m}\langle\bm{\lambda}|U_{l_1}(u)|\bm{\lambda}_1\rangle\langle\bm{\lambda}_1|U_{l_2}(u)|\bm{\lambda}_2\rangle\cdots\\
&\times\langle\bm{\lambda}_{m-1}|U_{l_m}(u)|\bm{\mu}\rangle
\end{aligned}
\end{equation*}
Combining the preceding formulas gives
\begin{equation*}
\begin{aligned}
&\langle\bm{\lambda}_1|\otimes\langle\bm{\lambda}_2|\prod_{k=0}^{\substack{\infty\\\rightarrow}}\text{Ad}_{\mc{L}}^k\prod_{\mu\in[0,1)}^{\rightarrow}(\text{Id}+U_{\mbf{m}+\mu\bm{\theta}}^+(u))|\bm{\mu}_1\rangle\otimes|\bm{\mu}_2\rangle\\
=&\delta_{\bm{\lambda},\bm{\mu}}+\sum_{m=1}^{\infty}\sum_{1\leq l_1,\cdots,l_m\leq n}\sum_{\substack{\square_i\in\bm{\lambda}_i,i=1}}^{m-1}\sum_{\substack{\square\in\bm{\lambda}\\\blacksquare\in\bm{\mu}}}\frac{1}{(1-\frac{\chi_{\square}}{\chi_{\square_1}})(1-\frac{\chi_{\square_1}}{\chi_{\square_2}})\cdots(1-\frac{\chi_{\square_{m-1}}}{\chi_{\blacksquare}})}\langle\bm{\lambda}|U_{l_1}(u)U_{l_2}(u)\cdots U_{l_m}(u)|\bm{\mu}\rangle
\end{aligned}
\end{equation*}

Since $U_{l_1}(u)\cdots U_{l_m}(u)$ is a Laurent polynomial in $u$, the required rationality follows.

\end{document}